\documentclass[11pt,leqno]{amsart}
\usepackage{amsmath,stmaryrd}
\usepackage{amssymb}
\usepackage{amsthm}
\usepackage{epsf}
\usepackage{graphicx}
\usepackage{esint} 
\usepackage{dsfont}
\usepackage[pagebackref]{hyperref} 
\hypersetup{pdfpagemode=FullScreen, colorlinks=true, citecolor={blue}} 
\usepackage{enumerate} 
\usepackage{xcolor,bbm,easybmat,bm}
\usepackage{tikz-cd}
\usepackage{centernot}
\usepackage{comment}
\usepackage{mathrsfs}
\usepackage{cancel}

\numberwithin{equation}{section}

\newcommand\footnoteref[1]{\protected@xdef\@thefnmark{\ref{#1}}\@footnotemark}
\makeatother

\DeclareMathOperator*{\osc}{osc}
\DeclareMathOperator{\Jac}{Jac}

\DeclareMathOperator{\diver}{div}

\DeclareMathOperator{\dist}{dist}
\DeclareMathOperator{\supp}{supp}
\DeclareMathOperator{\diam}{diam}

\newcommand{\1}{{\mathds 1}}
\newcommand{\ds}{\displaystyle}
\newcommand{\ms}{\medskip}

\newcommand{\R}{\mathbb R}
\newcommand{\bN}{\mathbb N}

\newcommand{\bp}{\noindent {\em Proof: }}
\newcommand{\ep}{\hfill $\square$ \medskip}

\newcommand{\A}{\mathcal A}
\newcommand{\N}{\mathcal N}

\newcommand{\norm}[1]{\left\Vert#1\right\Vert}

\newcommand{\br}[1]{\left(#1\right)}
\newcommand{\set}[1]{\left\{#1\right\}}

\newcommand{\om}{\Omega}
\newcommand{\pom}{\partial\Omega}

\def\YYint#1#2#3{{\setbox0=\hbox{$#1{#2#3}{\iint}$}
		\vcenter{\hbox{$#2#3$}}\kern-.51\wd0}}
 
\newcommand{\ffbox}[2]{\fbox{\begin{minipage}{3cm} \begin{center}  #1 \line(1, 0){\textwidth} \linebreak #2 \end{center} \end{minipage}}}

\newcommand{\sfbox}[2]{\fbox{\begin{minipage}{2.2cm} \begin{center}  #1 \line(1, 0){\textwidth} \linebreak #2 \end{center} \end{minipage}}}

\newcommand{\fffbox}[2]{\fcolorbox{black}{black!20}{\begin{minipage}{3.2cm}  \begin{center}  #1 \line(1, 0){\textwidth} \linebreak #2 \end{center} \end{minipage}}}

\newcommand{\sffbox}[2]{\fcolorbox{black}{black!20}{\begin{minipage}{3cm}  \begin{center}  #1 \line(1, 0){\textwidth} \linebreak #2 \end{center} \end{minipage}}}

\theoremstyle{plain}
\newtheorem{theorem}[equation]{Theorem}
\newtheorem{theoremc}[equation]{Theorem (Conjecture)}
\newtheorem{lemma}[equation]{Lemma}
\newtheorem{corollary}[equation]{Corollary}
\newtheorem{proposition}[equation]{Proposition}

\newtheorem{definition}[equation]{Definition}

\theoremstyle{definition}

\theoremstyle{remark}
\newtheorem{remark}[equation]{Remark}

\begin{document}
\title[Stability of $L^p$ Poisson-Neumann problem ]{Stability of the $L^p$ Poisson-Neumann problem under small perturbations of the coefficients and the domain}

\author[Feneuil]{Joseph Feneuil}
\address{Joseph Feneuil. Laboratoire de math\'ematiques d'Orsay, Universit\'e Paris-Saclay, France}
\email{joseph.feneuil@universite-paris-saclay.fr}

\author[Li]{Linhan Li} 
\address{Linhan Li. School of Mathematics, The University of Edinburgh, Edinburgh, UK }
\email{linhan.li@ed.ac.uk}

\maketitle

\begin{abstract} 
In this article, we obtain three stability results for the Poisson-Neumann problem: (1) the stability of the $L^p$ Poisson-Neumann problem under small Carleson perturbations of the coefficients on 1-sided Chord Arc domains, (2) the preservation of the $L^p$ Poisson-Neumann problem when transitioning from the Laplacian to operators in a special form with oscillations satisfying a small Carleson condition (``small special DKP operators''), and 
(3) the stability of the $L^p$ Poisson-Neumann problem under small bi-Lipschitz transformations of domains.

In contrast to the Neumann problem, the Poisson-Neumann problem is notably more self-contained, as results pertaining to it do not require assumptions regarding the 
$L^{p'}$ Dirichlet problem, making it a more suitable object for a perturbation theory.

As intermediate steps, we establish estimates on the square functional for solutions to the Poisson problem and their gradients on chord arc domains, which are of independent interest for specialists.
\end{abstract}

\ms\noindent{\bf Keywords: }  $L^p$ Poisson-Neumann problem, Laplacian, bi-Lipschitz transformation, Lipschitz domains, Carleson perturbations, 1-sided Chord Arc Domain, DKP operators, Chord Arc Domains, Green function, $A_\infty$-absolute continuity of the elliptic measure. 

\ms\noindent{\bf AI Disclaimer: } 
The present manuscript has been written without the use of AI for the development of mathematical concepts, ideas, organization, and the written content. The only exception is a limited application of AI for minor rephrasing at the beginning of the introduction (i.e. before Subsection 1.1), and for an additional proofreading when the full manuscript was written. All mathematical developments are the original work of the authors.

\medskip

%

\ms\noindent

\tableofcontents

\section{Introduction} 

Since the pioneer work of Dahlberg in the 1970s (\cite{Da77}), numerous mathematicians studied the solvability of the Dirichlet and Neumann problems with boundary data in $L^p$. The research area started from the Laplacian in Lipschitz domains, and continued with the goal of finding the necessary and sufficient conditions on the uniformly elliptic operator $L=-\diver A \nabla$, the domain $\Omega$, and the parameter $p\in (1,\infty)$ to ensure the $L^p$ solvability of boundary value problems.

Early in the 1980s, Jerison and Kenig derived a Rellich-type identity in their seminal papers \cite{JK81a,JK81b} that guarantees the $L^2$ solvability of both the Dirichlet and Neumann problems for the Laplacian in Lipschitz domains. This identity also extends to $t$-independent operators (see the introduction of \cite{KKPT00}).

For the Dirichlet problem, subsequent works have led to a characterization of the solvability of the $L^p$ Dirichlet problem for {\em some} $p\in (1,\infty)$, in terms of geometry and topology of the domain; see  \cite{HM14,HMU14,AHM+20} for the principal developments.

In contrast to the Dirichlet setting, the theory for the $L^p$
  Neumann problem remains fairly under–developed. Aside from the foundational results for Lipschitz domains and 
$t$-independent coefficients, progress has largely been limited to ``small-perturbation'' results; see for instance \cite{KP93, KP95, AAAHK11, DPR17, HMT10, MT??}, with some exceptions in the planar case $n=2$, see \cite{KR09, DHP23}. Even those perturbation results require {\em a priori} solvability of the $L^{p'}$ Dirichlet problem. Consequently, they do not provide genuine perturbation theorems for the $L^p$ Neumann problem itself; but they rather establish the stability of the equivalence \[ \|\partial_v u\|_{L^p(\partial \Omega)} \approx \|\nabla_t u\|_{L^p(\partial \Omega)},\]
where $\partial_v u$ and $\nabla_t u$ are respectively the normal and tangential derivative\footnote{The two derivatives need to be properly defined, depending on the geometric and elliptic setting in the considered article.} of $u$ at the boundary.

\medskip

Recently, the authors together with J. Zhuge looked for a characterization of the solvability of the $L^p$ Dirichlet problem for {\em all} $p\in (1,\infty)$, which ultimately led us to prove a stability result for the $L^p$ Dirichlet problem under small bi-Lipschitz transformations of domains (\cite{FLZ??}).

The aim of this paper is to establish an analogous stability result for the $L^p$ Neumann problem, ultimately paving the way for a geometric characterization of its solvability. While the methods we developed in 
\cite{FLZ??} are already quite involved, the Neumann problem poses significant additional challenges, necessitating a substantial amount of new computations and new ideas. A detailed discussion of our approach can be found in Section~\ref{S.mthd}. Notably, as the 
$L^p$ Neumann and the $L^{p'}$ Dirichlet problems are essentially independent (see \cite[Theorems 3.1 and 3.4]{KP95}), we deliberately refrain from assuming solvability of the $L^{p'}$ Dirichlet problem. To the best of our knowledge, the result presented here constitutes the first stability result for a Neumann problem that is obtained without such a hypothesis.

\medskip

Instead of proving the $L^p$ solvability of the Neumann problem, we show stability of the slightly stronger $L^p$ Poisson-Neumann problem. 
In \cite{FL-PN}, we introduced the Poisson-Neumann problem and established that solvability of the 
$L^p$ Poisson-Neumann problem implies that of the $L^{p'}$
  Neumann problem, and is itself implied by both the $
L^p$ Dirichlet and $L^{p'}$
 Neumann problems (see Theorem \ref{ThPNvsN} below).
The counterexamples in \cite[Theorems 3.1 and 3.4]{KP93} show that the combination of the $L^p$ Dirichlet and $L^{p'}$ Neumann problems is strictly stronger than the $L^{p}$ Neumann problem, but we do not know at the moment an example where the $L^{p'}$ Neumann problem is solvable while the $L^p$ Poisson-Neumann is not, making the Poisson-Neumann problem a close and fair substitute of the Neumann problem.
Moreover, we showed in \cite{FL-PN} that some basic properties on the $L^p$ Poisson-Neumann problem---interpolation and extrapolation---do not require assuming the solvability of the $L^p$ Dirichlet problem (contrary to the Neumann problem), making the Poisson-Neumann a better candidate for study.

\medskip

The article establishes the stability (for the Laplacian in Lipschitz domains) of the $L^p$ Poisson-Neumann problem under small bi-Lipschitz transformations of domains, without assuming the solvability of the $L^p$ Dirichlet problem (see Theorem~\ref{ThLipPert}). As part of the proof, we establish two other results of significant independent interest: 
\begin{enumerate}
    \item The stability of the $L^p$ Poisson-Neumann problem under small Carleson perturbations of the coefficients on 1-sided Chord Arc Domains (Theorem~\ref{ThCarlesonPert}). This improves the 30-year old result of Kenig and Pipher \cite[Theorem 2.2]{KP95}. 
    \item The $L^p$ solvability of the Poisson-Neumann problem for the Laplacian in Chord Arc Domains transfers to ``small special DKP operators'' (Theorem~\ref{MainTh}), which is a completely new result. The result closest to it is \cite[Theorem 7.1]{DPR17}, which establishes the $L^p$ solvability of the Neumann problem for small DKP operators in the half space.
\end{enumerate} Precise statements are given below.

\subsection{Main result}

\subsubsection{Domain.} Our main results hold when the domain $\Omega \subset \R^n$, $n\geq 2$, is bounded and either a Lipschitz domain, a 1-sided Chord Arc Domain (1-CAD) or a Chord Arc Domain (CAD), see Definition \ref{defCAD} for the later two notions.

\subsubsection{Operator and solutions.} The considered operator will be either the Laplacian, or a {\bf uniformly elliptic operator} in divergence form, that is, an operator $L=-\diver A\nabla$ for which the matrix function $A$ satisfies the ellipticity and boundedness conditions
\begin{equation} \label{defelliptic} 
A(X) \xi \cdot \xi \geq C^{-1} |\xi|^2 \text{ and } |A(X) \xi \cdot \zeta| \leq C |\xi||\zeta| \quad \text{ for a.e. $X\in \Omega$, and for $\xi, \zeta \in \R^n$.} \end{equation}
If $D\subset \Omega$ is an open set and if $h \in L^2_{loc}(D)$, $\vec F \in L^2_{loc}(D)$, we say that $u$ is a {\bf weak solution to $Lu = h - \diver \vec F$  in $D$} if $u\in W^{1,2}_{loc}(D)$ and
\begin{equation}  \label{defweaksol}
\int_{D} A\nabla u \cdot \nabla \varphi \, dX = \int_{D} [h \varphi + \vec F \cdot \nabla \varphi] \, dX \qquad \text{ for } \varphi \in C^\infty_0(D).
\end{equation}
We say that $u$ is a {\bf weak solution to $Lu = h - \diver \vec F$  in $D$ with zero Neumann data} if $u\in W^{1,2}_{loc}(D)$ and
\[\int_{D} A\nabla u \cdot \nabla \varphi \, dX = \int_{D} [h \varphi + \vec F \cdot \nabla \varphi] \, dX \ \text{ for } \varphi \in C^\infty(\overline{D}), \, \varphi \equiv 0 \text{ on a neighborhood of $\partial D \setminus \partial \Omega$}.\]
The most classical scenario of solutions with zero Neumann data will be when $D= \Omega$. In this case, a weak solution to $Lu = h - \diver \vec F$  in $\Omega$ with zero Neumann data exists in $W^{1,2}(\Omega)$ if and only if $h$ has zero average, that is, $\int_{\om}h\,dX=0$. The solution is unique up to a constant, and the (unique) solution with zero average can be represented with the help of a Neumann function (see Theorem \ref{ThNeumann}).

\subsubsection{Non-tangential and area functionals.} The {\bf distance to the boundary $\pom$} is denoted $\delta_\Omega$, or $\delta$ when the domain is clear, and for $X\in \Omega$, we write $B_X$ for the {\bf Whitney ball} $B(X,\delta_\Omega(X)/4)$. We define the cone with vertex at $x\in \partial\Omega$ as
\[ \Gamma(x):= \{X\in \Omega, \, |X-x| \leq 2 \delta_\Omega(X)\},\]
the averaged non-tangential function $\N^2$ as
\begin{equation} \label{defN2}
\N^2(u)(x) := \sup_{X\in \Gamma(x)} \left( \fint_{B_X} |u|^2 \, dY\right)^\frac12 \quad \text{ for } x\in \partial \Omega \text{ and } u\in L^2_{loc}(\Omega),
 \end{equation}
and the averaged area function as
\begin{equation} \label{defA12}
\A_1^2(u)(x) := \int_{\Gamma(x)}  \left( \fint_{B_X} |u|^2 \, dY\right)^\frac12 \frac{dX}{\delta_\Omega^n(X)} \quad \text{ for } x\in \partial \Omega \text{ and } u\in L^2_{loc}(\Omega).
\end{equation}
Many properties of the functionals can be found in Subsection \ref{SStents}.

\subsubsection{Poisson-Neumann problem.} The focus of the article is on the Poisson-Neumann problem, that is the Poisson problem in tent spaces for solution with zero Neumann data. The precise definition is given below.

\begin{definition}
Let $\Omega \subset \R^n$ be a bounded 1-sided Chord Arc Domain (equipped with its Ahlfors regular measure $\sigma$ on $\partial \Omega$), $L:=-\diver A \nabla$ be a uniformly elliptic operator on $\Omega$, and $p\in (1,\infty)$.
We say that the $L^p$ Poisson-Neumann problem is solvable--- \eqref{PNp} for short, and \eqref{PNp}$_{L}$, \eqref{PNp}$_{L,\Omega}$ when we want--- if there exists a constant $C_p$ such that for any $\vec F \in L^\infty_c(\Omega,\R^n)$, the weak solution $u \in W^{1,2}(\Omega)$ to 
$Lu=-\diver \vec F$ in $\Omega$ with zero Neumann data and zero average verifies
\begin{equation} \label{PNp} \tag{PN$_p$} \|\N^2(u)\|_{L^p(\partial \Omega,\sigma)} \leq C_p \|\A^2_1(\delta \vec F)\|_{L^p(\partial \Omega,\sigma)}.
\end{equation}
\end{definition}

We will not introduce precisely in this paper what we mean by solvability of the Dirichlet and Neumann problem. They will not be mentioned aside from the two next theorems, where they are used to position the Poisson-Neumann problem among the other boundary value problems. We refer to \cite[Definitions 1.1 and 1.8]{FL-PN} for the definitions which are used to prove the next theorem.

\begin{theorem}[{\cite[Proposition 1.8]{FL-PN}}]\label{ThPNvsN}
Let $\Omega \subset \R^n$ be a bounded 1-sided Chord Arc Domain, $L:=-\diver A \nabla$ be a uniformly elliptic operator on $\Omega$, and $p\in (1,\infty)$. Then:
\begin{itemize}
\item If the $L^p$ Dirichlet problem for $L$ and the $L^{p'}$ Neumann problem for $L^*$ are solvable, then the $L^p$ Poisson-Neumann problem for $L$ is also solvable.
\item  If the $L^p$ Poisson-Neumann problem for $L$ is solvable, then the $L^{p'}$ Neumann problem for $L^*$ is also solvable.
\end{itemize}
\end{theorem}

As mentioned in the introduction, both the Dirichlet and Neumann problem for the Laplacian in Lipschitz domains are solvable in $L^2$ (see \cite{JK81a,JK81b}). So, combined with Theorem \ref{ThPNvsN}, we have

\begin{theorem} \label{ThRellich}
Let $\Omega \subset \R^n$ be a bounded Lipschitz domain. Then the Poisson-Neumann problem {\rm (PN$_2$)}$_{-\Delta,\Omega}$ is solvable.
\end{theorem}

But we also know that the range of $p$ for which the $L^p$ Poisson-Neumann is solvable is an open interval that goes to infinity, see \cite[Corollary 1.21, Proposition 1.17]{FL-PN}. Therefore, we have

\begin{corollary} \label{CorRellich}
Let $\Omega$ be a bounded Lipschitz domain. There exists $p_-\in (1,2)$ such that the the Poisson-Neumann problem \eqref{PNp}$_{-\Delta,\Omega}$ is solvable if and only if $p\in (p_-,\infty)$.
\end{corollary}

\subsubsection{Stability under small bi-Lipschitz transformations} By Corollary \eqref{CorRellich}, for any Lipschitz domain $\Omega$, we can find a critical exponent $p_\Omega \in (1,2)$ such that the the Poisson-Neumann problem \eqref{PNp}$_{-\Delta,\Omega}$ is solvable if and only if $p\in (p_\Omega,\infty)$. However, we do not have any more precision about those value: we do not even know that $p_{\Omega_1}$ and $p_{\Omega_2}$ are close when $\Omega_1$ and $\Omega_2$ are ``close''.  
\smallskip

The question for the Dirichlet problem instead was solved recently by the authors together with J. Zhuge in \cite{FLZ??}: we established that 
the solvability of the $L^p$ Dirichlet problem is  stable under small bi-Lipschitz diffeomorphisms/transformations. 
The main goal of the article is to prove the analogous stability result for the Poisson-Neumann problem:

\begin{theorem} \label{ThLipPert}
Let $\Omega \subset \R^n$ be a bounded Lipschitz Domain. Let $p\in (1,\infty)$ be such that the Poisson-Neumann problem \eqref{PNp}$_{-\Delta,\Omega}$ is solvable. 

Then there exists $\epsilon>0$ (depending on $n$, $p$, the Lipschitz constant of $\Omega$, and the constant in the $L^p$-Poisson-Neumann problem) such that, if $\Phi:\, \R^n\to \R^n$ is a bi-Lipschitz change of variable whose Jacobian matrix $\Jac(\Phi)$ verifies
\begin{equation} \label{NPhi-I<e}
\|\Jac(\Phi) - I\|_{L^\infty} \leq \epsilon,
\end{equation}
then the Poisson-Neumann problem  \eqref{PNp}$_{-\Delta,\Phi(\Omega)}$ is also solvable.
\end{theorem}



\subsection{Our methods}\label{S.mthd}
Our proof relies on the ``better'' change of variable\footnote{Some modifications were needed to account for the fact that the domains has to be bounded, see Subsection \ref{SSMainThProof}.} $\rho$ constructed in \cite{FLZ??} that was used to prove the stability of the $L^p$ Dirichlet problem under small bi-Lipschitz transformations. This change of variable $\rho$ maps $\Omega$ to $\Phi(\Omega)$ --- like $\Phi$ --- and is such that $\rho(x) = \Phi(x)$ when $x\in \partial \Omega$. Since \eqref{PNp}$_{-\Delta,\Phi(\Omega)}$ is equivalent to \eqref{PNp}$_{L_\rho,\Omega}$, where $L_\rho$ is the conjugate of $-\Delta$ by $\rho^{-1}$, the map $\rho$ changes the question from the stability under domain perturbations to the stability under operator perturbations.

We will have three operators on $\Omega$ --- the Laplacian, $L_\rho$, and an intermediate operator $L_0$ --- such that $L_0$ is a small ``special DKP operator'', and $L_\rho$ is a small Carleson perturbation of $L_0$. It raises the need of two different perturbation theories for the $L^p$ Poisson-Neumann problem: one that takes care of the transition from the Laplacian to a small special DKP operator, and
another that shows the stability of the 
$L^p$ Poisson-Neumann problem under small Carleson perturbations.

 The stability result for the Dirichlet problem hinges crucially on the technical change of variable $\rho$ that we introduced in \cite{FLZ??}. Once this $\rho$ is carefully constructed, the proof unfolds relatively easily. In contrast, the Poisson-Neumann problem requires a substantial amount of additional computations and new ideas. While the approach in \cite{FLZ??} was characterized by a elegant geometric argument, the present article complements this theory with an additional analytical framework, often referred to as ``$S<N$ estimates". 

\subsubsection{Carleson measures.} Perturbations of operators, that preserves the solvability of the $L^p$ Dirichlet and Neumann problems, can take two forms for us, both described in terms of Carleson measures.

\medskip

First, we say that $f$ satisfies the {\bf Carleson measure condition} $f\in CM_1^\infty(M)$ if 
\begin{equation} \label{defCM1}
\sup_{x\in \partial \Omega, \, r\in (0,\diam \Omega)} r^{1-n} \int_{B(x,r) \cap \Omega}  \left(\sup_{B_X} |f| \right) \,\frac{dX}{\delta_\Omega(X)} \leq M,
\end{equation}
and $f\in CM_2^\infty(M)$ if 
\begin{equation} \label{defCM2}
\sup_{x\in \partial \Omega, \, r\in (0,\diam \Omega)} \left(r^{1-n} \int_{B(x,r) \cap \Omega}  \left(\sup_{B_X} |f|^2 \right) \,\frac{dX}{\delta_\Omega(X)} \right)^\frac12 \leq M.
\end{equation}
We write $f\in CM_1^\infty$ or $f\in CM_2^\infty$ when we do not need to specify the constant, and $\|f\|_{CM_1^\infty}$ or $\|f\|_{CM_2^\infty}$ for the left hand side of \eqref{defCM1}--\eqref{defCM2}. We let the reader check that
\[ \|f\|_{L^\infty} \lesssim \|f\|_{CM_2^\infty} \lesssim \|f\|_{CM_1^\infty}.\]

\subsubsection{$A_\infty$-absolute continuty of measures.}  When the domain is Lipschitz, the $L^2$ solvability of the Dirichlet problem for an operator $L$ corresponds to a $L^2$ reverse H\"older bound on the Poisson kernel $\frac{d\omega^X}{d\sigma}$, which in turn implies the $A_\infty$-absolute continuity of the elliptic measure with respect to the surface measure. When we work on non-Lipschitz domains, we will nonetheless require this later condition as a hypothesis.

The definition of $A_\infty$-absolute continuity is given here:

\begin{definition} \label{defAinfty}
Let $\mu,\nu$ be two measures supported on $E  \subset \R^n$. Then we say that $\nu$ is $A_\infty$-absolutely continuous with respect to $\mu$ --- or  $\nu \in A_\infty(\mu)$ for short --- if there exists $\alpha>0$ and $C>0$ such that for any $x\in E$, any $r>0$, and any Borel set $S\subset B(x,r) \cap \partial \Omega$, one has
\[ \frac{\mu(S)}{\mu(B(x,r)\cap E)} \leq  C \left( \frac{\nu(S)}{\nu(B(x,r)\cap E)}\right)^\alpha.\]
Note that despite the notation, $\nu \in A_\infty(\mu)$ is an equivalence relationship, meaning that 
\[\nu \in A_\infty(\mu) \iff \mu \in A_\infty(\nu).\]
\end{definition}

\subsubsection{Carleson perturbations.} The stability the $L^p$ boundary value problems under Carleson perturbations of the coefficients is well known. The stability for the $L^p$ Dirichlet problem is established in \cite{FKP91}, and the stability of the combination of the $L^p$ Neumann problem and $L^{p'}$ Dirichlet problem under small Carleson perturbations is proved in \cite[Theorem 2.2]{KP95}. 

Since the solvability of the $L^{p'}$-Dirichlet problem implies the $A_\infty$ absolute continuity of the elliptic measure with respect to the surface measure, see \cite[Theorem 1.7.3]{Ken94}, the next theorem is an improvement of  \cite{KP95}'s result:

\begin{theorem} \label{ThCarlesonPert}
Let $\Omega \subset \R^n$ be a bounded 1-sided Chord-Arc Domain, with the Ahlfors-regular measure $\sigma$. Let $L_0:=-\diver A_0 \nabla$ be a uniformly elliptic operators such that  $\omega_{L_0} \in A_\infty(\sigma)$, where $\omega_L$ is the elliptic measure of $L$ with pole at center defined in Definition \ref{defomG}. Let finally $p\in (1,\infty)$ be such that the the Poisson-Neumann problem \eqref{PNp}$_{L_0,\Omega}$ is solvable. 

Then there exists $\epsilon>0$ (depending on $n$, $p$, the elliptic and boundedness constants of $L_0$, the constant in \eqref{PNp}$_{L_0,\Omega}$, and the 1-sided CAD constants of $\Omega$) such that, if $L_1=-\diver A_1\nabla$ is an operator satisfying 
\begin{equation} \label{A1A0isCM}
A_0-A_1 \in CM_2^\infty(\epsilon)
\end{equation}
then the Poisson-Neumann problem \eqref{PNp}$_{L_1,\Omega}$ is also solvable. 

Moreover, for any $\delta>0$, we can choose $\epsilon$ small enough so that we can also make the constant $C_p = C_p(L_1,\Omega)$ in \eqref{PNp}$_{L_1,\Omega}$ to be smaller than $C_p(L_0,\Omega)+\delta$.
\end{theorem}

\subsubsection{Dahlberg-Kenig-Pipher (DKP) operators.} We say that a uniformly elliptic operator $L=-\diver A \nabla$ in $\Omega$ is a {\bf DKP operator} (for Dahlberg, Kenig and Pipher) if $\delta \nabla A \in CM_2^\infty$. Those operators were most notably studied in \cite{KP01} ($L^p$ solvability of the Dirichlet problem when $\|\delta \nabla A\|_{CM_2^\infty}$ is large), \cite{DPP07}  ($L^p$ solvability of the Dirichlet problem when $\|\delta \nabla A\|_{CM_2^\infty}$ is small), and \cite{DPR17} ($L^p$ solvability of the Neumann problem  when $\|\delta \nabla A\|_{CM_2^\infty}$ is small). 

When the domain is the half space $\R^n_+= \{(x,t) \in \R^{n-1}\times (0,\infty)\}$, there is a special class of DKP operators which are significantly easier to deal with: DKP operators for which the last row of $A$ is $(0,\dots,0,1)$.  We can reduce DKP operators to this class of DKP operators with the help of drift terms, like in \cite{KP01}, or with Carleson perturbation, like in \cite{FenDKP}. 

For general domains, we introduce the notion of {\bf special DKP operators} as those DKP operators that satisfy the condition \[(A-I)\nabla G_{-\Delta} = 0,\] where $G_{-\Delta}$ is the Green function for the Laplacian with a pole far from the boundary. This definition extends the previous one, which was tailored to the specific case of $\Omega = \R^n_+$, where the last row of $A$ 
is $(0,\dots,0,1)$, and 
$t$ is the Green function with a pole at infinity (see Definition \ref{defomG}).

A key component of our proof, and a notable result in its own right, is that solvability of the $L^p$ Poisson-Neumann problem for the Laplacian implies solvability for special DKP operators, provided the constants are sufficiently small. Specifically, we establish the following result.
\begin{theorem} \label{MainTh}
Let $\Omega \subset \R^n$ be a bounded Chord Arc Domain. Take $G$ to be the Green function for the Laplacian with pole at center (see Definition~\ref{defomG}). Let $p\in (1,\infty)$ be such that the Poisson-Neumann problem \eqref{PNp}$_{-\Delta,\Omega}$ is solvable.  Then there exists $\epsilon>0$ --- depending on $n$, $p$, the CAD constants of $\Omega$ and the constant in  \eqref{PNp}$_{-\Delta,\Omega}$ --- such that if $L=-\diver A \nabla$ is a uniformly elliptic operator satisfying
\begin{enumerate}[(i)]
\item $\|A-I\|_{L^\infty(\Omega)} \leq \epsilon$,
\item $\delta \nabla A \in CM_2^\infty(\epsilon)$, 
\item $(A-I)\nabla G = 0$ on $\Omega \setminus 2B_\Omega$,
\end{enumerate}
then the Poisson-Neumann \eqref{PNp}$_{L,\Omega}$ is also solvable.

Moreover, for any $\delta>0$, we can choose $\epsilon$ small enough so that we can also make the constant $C_p = C_p(L,\Omega)$ in \eqref{PNp}$_{L,\Omega}$ is smaller than $C_p(-\Delta,\Omega)+\delta$.
\end{theorem}

\subsubsection{Estimates on the square functional.} For the proof of the three main theorems, we established $L^p$ bounds on the square functional of solutions to the Poisson problem and their gradients, which the reader can find as Theorems \ref{ThS<N11} and \ref{ThS<N21}. Those estimate are cornerstones of our argument, and are used in combination with a fixed point arguments to prove Theorems \ref{ThCarlesonPert} and \ref{MainTh}.

\subsection{Plan of the proof} Below is a diagram to show the implications between the different results of the manuscript. 

Theorems \ref{ThS<N1} and \ref{ThS<N2} are placed on the same level because their proofs are similar, and both rely on an integration parts argument, with Theorem \ref{ThS<N2} being the more complicated version of Theorem \ref{ThS<N1}. Combined with a good-$\lambda$ argument, they self-improve into Theorems \ref{ThS<N11} and \ref{ThS<N21} respectively.

Theorems \ref{ThCarlesonPert} and \ref{MainTh} are also placed on the same line, because they come from respectively Theorems \ref{ThS<N11} and \ref{ThS<N21} coupled with a fixed point argument. As before,  Theorem \ref{ThCarlesonPert} can be seen as the baby version of Theorem \ref{MainTh}.

Theorems \ref{ThCarlesonPert} and \ref{MainTh} are used in combination with the change of variables from \cite{FLZ??} to prove the last main theorem, Theorem \ref{ThLipPert}.

The three results greyed out in the forth column are important results proved in earlier articles, two of them being from the authors themselves.

The two lemmas at the top are the most important lemma needed for the theorem below them. None of them are fundamentally complicated: the first one is proven by using the Caccioppoli argument twice, and the second one explains how to do ``add derivatives'' in the integration by part argument when the domain is not the half space. Note however that the second one --- Lemma \ref{IPPwithNG} --- is used twice at key moments, might be the most original idea of the article, and should also be of independent interest to other researchers in the field.

\pagebreak

\begin{tikzcd}[column sep=small]
\ffbox{Lemma \ref{lemDCaccio} }{Double Caccioppoli inequality} \arrow[dd,Rightarrow] &  & \ffbox{Lemma \ref{IPPwithNG}}{Integration by parts w.r.t. $\nabla G$} \arrow[dd,Rightarrow] \arrow[dddddd,Rightarrow, shift right = 0.7cm,  bend right = 22] &  \\ \\
\ffbox{Theorem \ref{ThS<N1}}{Local Poisson $S<N$ estimates} \arrow[dd,Rightarrow] & \sfbox{Theorem \ref{Thgoodlambda}}{Good-$\lambda$} \arrow[ddr,Rightarrow, crossing over, bend right=25] \arrow[ddl,Rightarrow, bend left=25] 
& \ffbox{Theorem \ref{ThS<N2}}{Local Poisson-regularity $S<N$ estimates} \arrow[dd,Rightarrow] &  
\\ & & &  \\
\ffbox{Theorem \ref{ThS<N11}}{$L^p$ Poisson $S<N$ estimates for all $p$} \arrow[dd,Rightarrow]  \arrow[ddrr,Rightarrow, bend right = 15] &  & \ffbox{Theorem \ref{ThS<N21}}{$L^p$ Poisson-regularity $S<N$ estimates for all $p$} \arrow[dd,Rightarrow] & \fffbox{Theorem \ref{ThUR=>Ainfty} and Theorem \ref{ThonG}}{\footnotesize CAD $\implies \omega \in A_\infty(\sigma)$ $ + \delta^2 \nabla^2 G/G \in CM^\infty_2$} \arrow[ddl,Rightarrow, bend left = 15] \arrow[uul,Rightarrow, bend right = 15] \\ \\
\ffbox{Theorem \ref{ThCarlesonPert}}{Stability of \eqref{PNp} under small Carleson perturbations} \arrow[ddrr,Rightarrow, bend right = 15] &  & \ffbox{Theorem \ref{MainTh}}{Stability of \eqref{PNp} for small special DKP operators} \arrow[dd,Rightarrow] & \sffbox{Theorem \ref{ThLrho}}{Change of variables} \arrow[ddl,Rightarrow, bend left = 15] \\ \\
 &  & \ffbox{Theorem \ref{ThLipPert}}{Stability of \eqref{PNp} under small bi-Lipschitz transformations} &   \\ \\
\end{tikzcd}

\subsection{Extension to unbounded domains.} The main three theorems, Theorems \ref{ThLipPert}, \ref{ThCarlesonPert}, and \ref{MainTh}, are stated with the condition that the domain $\Omega \subset \R^n$ is bounded, and it is natural to ask whether this could be extended to unbounded domains $\Omega$.

We claim that the answer is positive, as long as $n\geq 3$. Indeed, in our article, we need to cancel/remove the average of functions in $\Omega$, like in \eqref{defhatW12} or \eqref{eq.RieszNmn'}, to have a suitable Hilbert space for the application of the Lax-Milgram theorem. This average of functions on $\Omega$ is then kept for the definition of the (Poisson-)Neumann problem and the Neumann function.

When $n \geq 3$ and the domain supports suitable Poincar\'e inequalities, one can define a Hilbert space $W(\Omega)$ equipped with the norm $\|\nabla .\|_{L^2}$ by either defining it as $W(\Omega) := \dot W^{1,2}(\Omega) \cap L^{2n/(n-2)}(\Omega)$, or $W(\Omega) := \{ f\in \dot W^{1,2}(\Omega), \, \lim_{N\to \infty} \fint_{B(X_\Omega,N)\cap \Omega} f \, dX = 0\}$. This constructed Sobolev space is suitable for the use of the Lax-Milgram theorem, and allows for a construction of the Neumann function, see \cite{HS25} for details. From there, one needs to extend the properties of the $L^p$ Poisson-Neumann from \cite{FL-PN} and then the current results to unbounded 1-sided Chord Arc Domains in $\R^n$, $n\geq 3$, but we claim that the modifications are minor\footnote{With unbounded domains in $\R^n$, $n\geq 3$, one do not need to remove/cancel the averages of function, which makes all the computations actually easier.}.

For unbounded planar domains $\Omega \subset \R^2$, we do not know if/how the Poisson-Neumann problem can be defined, and if/how our theory can be adapted. 

\section{Preliminary results}

\subsection{Topology conditions on the domain.}

\begin{definition} \label{defAR}
We say that a measure $\sigma$ on $\R^n$ is a {\bf $d$-Ahlfors regular measure} if there exists $C>0$ such that 
\[ C^{-1} r^{d} \leq \sigma(B(X,r)) \leq C r^d \qquad \text{ for } X\in \supp \sigma, \, r\in (0, \diam E).\]
We say that $E \subset \R^n$ is a $d$-Ahlfors regular set if there exists a $d$-Ahlfors regular measure whose support is $E$. In this case, we say that $\sigma$ is an Ahlfors regular measure on $E$.
\end{definition}

\begin{proposition}
Let $\sigma$ a measure on $\R^n$. If $\sigma$ is $d$-Ahlfors regular, then there is a constant $C>0$ such that 
\[ C^{-1} \mathcal H_{|\supp \sigma}^d \leq \sigma \leq C \mathcal H_{|\supp \sigma}^d\]
where $\mathcal H_{|\supp \sigma}^d$ is the $d$-Hausdorff measure on the support of $\sigma$. 
\end{proposition}

\bp
The result is simple, it suffices to cover of $B(x,r) \cap \supp \sigma$ by (a uniformly finite number of) balls with smaller radius. The proof can be found as \cite[Lemma 1.8]{FenHDR}.
\ep

\begin{definition} \label{defCP}
For $x\in \partial \Omega$, $r>0$ and $\epsilon >0$, we say that $X\in \Omega$ is a {\bf $\epsilon$-Corkscrew point} for $(x,r)$ if $X\in B(x,r)$ and $B(X,\epsilon r) \subset \Omega$.

We say that an open set $\Omega \subset \R^n$ satisfies the {\bf Corkscrew point condition} if there exists $\epsilon_{CP}>0$ such that there exists an $\epsilon_{CP}$-Corkscrew point for any couple $(x,r) \in \partial \Omega \times (0,\diam \Omega)$. When the domain satisfies the Corkscrew point condition, we say that $X$ is a Corkscrew point for $(x,r) \in \partial \Omega \times (0,\diam \Omega)$ if  $X$ is a $\epsilon_{CP}$-Corkscrew point  $(x,r) \in \partial \Omega \times (0,\diam \Omega)$.

When $\Omega$ is bounded and satisfies the Corkscrew point condition, there is a $X_\Omega$ satisfying 
\[B(X_\Omega,\epsilon_{CP}\diam(\Omega)) \subset \Omega.\] 
We call $X_\Omega$ {\bf the center of $\Omega$} and $B_\Omega := B(X_\Omega,\delta_\Omega(X_\Omega)/16)$ the {\bf center ball of $\Omega$}.
\end{definition}

\begin{definition} \label{defHC}
We say that $\Omega$ satisfies the Harnack chain condition if for any $\Lambda \geq 1$, there exists $N_\Lambda\in \bN$ such that the following holds. For any $X,Y\in \Omega$ verifying
\[ \dfrac{|X-Y|}{\min\{\delta_\Omega(X),\delta_\Omega(Y)\}} \leq \Lambda,\]
there exist $N_\Lambda+1$ points $\{X_0,\dots,X_{N_\Lambda}\}$ in $\Omega$ that verify
\begin{enumerate}
\item $X_0=X$ and $X_{N_\Lambda} = Y$;
\item for each $1\leq i \leq N_\lambda$, $|X_{i-1} - X_i| \leq \delta_\Omega(X_i)/4$.
\end{enumerate}
Here, as before, $\delta_\Omega(Z)$ stands for $\dist(Z,\partial \Omega)$. We say that the collection of points $X_0,\dots,X_{N_\lambda}$ is a Harnack chain of points linking $X$ to $Y$. 
\end{definition}

If $\Omega$ satisfies both the Harnack chain condition and the Corkscrew point condition if and only if $\Omega$ is uniform, which means that we have the following proposition.

\begin{proposition}
Let $\Omega \subset \R^n$ be open. The domain $\Omega$ satisfies the Corkscrew point condition and the Harnack chain condition if and only if $\Omega$ is uniform, i.e. there exist $\epsilon,C>0$ such that the following holds. For any $X,Y\in \overline{\Omega}$, there exists a rectifiable curve $\gamma: [0,L] \to \overline{\Omega}$ joining $X$ and $Y$ and parametrized by arc length such that 
\begin{enumerate}[(i)]
\item the length $L$ of the rectifiable curve is smaller than $C|X-Y|$,
\item for any $t\in (0,L)$, 
\[ \delta_\Omega(\gamma(t)) \geq \epsilon \min\{t,L-t\}. \]
\end{enumerate}
\end{proposition}

\bp
The result is classical, although the authors cannot pinpoint where it was proved. The proposition is only given to help the readers relate with the different characterizations of 1-sided Chord Arc Domain. The proof is sketched below.

\medskip

\noindent $(\Rightarrow)$ Let $\Omega$ satisfy the Corkscrew point and Harnack chain condition, then for $X,Y\in \overline{\Omega}$, take $x,y\in \partial \Omega$ such that $|X-x| = \delta_\Omega(X)$ and $|Y-y| = \delta_\Omega(Y)$. For $i\in \mathbb Z$ satisfying $i_- := \ln_2(\delta_\Omega(X)) < i \leq i_+:= \ln_2(|X-Y|)$, let $X_i$ to be a Corkscrew point for $(x,2^i)$ and  for integers $j\in \mathbb Z$ satisfying $j_- < j \leq j_+:= \ln_2(|X-Y|)$, let $Y_i$ be a Corkscrew point for $(y,2^j)$. 
With the Harnack chain condition, we can find a uniform number $N\in \mathbb N$,  such that each of the following couple: $(X,X_0)$, $(X_i, X_{i+1})$ for $i_- < i \leq i_+-1$, $(X_{i_0},Y_{j_0})$, $(Y_{j},Y_{j-1})$ for $j_-+1 < j \leq j_0$, and finally $(Y_{0},Y)$, can be linked by an Harnack chain of points of length at most $N$. Juxtaposing the chains gives a chain $\{Z_k\}_{k\in K\subset \mathbb Z}$ --- of infinite length when either $\delta_\Omega(X)=0$ or $\delta_\Omega(Y) = 0$ --- linking $X$ and $Y$. The rectifiable curve is obtained by concatenating the segments $[Z_k,Z_{k+1}]$, and we let the reader check that it has the desired properties.

\medskip

\noindent  $(\Leftarrow)$ Let $\Omega$ be uniform, and $X,Y \in \overline{\Omega}$. Let $\gamma:\, [0,L]$ be the rectifiable curve linking $X$ and $Y$ as in the proposition.  We construct the collection of points $X_i$ by induction such that $X_0= X$, and $X_{i+1}\in \gamma([0,L])$ be such that $|X_{i+1}-X_i| = \delta_\Omega(X_i)/4$ or $X_{i+1} = Y$. We verify that this collection of points forms a Harnack chain of points between $X$ and $Y$.

If $x\in \partial \Omega$ and $r\in (0,\diam \Omega)$, we pick $Y\in \overline{\Omega} \setminus B(x,r)$ and let $\gamma;\, [0,L]$ be the rectifiable curve linking $x$ to $Y$. We then let the reader check that $\gamma(r/2)$ is a Corkscrew point for $(x,r)$.
\ep

\begin{definition} \label{defCAD}
Let $\Omega \subset \R^n$. We say that $\Omega$ is a 1-sided Chord Arc Domain (1-CAD) if $\Omega$ satisfies the Corkscrew point condition and the Harnack chain condition, and if $\partial \Omega$ is $(n-1)$-Ahlfors regular.

We say that $\Omega$ is a Chord Arc Domain (CAD) if $\Omega$ is 1-CAD and $\R^n \setminus \overline{\Omega}$ satisfies the Corkscrew point condition.

The CAD and 1-CAD constants are the constants that appears in Definitions \ref{defAR}, \ref{defCP} and \ref{defHC}.
\end{definition}

Chord Arc Domains are 1-CAD with uniformly rectifiable boundaries, as shown in \cite[Theorem 1.2]{AHMNT17}

\begin{proposition} \label{prUR<=>CAD}
Let $\Omega \subset \R^n$ be a 1-sided Chord Arc Domain. Are equivalent:
\begin{enumerate}[(i)]
\item $\R^n \setminus \overline{\Omega}$ satisfies the Corkscrew point condition, i.e. $\Omega$ is a Chord Arc Domain;
\item $\partial \Omega$ is uniformly rectifiable, i.e. there exists $\epsilon>0$ such that for any $x\in \partial \Omega$ and any $r\in (0,\diam \partial \Omega)$, there exists a Lipschitz function $f_{x,r}:\, \R^{n-1} \to \R^n$ such that 
\[ \mathcal H^{n-1}_{|\partial \Omega}(B(x,r) \cap f_{x,r}(\R^{n-1})) \geq \epsilon r^{n-1}.\]
\end{enumerate}
\end{proposition}

The setting of 1-sided Chord Arc Domains allows for a nice elliptic theory, in particular a boundary Poincar\'e inequality, and --- if a uniformly elliptic operator $L=-\diver A \nabla$ on $\Omega$ is considered --- the Caccioppoli and Moser inequalities (inside the domain and at the boundary), as well as the existence of an elliptic measure and Green function. The presentation of the elliptic theory and its proof can be found for instance in \cite{DFMmixed}.

\subsection{Elliptic measure and Green function with pole at center/infinity}

\begin{proposition} \label{propGinfty}
Let $\Omega$ be an unbounded 1-CAD and $L:=-\diver A \nabla$ be a uniformly elliptic operator on $\Omega$. Then there exists a positive solution $G_{L^*} \in C^0(\overline \Omega) \cap W^{1,2}_{loc}(\Omega)$ to $L^* u= -\diver A^T \nabla u = 0$ in $\Omega$ with $G_{L^*}\equiv 0$ on $\partial \Omega$. Moreover such $G_{L^*}$ is unique up to multiplication by a constant.

If $G_{L^*} \in C^0(\overline \Omega) \cap W^{1,2}_{loc}(\Omega)$ is as above, then there exists a positive measure $\omega_{L}$ on $\partial \Omega$ such that 
\[\int_{\Omega} \nabla G_{L^*} \cdot \nabla \varphi \, dX = - \int_{\partial \Omega} \varphi d\omega_{L} \qquad \text{ for } \varphi \in C_0^\infty(\R^n).\]
\end{proposition}

\bp 
See \cite[Lemma 3.2]{DFMpert}, and earlier proofs in \cite[Lemma 3.7, Corollary 3.2]{KT99} and \cite[Lemma 6.5]{DEM21}.
\ep

\begin{definition} \label{defomG}
Let $\Omega$ be a 1-CAD.
\begin{enumerate}[(i)]
\item If $\Omega$ is bounded, let $B_\Omega$ be the center ball introduced in the Corkscrew point condition (Definition \ref{defCP}). Then we define the {\bf Green function (for $L^*$) with pole at center} as
\[ G_{L^*} (X) := \fint_{B_{\Omega}/2} G_{L^*}(X,Y) dY,\]
where $G_{L^*}(.,Y) $ is the Green function (for $L^*$) with pole at $Y$, and the {\bf elliptic measure (for $L$)  with pole at center} as $\omega_L := \omega_L^{X_0},$ where $\omega_L^Y$ is the elliptic measure for $L$ with pole at $Y$.

\item If $\Omega$ is unbounded, then we say that $(G_{L^*},\omega_L)$ are the {\bf Green function (for $L^*$) and the elliptic measure (for $L$) with pole at infinity} if the couple satisfies Proposition \ref{propGinfty}. They are thus defined up to a constant.
\end{enumerate}
To avoid the distinction of cases, we call {\bf Green function or elliptic measure with pole at infinity/center} when we do not know whether the domain $\Omega$ is bounded or unbounded.
\end{definition}

\begin{theorem} \label{ThCP}
Let $\Omega$ be a 1-CAD and $L:=-\diver A \nabla$ be a uniformly elliptic operator on $\Omega$. Let $(G_{L^*},\omega_L)$ be respectively the Green function for $L^*$ and elliptic measure for $L$ with pole at infinity/center. Then there exists $C>0$ (depending only on the 1-CAD constants and the ellipticity/boundedness of the coefficients of $L$) such that, for any $(x,r) \in \partial \Omega \times (0,\diam \Omega)$, and any Corkscrew point $X$ for $(x,r)$, we have
\[ C^{-1} r^{n-2} G_{L^*}(X)  \leq \omega_L(B(x,r)\cap \partial \Omega) \leq C r^{n-2} G_{L^*}(X).\]
\end{theorem}

\bp
The above estimate is the Caffarelli-Fabes-Mortola-Salsa estimate applied to the Green function and harmonic measure with pole at center/infinity. The proof can be found as \cite[Lemma 3.3]{DFMpert} when $\Omega$ is unbounded, and the earlier CFMS estimates on the harmonic measure and Green function with finite pole are found in \cite[Lemma 15.28]{DFMmixed} and \cite{CFMS81}.  
 \ep

\subsection{Tents spaces} \label{SStents}

When $X\in \Omega$, we recall that we write $\delta(X)$ for $\delta_\Omega(X):= \dist(X,\Omega)$ and $B_X$ for $B(X,\delta(X)/4)$.
In the next results, we always assume that $\Omega \subset \R^n$ is a 1-sided Chord Arc Domain because it is the setting of our article, and not because it is optimal in any way. For instance, the connectedness condition (Harnack chain condition) has actually no play in any of the following statements where the harmonic measure is not involved.

\begin{definition}
Let $\Omega \subset \R^n$. 
\begin{itemize}
\item For $\gamma >1$ and $x\in \partial \Omega$
\[ \Gamma_{(\gamma)}(x) := \{ X\in \Omega, \, |X-x| \leq \gamma \delta(X)\}.\]
\item For $\alpha\in [1,\infty)$, $\gamma >1$, $\beta \in [1,\infty]$ and $\kappa \in (0,4)$, we define the quantity $\A_{\alpha,(\gamma,\kappa)}^\beta$ as
\[ \A_{\alpha,(\gamma,\kappa)}^\beta(f)(x):= \left( \int_{\Gamma_{(\gamma)}(x)}  \|f\|_{\text{\L}^\beta(\kappa B_X)}^{\alpha} \, \frac{dX}{\delta(X)^n} \right)^\frac1\alpha \qquad \text{ for $f\in L^\beta_{loc}(\Omega)$ and $x\in \partial \Omega$},\]
where 
\[\|f\|_{\text{\L}^\beta(\kappa B_X)}
:=\br{\frac1{|\kappa B_X|}\int_{\kappa B_X}|f(Z)|^\beta\,dZ}^{1/\beta} \quad\text{when }\beta\in[1,\infty)
\] 
and $\|f\|_{\text{\L}^\beta(\kappa B_X)}
:=\displaystyle\sup_{\kappa B_X}|f|$ when $\beta=\infty$. 
\item When $\alpha = \infty$, we like to use the notation $\N$, and so for $\gamma >1$, $\beta \in [1,\infty]$ , and $\kappa \in (0,4)$, we define
\[ \A^{\beta}_{\infty,{(\gamma,\kappa)}}(f)(x) = \N_{(\gamma,\kappa)}^\beta(f)(x):= \sup_{X\in \Gamma_{(\gamma)}(x)}\norm{f}_{\text{\L}^{\beta}(\kappa B_X)}\quad \text{ for $f\in L^\beta_{loc}(\Omega)$ and $x\in \partial \Omega$}.\]
\item As we shall see in (i) of Lemma \ref{lemcompNA}, the parameters $\gamma$ and $\kappa$ are not relevant in our theory. We write then $\Gamma(x)$, $\A_\alpha^\beta$, and $\N^\beta$ for $\Gamma_{(2)}(x)$, $\A_{\alpha,(2,1)}^\beta$, and $\N^\beta_{(2,1)}$ respectively.
\item Finally, $\A$ is $\A_2^2$ and $\N$ is $\N^\infty$.
\end{itemize}
\end{definition}

\begin{definition}
A measure $\mu$ on $\partial \Omega$ is a doubling  if
\[ C_{\mu}:= \sup_{x\in \partial \Omega, r>0} \frac{\mu(B(x,2r))}{\mu(B(x,r))} <+\infty,\]
and $C_\mu$ is called the doubling constant. One can easily check that Ahlfors regular measures are doubling.
\end{definition}

\begin{lemma} \label{lemFubini}
Let $\Omega$ be a $1$-sided CAD, with Ahlfors regular measure $\sigma$. 
\begin{enumerate}[(i)]
\item Then, for any $p\in [1,\infty]$, and any $\gamma >1$,
\begin{equation} \label{FubiniforA}
\|\mathcal A_{p,(\gamma,\kappa)}^p(f)\|_{L^{p}(\partial \Omega,\sigma)}^p  \approx  \int_\Omega |f(X)|^p \, \frac{dX}{\delta(X)},
\end{equation}
where the constant depends on $p$, $\gamma$, and the doubling constant of $\sigma$.
\item If $L:=-\diver A \nabla$ is a uniformly elliptic operator on $\Omega$, and if $(G_{L^*},\omega_L)$ are the Green function for $L^*$ and the elliptic measure for $L$ with pole at center/infinity, then 
\begin{equation} \label{FubiniAtoG}
\|\mathcal A_{p,(\gamma,\kappa)}^p(f)\|_{L^{p}(\partial \Omega,\omega_L)}^p  \approx  \int_\Omega |f|^p\, \frac{G_{L^*}}{\delta^2} \, dX,
\end{equation}
where the constant depends on $p$, $\beta$, $\gamma$, $\kappa$, and the 1-sided CAD constants of $\Omega$, and the ellipticity/boundedness constant of the coefficients of $L$
\end{enumerate}
\end{lemma}

\bp
The first part is a simple consequence of the Fubini's theorem.

The second part, on top of using Fubini's theorem, also use the comparison between the Green function and the harmonic measure given in Theorem \ref{ThCP}. The case $p=2$ is found as Equation (5.15) in \cite{FP23}, and the general case --- which is not more difficult --- is left to the reader.
\ep

\begin{lemma} \label{lemcompNA}
Let $\Omega$ be a $1$-CAD, and let $\mu$ be a doubling measure on $\partial \Omega$.
\begin{enumerate}[(i)]
\item For any $\alpha,\beta \in [1,\infty]$, $\gamma >1$, $\kappa \in (0,4)$, and $p\in (0,+\infty)$, we have for any $u\in L^\infty_c(\Omega)$,
\[ \big\| \mathcal A_\alpha^\beta(u) \big\|_{L^p(\partial \Omega,\mu)} \approx  \big\| \mathcal A_{\alpha,(\gamma,\kappa)}^\beta(u) \big\|_{L^p(\partial \Omega,\mu)} , \qquad \big\| \mathcal N^\beta(u) \big\|_{L^p(\partial \Omega,\mu)} \approx  \big\| \mathcal N_{(\gamma,\kappa)}^\beta(u) \big\|_{L^p(\partial \Omega,\mu)},\]
where the constants depends only on the relevant parameters (the doubling constant of $\mu$, $\alpha$, $\beta$, $\kappa$, $\gamma$, and $p$), and where the finiteness of one side guarantees the finiteness of the other side.

\medskip

\item  If $1\leq \alpha' \leq \alpha \leq \infty$ and $1 \leq \beta \leq \beta' \leq \infty$, then for any $p\in (0,\infty)$, 
\[ \| \mathcal A_\alpha^\beta(u) \big\|_{L^p(\partial \Omega,\mu)} \lesssim \| \mathcal A_{\alpha'}^{\beta'}(u) \big\|_{L^p(\partial \Omega,\mu)}  \qquad \text{ for $u\in L^{\beta'}_{loc}(\Omega)$.}\]
\end{enumerate}
\end{lemma}

\begin{definition} \label{defTp}
For $p\in (1,\infty)$, $\alpha,\beta \in [1,\infty]$ and $s\in \mathbb R$, we define the tent spaces $T^{p,\alpha,\beta}_{s}(\Omega)$ as
\[ T^{p,\alpha,\beta}_{s}(\Omega) := \{h \in L^2_{loc}, \, \|h\|_{T^{p,\alpha,\beta}_s(\Omega)}:= \|\mathcal A^\beta_\alpha(\delta^{-s}h)\|_{L^p(\partial \Omega)} < +\infty\}.\]
Lemma \ref{lemcompNA} shows that the spaces $T^{p,\alpha,\beta}_{s}(\Omega)$ do not depend on the choice of $\gamma$ and $\kappa$.
\end{definition}

\noindent{ \em Proof of Lemma \ref{lemcompNA}: }
Those results are fairly classical, although it will be hard to pinpoint a literature proving the result in this form, as they are mostly proved in the Euclidean setting (but the proofs extend easily to our cases). It is partially found in \cite[Paragraph II.2.5]{Stein93}, as a combination of \cite[Lemmas 2.1 and 2.2]{MPT} and \cite[Theorem 3]{CMS85}, or in the Appendix of \cite{Mando26}. Let us give some hints on how the proof would work, if one wants to prove it by themselves. For (i), notice that
\begin{itemize}
\item The parameter $\kappa$ is easy to handle.  When $\kappa \leq \kappa'$, we easily have
\[\mathcal A_{\alpha,(\gamma,\kappa)}^\beta(u)(x)  \lesssim \mathcal A_{\alpha,(\gamma,\kappa')}^\beta(u)(x),\]
and when $\kappa \geq \kappa'$, we cover $\kappa B_X$ by a uniformly finite number of balls of the form $\kappa'B_Y$, with $Y \in \Gamma_{(\gamma')}(x)$ and $\gamma'$ large enough independent of $x$, and then we have 
\[\mathcal A_{\alpha,(\gamma,\kappa)}^\beta(u)(x)  \lesssim \mathcal A_{\alpha,(\gamma',\kappa')}^\beta(u)(x).\]
So by increasing $\gamma$, we can change $\kappa$ to $\kappa'$.
\item When $p=\alpha$, the equivalence for the different value of $\gamma$ is a consequence of Fubini's theorem. This is a variant of Lemma \ref{lemFubini} (i), where we consider general doubling measures $\mu$ instead of the only Ahlfors regular measures.
\item When $\alpha \neq p$ is  the equivalence for the different value of $\gamma$ is obtained from the case $\alpha = p$ and a good-lambda argument (which, we admit, is not trivial but classical). The proof of Theorem \ref{Thgoodlambda} below --- which is a good-lambda argument type of proof  --- does not rely on changing the aperture of the cones $\Gamma(x)$, and an interested reader can adapt the proof of Theorem \ref{Thgoodlambda} to prove the equivalence for the different values of $\gamma$.
\end{itemize}
Finally, for (ii), the key is to remark that there are  $1 < \gamma_1< \gamma_2$ and $0 < \kappa_1 < \kappa_2 < 4$ and a constant $C$ such that  
\begin{equation} \label{defbarAba}
 C^{-1} \A_{\alpha,(\gamma_1,\kappa_1)}^\beta \leq \bar \A_\alpha^\beta \leq C^{-1} \A_{\alpha,(\gamma_2,\kappa_2)}^\beta,
 \end{equation}
where
\[ \bar \A_\alpha^\beta(u)(x) := \left( \sum_{k\in \mathbb Z} \norm{u}_{\text{L}^\beta(\Gamma^k(x))}^{\alpha} \right)^\frac{1}\alpha,\]
$\Gamma^k(x) := \{X\in \Gamma(x), \, 2^{k-1} \leq  \delta(X)  \leq 2^{k+1}\}$.
The quantity $\bar \A_\alpha^\beta(u)(x)$ clearly increases with $\beta$ (H\"older inequality) and decreases with $\alpha$ (comparison of $\ell^\alpha(\mathbb Z)$ spaces). With the first part of the Lemma and \eqref{defbarAba}, it means that 
\[\|\A_{\alpha}^\beta(u)\|_{L^p(\partial \Omega,\mu)} \approx \|\bar \A_{\alpha}^\beta(u)\|_{L^p(\partial \Omega,\mu)} \]
also  increases with $\beta$ and decreases with $\alpha$, as desired.
\ep

\begin{proposition} \label{prDualityNA1}
Let $\Omega \subset \R^n$ be a 1-CAD, and let $\mu$ be a doubling measure on $\partial \Omega$. 
\begin{enumerate}[(i)]
\item If $(\alpha,\alpha_1,\alpha_2) \in [1,\infty]^3$, $(\beta,\beta_1,\beta_2) \in [1,\infty]^3$, and $(p,p_1,p_2) \in [1,\infty]^3$ are such that 
\[ \frac1\alpha = \frac1{\alpha_1} + \frac1{\alpha_2}, \qquad \frac1\beta = \frac1{\beta_1} + \frac1{\beta_2}, \quad \text{ and } \quad \frac1p = \frac1{p_1} + \frac1{p_2}\]
then, for any $u\in L^{\beta_1}_{loc}(\Omega)$ and any $v\in L^{\beta_2}_{loc}(\Omega)$
\[  \|\A_\alpha^\beta(uv)\|_{L^p(\partial \Omega,\mu)} \leq \|\A_{\alpha_1}^{\beta_1}(u)\|_{L^{p_1}(\partial \Omega,\mu)} \|\A_{\alpha_2}^{\beta_2}(v)\|_{L^{p_2}(\partial \Omega,\mu)}.\]
\end{enumerate}
Let us record two special cases, 
\begin{enumerate}[(i)] \addtocounter{enumi}{1}
\item If $\sigma$ is an Ahlfors-regular measure on $\partial \Omega$, then there exists $C_\sigma>0$ such that  for any $(\alpha_1,\alpha_2) \in [1,\infty]^2$, any $(\beta_1,\beta_2) \in [1,\infty]^2$, and any $(p_1,p_2) \in [1,\infty]^2$ verifying 
\[ 1 = \frac1{\alpha_1} + \frac1{\alpha_2} = \frac1{\beta_1} + \frac1{\beta_2} = \frac1{p_1} + \frac1{p_2},\]
and then for any $u\in L^{\beta_1}_{loc}(\Omega)$ and any $v\in L^{\beta_2}_{loc}(\Omega)$, we have
\[  \left| \int_\Omega uv\, \frac{dX}{\delta} \right| \leq C_\sigma \|\A_{\alpha_1}^{\beta_1}(u)\|_{L^{p_1}(\partial \Omega,\sigma)} \|\A_{\alpha_2}^{\beta_2}(v)\|_{L^{p_2}(\partial \Omega,\sigma)},\]
\item If $L:=-\diver A \nabla$ is a uniformly elliptic operator on $\Omega$, and if $(G_{L^*},\omega_L)$ are the Green function for $L^*$ and the elliptic measure for $L$ with pole at center/infinity, then there exists $C>0$ such that  for any $(\alpha_1,\alpha_2) \in [1,\infty]^2$, any $(\beta_1,\beta_2) \in [1,\infty]^2$, and any $(p_1,p_2) \in [1,\infty]^2$ verifying 
\[ 1 = \frac1{\alpha_1} + \frac1{\alpha_2} = \frac1{\beta_1} + \frac1{\beta_2} = \frac1{p_1} + \frac1{p_2},\]
and then for any $u\in L^{\beta_1}_{loc}(\Omega)$ and any $v\in L^{\beta_2}_{loc}(\Omega)$, we have
\[  \left| \int_\Omega uv\, \frac{G_{L^*}}{\delta^2} \, dX \right| \leq C \|\A_{\alpha_1}^{\beta_1}(u)\|_{L^{p_1}(\partial \Omega,\omega_L)} \|\A_{\alpha_2}^{\beta_2}(v)\|_{L^{p_2}(\partial \Omega,\omega_L)}.\]
\end{enumerate}
\end{proposition}

\bp 
Point (i) is only an application of the H\"older inequality. Points (ii) and (iii) are the combination of (i) and Lemma \ref{lemFubini}.
\ep

\subsection{Carleson inequalities}

We need several Carleson inequalities, in various settings. The first ones that we will need are

\begin{theorem}[Carleson inequality] \label{ThCarleson0}
Let $\Omega \subset \R^n$ be a 1-CAD. There exists $C>0$ such that for any  $f \in CM_1^\infty$ and $g\in L^1_{loc}(\Omega)$, we have
\begin{equation} \label{ThCarleson0a}
 \int_\Omega |f g| \, \frac{dX}{\delta} \leq C \|f\|_{CM_1^\infty} \|\N^1(g)\|_{L^1(\partial \Omega,\sigma)}.
 \end{equation}
Equivalently, there exists $C>0$ such that for any  $f \in CM_2^\infty$ and $g\in L^2_{loc}(\Omega)$, we have
\begin{equation} \label{ThCarleson0b}
 \|\A(fg)\|_{L^2(\partial \Omega,\sigma)} \leq C \|f\|_{CM_2^\infty} \|\N^2(g)\|_{L^2(\partial \Omega,\sigma)}.
  \end{equation}
\end{theorem}

\begin{theorem}[Carleson inequalities involving elliptic measures] \label{ThCarleson1}
Let $\Omega \subset \R^n$ be a 1-CAD and $L:=-\diver A \nabla$ be a uniformly elliptic operator on $\Omega$. Write $G_{L^*}$ and $\omega_L$ for respectively the Green function for $L^*$ and the elliptic measure for $L$ with pole at center/infinity as defined in Definition \ref{defomG}, and assume 
\[\omega_L \in A_\infty(\sigma).\]
Then there exists $C>0$ such that for any  $f \in CM_1^\infty$ and $g\in L^1_{loc}(\Omega)$, we have
\begin{equation} \label{ThCarleson1a}
 \int_\Omega |f g| \, G_{L^*}\, \frac{dX}{\delta^2} \leq C \|f\|_{CM_1^\infty} \|\N^1(g)\|_{L^1(\partial \Omega,\omega_L)}.
 \end{equation}
Equivalently, there exists $C>0$ such that for any  $f \in CM_2^\infty$ and $g\in L^2_{loc}(\Omega)$, we have
\begin{equation} \label{ThCarleson1b}
 \|\A(fg)\|_{L^2(\partial \Omega,\omega_L)} \leq C \|f\|_{CM_2^\infty} \|\N^2(g)\|_{L^2(\partial \Omega,\omega_L)}.
  \end{equation}
\end{theorem}

\noindent{\em Proof of Theorems \ref{ThCarleson0} and \ref{ThCarleson1}.}
Let us first note that \eqref{ThCarleson0b} and \eqref{ThCarleson1b} are immediate consequences of Lemma \ref{lemFubini} and the $L^2$ version of \eqref{ThCarleson0a} and \eqref{ThCarleson1a} respectively. The reverse implication \eqref{ThCarleson1b} $\implies$  \eqref{ThCarleson1a} are \eqref{ThCarleson0b} $\implies$  \eqref{ThCarleson0a} similarly simple by taking square roots.

Let us prove \eqref{ThCarleson0a} and \eqref{ThCarleson1a}. The classical Carleson inequality in $\R^n_+$ with respect to the Lebesgue measure can be found in \cite[II.2.2, Theorem 2]{Stein93}, and the arguments extends to the case where the measure $\mu$ at the boundary is doubling. More precisely, the Carleson inequality that we obtain that way is
\begin{multline} \label{CarlesonP1} \int_{\partial \Omega} \int_{\Gamma(z)} |uv(X)| \frac{dX}{\delta(X)^n} d\mu(z) \\ 
\leq C_\mu \left( \sup_{x\in \partial \Omega, r>0} \int_{\partial \Omega} \int_{\Gamma(z) \cap B(x,r)} |u(X)|  \frac{dX}{\delta(X)^n} d\mu(z)  \right) \left( \int_{\partial \Omega} \N(v) d\mu \right).\end{multline}

\medskip

For \eqref{ThCarleson0a}, we apply the Carleson inequality \eqref{CarlesonP1} with $\mu = \sigma$, 
\[ u(X) := \sup_{B_X} |f| \quad \text{ and } \quad v(X) := \fint_{B_X} |g(Y)|dY,\]
and we obtain
\begin{multline*}
\int_\Omega |f g| \, \frac{dX}{\delta} \approx \int_{\partial \Omega} \A_1^1(fg) \, d\sigma(z)   \leq \int_{\partial \Omega} \int_{\Gamma(z)} |uv(X)| \frac{dX}{\delta(X)^n} d\sigma(z) \\
\lesssim \Big( \sup_{x\in \partial \Omega, r>0} \underbrace{\int_{\partial \Omega} \int_{\Gamma(z)} |u(X)| \1_{B(x,r)}  \frac{dX}{\delta(X)^n} d\sigma(z)}_{ \approx \int_{\Omega \cap B(x,r)} |u| \frac{dX}{\delta} \text{ by Lemma \ref{lemFubini}}}  \Big)  \|\N^1(g)\|_{L^1(\partial \Omega,\omega_L)} \\ \approx \|f\|_{CM_1^\infty} \|\N^1(g)\|_{L^1(\partial \Omega,\omega_L)}.
\end{multline*}

For \eqref{ThCarleson1a}, we apply the Carleson inequality \eqref{CarlesonP1} with $\mu = \omega_L$, and the same $u,v$ as above.
Combined with (ii) of Lemma \ref{lemFubini},
\begin{multline*}
\int_\Omega |f g| \, \frac{G_{L^*}}{\delta^2} dX \approx \int_{\partial \Omega} \A_1^1(fg) \, d\omega_L(z)   \leq \int_{\partial \Omega} \int_{\Gamma(z)} |uv(X)| \frac{dX}{\delta(X)^n} d\omega_L(z) \\
\lesssim \Big( \sup_{x\in \partial \Omega, r>0} \underbrace{\int_{\partial \Omega} \int_{\Gamma(z)} |u(X)| \1_{B(x,r)}  \frac{dX}{\delta(X)^n} d\omega_L(z)}_{ \approx \int_{\Omega \cap B(x,r)} |u| \frac{G_{L^*}}{\delta^2} dX \text{ by Lemma \ref{lemFubini}}}  \Big)  \|\N^1(g)\|_{L^1(\partial \Omega,\omega_L)}.
\end{multline*}
However, \cite[Theorem 1.29]{FP23} --- or more precisely \cite[Lemma 3.30]{FP23} for the equivalence of norms --- shows that if $\omega \in A_\infty(\sigma)$, then the Carleson norms of $f$ with respect to $\omega_L$ or $\sigma$ are equivalent, i.e.  
\[ \sup_{x\in \partial \Omega, r>0} \int_{\Omega \cap B(x,r)} |u| \frac{G_{L^*}}{\delta^2} dX \approx \|f\|_{CM_1^\infty}.\]
The theorems follow.
\ep

The next Carleson inequality that we require is:

\begin{theorem} \label{ThCarleson2}
Let $\Omega \subset \R^n$ be a 1-CAD. Then the following holds. For any $p\in (1,\infty)$, there exists $C_p$ such that, for any $f\in CM_2^\infty$ and any $g\in L^2_{loc}(\Omega)$, we have 
\[\| \A_1^2(fg)\|_{L^p(\partial \Omega,\sigma)} \leq C_p \|f\|_{CM_2^\infty} \|\A(g)\|_{L^p(\partial \Omega,\sigma)}.\]
\end{theorem}

\bp Our base Carleson inequality this time is
\begin{equation} \label{CarlesonP2}
 \int_{\partial \Omega} \int_{\Gamma(x)} |uv(X)| \frac{dX}{\delta(X)^n} \, d\sigma(x) \leq C \int_\Omega \bar {\mathcal C}(u)(x) \bar \A(v)(x) \, d\sigma(x),
 \end{equation}
where 
\[ \bar {\mathcal C}(u)(x) := \sup_{r>0} \left( r^{1-n} \int_{B(x,r)\cap \Omega} |u|^2  \frac{dX}{\delta}\right)^\frac12 \quad \text{ and } \quad  \bar {\mathcal A}(u)(x) :=  \left( \int_{\Gamma(x)} |u|^2  \frac{dX}{\delta^{n}}\right)^\frac12. \]
The inequality \eqref{CarlesonP2} is proved in the flat case, i.e. $\Omega = \R^n_+$, in \cite[Theorem 1]{CMS85}. However, the proof does not rely on the Euclidean structure: it uses a stopping time argument and can be easily extended to our scenario where the boundary is Ahlfors regular.

\medskip

Take now $f\in CM_2^\infty$ and $g,h\in L^2_{loc}(\Omega)$. Similarly as the previous lemma, we apply \eqref{CarlesonP2} with
\[ u(X) := \sup_{B_X} |f| \quad \text{ and } \quad v(X) := \fint_{B_X} |g(Y)h(Y)|dY.\]
Together with (i) of Lemma \ref{lemFubini}
\begin{multline} \label{CarlesonP3}
\int_\Omega |fgh| \, \frac{dX}{\delta} \approx \int_{\partial \Omega} \A^1_1(fgh) \, d\sigma \leq \int_\Omega \int_{\Gamma(z)} |uv(X)| \frac{dX}{\delta(X)^n} \, d\sigma(z) \\
\lesssim \int_\Omega \bar{\mathcal C}(u)(x) \bar \A(v)(x) \, d\sigma(x) \leq \underbrace{\big\| \bar{\mathcal C}(u)(x)\big\|_{L^\infty(\partial \Omega,\sigma)} }_{=\|f\|_{CM_2^\infty}} \int_{\partial \Omega} \A_2^1(gh) \, d\sigma(x) \\
\leq \|f\|_{CM_2^\infty} \int_{\partial \Omega} \A(g)\N^2(h) \, d\sigma(x) \leq \|f\|_{CM_2^\infty} \|\A(g)\|_{L^p(\partial \Omega,\sigma)} \|\N^2(h)\|_{L^{p'}(\partial \Omega,\sigma)}
\end{multline}
by using the H\"older inequality twice.

\medskip

We can assume that $g$ is compactly supported. Indeed, assuming that the case where $g$ is compactly supported is true, we can construct a sequence of compact sets $K_n \uparrow \Omega$, and then for a general $g$ we have
\[\| \A_1^2(fg\1_{K_n})\|_{L^p(\partial \Omega,\sigma)} \leq C_p \|f\|_{CM_2^\infty} \|\N^2(g\1_{K_n})\|_{L^p(\partial \Omega,\sigma)}.\]
Taking the limit as $n\to \infty$ gives the result for non-compactly supported $g$.

\medskip

We assume now that $g$ is compactly supported, which means that $ \|\A_1^2(fg)\|_{p} < \infty$. The result is now a simple consequence of the duality of tent spaces and the Carleson inequality. The duality in the case of the functional $\A_1^2$ is treated in \cite[lemma A.8]{Mando26} (see also \cite[Proposition 2.4]{MPT} for the duality the Carleson functional instead, whose $L^p$ norm is equivalent to the one of the square function $\A_1^2$)  and gives
\[ \| \A_1^2(fg)\|_{L^p(\partial \Omega,\sigma)}  \lesssim \sup_{h:\, \|\N^2(h)\|_{p'} \leq 1} \int_\Omega fgh \, \frac{dX}{\delta}.\]
The inequality \eqref{CarlesonP3} implies now
\begin{multline*} 
\| \A_1^2(fg)\|_{L^p(\partial \Omega,\sigma)} \lesssim  \sup_{h:\, \|\N^2(h)\|_{p'} \leq 1}   \|f\|_{CM_2^\infty} \|\A(g)\|_{L^p(\partial \Omega,\sigma)} \|\N^2(h)\|_{L^{p'}(\partial \Omega,\sigma)} \\
= \|f\|_{CM_2^\infty} \|\A(g)\|_{L^p(\partial \Omega,\sigma)}.
\end{multline*}
The theorem follows.
\ep

The last needed Carleson inequality pertains to cut-off functions, which we need in our proofs. We start by introducing the notion and properties of cut-off functions.

\begin{definition} \label{defcutoff}
We say that a function $\varphi \in C^\infty(\Omega)$ is a cut-off function if there exists $\gamma>1$, $\kappa \in (0,4)$, and $K>0$ such that
\begin{enumerate}[(i)]
\item $0\leq \varphi \leq 1$ on $\Omega$,
\item $\delta |\nabla \varphi| \leq K$,
\item $\A_{1,(\gamma,\kappa)}^\infty(\1_{\supp \nabla \varphi})(x) \leq K$.
\end{enumerate}
We say that $\varphi$ is a  $(\gamma,\kappa,K)$-cut-off function when we want to refer to the constants, and we call the triple $\gamma$, $\kappa$, and $K$ the cut-off constants of $\varphi$.
\end{definition}

\begin{lemma} \label{lemcutoff}
Let $\Omega$ be a 1-CAD. There exists $\gamma>1$, $\kappa\in (0,4)$, and $K>0$ such that, for any $\zeta \in C^\infty(\R,[0,1])$ satisfying $\zeta(t) = 0$ when $t\leq 1$ and $\zeta(t) =1$ when $t\geq 2$, and $|\zeta'(t)| \leq 2$, any $E\subset \partial \Omega$, and any $l>0$, the functions $\phi_E$ and $\phi^l$ defined as
\[ \phi_E(X):= \zeta\left( \frac{4\delta(X)}{\dist(X,E)} \right) \ \text{ and } \  \phi^l(X):= \zeta\left( \frac{\delta(X)}{l} \right) \]
are $(\gamma,\kappa,K)$-cut-off functions.
\end{lemma}

\bp Choose $\gamma = \frac54$ and $\kappa = \frac14$. We want to show that there exists $K>0$ such that $\delta |\nabla \phi| \leq K$ and $\A_{1,(\frac54,\frac14)}^\infty(\1_{\supp \varphi})(x) \leq K$ when $\phi = \phi_E$ or $\phi = \phi^l$. The proof is computational but straightforward, and so we will only show the estimates for $\phi_E$, and leave $\phi^l$ as an exercise for the reader.

First, note that 
\[ \delta(X) \nabla \phi_E(X) = \frac{4\delta(X)}{\dist(X,E)}\zeta'\left( \frac{4\delta(X)}{\dist(X,E)} \right)  \left[ \nabla \delta(X) - \frac{\delta(X)}{\dist(X,E)}\nabla \dist(X,E) \right].\]
But since $1 \leq \frac{4\delta(X)}{\dist(X,E)} \leq 2$ on $\supp \zeta'$ and $|\zeta'|\leq 2$,
\[|  \delta(X) \nabla \phi_E(X) | \leq 16 [ 1 +  \frac12 ] \leq 24.\]

We turn to the second bound. For any $X\in \Gamma_{(5/4)}(x)$ and $Y\in \frac14 B_X$, we have
\begin{equation} \label{randombdA}
\dist(Y,E) \leq |Y-x| \leq |X-x| + \frac1{16} \delta(X) \leq \frac{21}{16}\delta(X) \leq \frac{21}{15}\delta(Y).
\end{equation}
When $x\in \overline{E}$, we have then $4\delta(Y)/\dist(Y,E) \geq \frac{60}{21} \geq 2$ hence $\nabla \phi_E(Y) = 0$, and consequently 
$\A_{1,(\frac54,\frac14)}^\infty(\1_{\supp \varphi})(x) = 0$.

Take now $x\in \partial \Omega \setminus \overline{E}$, $X\in \Gamma_{(5/4)}(x)$ and $Y\in \frac14 B_X$. With \eqref{randombdA}, it forces $Y \in \Gamma_{(21/15)}(x)$, and if in addition $\nabla \phi_E(Y) \neq 0$, then 
\[|Y-x| \leq \frac{21}{15} \delta(Y) \leq \frac{21}{30}\dist(Y,E) \leq \frac{21}{30} [|Y-x| + \dist(x,E)]\]
and
\[ \dist(x,E) \leq \dist(Y,E) + |Y-x| \leq \frac{51}{30} \dist(Y,E) \leq \frac{204}{30} \delta(Y) .\]
We conclude that $\frac{30}{204} \dist(x,E) \leq \delta(Y) \leq \frac{30}{9}  \dist(x,E)$, hence that
\[ \frac18 \dist(x,E) \leq \frac{16 \times 30}{17 \times 204} \dist(x,E) \leq  \delta(X) \leq \frac{16 \times 30}{15 \times 9} \dist(x,E)  \leq 5 \dist(x,E).\]
The computations show that 
\[\A_{1,(\frac54,\frac14)}^\infty(\1_{\supp \varphi})(x)  \lesssim \frac1{\dist(x,E)^n} |\{X\in \Gamma(x), \, \frac18 \dist(x,E) \leq \delta(X) \leq 5\dist(x,E)\}| \lesssim 1\]
as desired for the lemma.
\ep

\begin{corollary} \label{corcutoff}
Let $\Omega$ be a 1-CAD. Then there exists $\gamma>1$, $\kappa \in (0,4)$ and $K>0$ such that, for any $0< \epsilon < l$ and any $E \subset \partial \Omega$, the functions $\phi_E(1-\phi^l)$ and $\phi_E\phi^\epsilon(1-\phi^l)$ are $(\gamma,\kappa,K)$-cut-off functions.

Moreover, if $x_0\in \partial \Omega$ and $\Delta(x_0,r) = B(x_0,r) \cap \partial \Omega$, the functions 
\[ \phi_{\Delta(x_0,r)} \phi^{1/r}(1-\phi^{r})\]
are $(\gamma,\kappa,K)$-cut-off functions in $C^\infty_0(\Omega)$ that pointwise converges to $\1_\Omega$.
\end{corollary}

\bp
The first part is a consequence of the facts that
\begin{itemize}
\item $(1-\phi)$ is a $(\gamma,\kappa,K)$-cut-off function whenever $\phi$ is.
\item $\phi\psi$ is a $(\gamma,\kappa,2K)$-cut-off function whenever $\phi$ and $\psi$ are $(\gamma,\kappa,K)$-cut-off functions.
\end{itemize}
The second part comes from the fact that for any $X\in \Omega$, we take $x \in \partial \Omega$ such that $|X-x| = \delta(X)$ and then when $r$ is large enough, we have  $\frac1r \leq |X-x| \leq r/2$ and $x\in \Delta(x_0,r)$.
\ep

We can now state the Carleson inequality involving cut-off functions.
\begin{theorem} \label{ThCarlesonCO}
Let $\Omega$ be a 1-sided CAD, and let $\mu$ be a doubling measure on $\partial \Omega$. Let $p\in [1,\infty)$, and $\alpha, \beta \in [1,\infty]$. For any $\gamma>1$, $\kappa \in (0,4)$ and $K>0$, there exists $C = C_{\gamma,\kappa,K}>0$ such that, for any $f\in L^\beta_{loc}(\Omega)$ and any $(\gamma,\kappa,K)$-cut-off function $\varphi$
\[ \| \A_\alpha^\beta(f\1_{\supp \nabla \varphi})\|_{L^p(\partial \Omega,\mu)} \leq C \|\N^\beta(f)\|_{L^p(\partial \Omega,\mu)}.\]

In particular, if $\sigma$ is an Ahlfors regular measure on $\partial \Omega$, then
\[ \int_\Omega |f| \1_{\supp \nabla \varphi} \frac{dX}{\delta} \leq C \|\N^1(f)\|_{L^1(\partial \Omega,\sigma)} ,\]
and if $L:=-\diver A \nabla$ is a uniformly elliptic operator on $\Omega$, and $(G_{L^*},\omega_L)$ are the Green function for $L^*$ and the elliptic measure for $L$ with pole at center/infinity, then
\[ \int_\Omega |f| \1_{\supp \nabla \varphi} \frac{G_{L^*}}{\delta^2}\, dX \leq C \|\N^1(f)\|_{L^1(\partial \Omega,\omega_L)},\]
\end{theorem}

\bp
The first inequality is a consequence of the fact that 
\begin{multline*} 
\| \A_\alpha^\beta(f\1_{\supp \nabla \varphi})\|_{L^p(\partial \Omega,\mu)} \lesssim \| \A_{1,(\gamma,\kappa)}^\beta(f\1_{\supp \nabla \varphi})\|_{L^p(\partial \Omega,\mu)} \leq \| \A_{1,(\gamma,\kappa)}^\infty(\1_{\supp \nabla \varphi}) \N^\beta_{(\gamma,\kappa)}(f) \|_{L^p(\partial \Omega,\mu)} \\ \leq K  \|\N^\beta_{(\gamma,\kappa)}(f) \|_{L^p(\partial \Omega,\mu)} \approx \|\N^\beta(f)\|_{L^p(\partial \Omega,\mu)},
\end{multline*}
where we used in order Lemma \ref{lemcompNA}, assumption (iii) of Definition \ref{defcutoff}, and Lemma \ref{lemcompNA} again.

The two last inequalities are consequences of the first one and Lemma \ref{lemFubini}.
\ep

\subsection{On the Green function}

We start with the following lemma providing local estimates on  the solutions, and those proof can be found in \cite[Lemma 4.47]{HMT17}.

\begin{lemma} \label{lembdonG}
Let $\Omega \subset \R^n$, and $L=-\diver A \nabla$ be a uniformly elliptic operator that satisfies $\delta |\nabla A| \leq K$. Then
\begin{enumerate}[(i)]
\item for any ball $B$ of radius $r_B$ that satisfies $2B \subset \Omega$, and for any solution $u$ to $Lu=0$, we have
 \begin{equation} \label{uW22}
 r_B \left( \fint_{B} |\nabla^2 u|^2 \, dX\right)^\frac12 \leq C_K \left( \fint_{2B} |\nabla u|^2\right)^\frac12
 \end{equation}
 \item If in addition, $u\geq 0$ on $2B$, then, for any $X,Y\in B$,
  \begin{equation} \label{uW1infty}
  \delta(X)|\nabla u(X)| \leq C_K u(Y)
   \end{equation}
\end{enumerate}
\end{lemma}

Assuming that $\Omega$ is a Chord Arc Domain, we have additional ``boundary'' control on the second derivative of the Green function.

\begin{theorem} \label{ThonG}
Let $\Omega \subset \R^n$ be a Chord Arc Domain. Let $G$ be the Green function (for $-\Delta$) with pole at infinity/center. Then there exists $K$ depending only on the CAD constants such that
\[ \frac{\delta^2 \nabla^2 G}{G} \1_{\Omega \setminus 2B_\Omega} \in CM_2^\infty(K),\]
where $B_\Omega = \emptyset$ when $\Omega$ is unbounded and $B_\Omega $ is the center ball  from Definition \ref{defCP} when $\Omega$ is bounded.
\end{theorem}

\bp The theorem is a (weaker) variant of \cite[Theorem 1.13]{FL23} (see also \cite[Theorem VI]{Azz??}) that states the existence of $K'$ such that,
\begin{equation} \label{CMbd1}
r^{1-n} \int_{\Omega \cap B(x,r)} \frac{|\nabla^2 G|^2}{G^2}\delta^3 \, dX  \leq K'
\end{equation}
whenever $B(x,r)$ is ball centered on $\partial \Omega$ that satisfies $B(x,2r) \cap B_\Omega = \emptyset$.

Take now $x\in \partial \Omega$ and $r>0$ such that  $B(x,4r) \cap B_\Omega = \emptyset$, and then $X\in B(x,r)$. One can check that $2B_X \subset B(x,2r) \cap \Omega \subset \Omega \setminus B_\Omega$, so we can use the Moser estimate on the harmonic functions $\partial_i\partial_j G$, $i,j\in \{1,\dots,n\}$, to obtain that
\[ \sup_{B_X} |\nabla^2 G| \lesssim  \left( \fint_{2B_X} |\nabla G|^2 \right)^\frac12.\]
Since $G$ is positive on $3B_X$, one can also use the Harnack inequality to further write
\[ \left( \sup_{B_X} \frac{|\nabla^2 G|}{G} \right)^2 \lesssim  \fint_{2B_X} \frac{|\nabla^2 G|^2}{G^2} dY .\]
It yields, when $x\in \partial \Omega$ and $r < \dist(x,B_\Omega)/4$, that
\begin{multline} \label{CMbd2}
I(x,r):= r^{1-n} \int_{X\in \Omega \cap B(x,r)} \left(\sup_{Y\in B_X} \frac{|\nabla^2 G(Y)|}{G(Y)} \1_{\Omega \setminus 2B_\Omega}(Y)\right)^2 \delta(X)^3 \, dX \\
\lesssim  r^{1-n} \int_{X\in \Omega \cap B(x,r)} \left(\fint_{Y\in 2B_X} \frac{|\nabla^2 G(Y)|^2}{G(Y)^2} dY\right)\,  \delta(X)^3 \, dX \\
\lesssim (2r)^{1-n} \int_{\Omega \cap B(x,2r)} \frac{|\nabla^2 G|^2}{G^2} \,  \delta^3 \, dY \leq K',
\end{multline}
where we used Fubini theorem for the second inequality, and then \eqref{CMbd1}. 

To prove the theorem, it then only remains to prove a bound the quantity $I(x,r)$ when $B(x,4r) \cap B_\Omega \neq \emptyset$. But if that is the case, it means that $\Omega$ is bounded and $r\approx \diam \Omega$, meaning that 
\[ I(x,r) \lesssim I_\Omega := \frac1{\sigma(\partial \Omega)} \int_{X\in \Omega} \left(\sup_{Y\in B_X} \frac{|\nabla^2 G(Y)|}{G(Y)} \1_{\Omega \setminus 2B_\Omega}(Y)\right)^2 \delta(X)^3 \, dX.\]
To bound $I_\Omega$, we cover $\partial \Omega$ by the balls $\{B_i:= B(x,r_x)\}_{x\in \partial \Omega}$, where $r_x = \dist(x,B_\Omega)/40$. By Vitali's lemma, we can find a non-overlapping subcollection $\{B(x_i,r_{x_i})\}_{i\in I}$ such that $\{5B_i\}_{i\in I}$ covers $\partial \Omega$, meaning that $\{10B_i\}_{i\in I}$ covers $\Omega_b:= \{X\in \Omega, \delta(X) \leq \inf_{i\in I} 5r_{x_i}\}$. It yields that 
\begin{multline*}
\sigma(\partial \Omega ) I_\Omega \leq \int_{X\in \Omega \setminus \Omega_b} \left(\sup_{Y\in B_X} \frac{|\nabla^2 G(Y)|}{G(Y)} \1_{\Omega \setminus 2B_\Omega}(Y)\right)^2 \delta(X)^3 \, dX \\ + \sum_{i\in I} \int_{X\in 10B_i} \left(\sup_{Y\in B_X} \frac{|\nabla^2 G(Y)|}{G(Y)} \1_{\Omega \setminus 2B_\Omega}(Y)\right)^2 \delta(X)^3 \, dX =: T_1 + T_2.
\end{multline*}
The term $T_2$ is simple. Since $10r_{x_i} \leq \dist(x_i,B_\Omega)/4$, we have by \eqref{CMbd2} that
\[ T_2 = \sum_{i\in I} (10r_{x_i})^{1-n} I(x_i,10r_{x_i}) \lesssim K' \sum_{i\in I} \sigma(B_i) \leq K' \sigma(\Omega).\]
As for the term $T_1$, we use the fact that $G$ is positive and all its derivatives are harmonic on $\Omega \setminus B_\Omega$ to get that $\delta^2|\nabla^2 G| \lesssim G$ on $\Omega \setminus 2B_\Omega$. It entails that 
\[ T_1 \lesssim \int_{\Omega \setminus \Omega_b} \frac{dX}{\delta(X)}.\]
But $r_{x_i} \approx \dist(x_i,B_\Omega) \approx \diam(\Omega)$, because $B_\Omega$ is the center ball of $\Omega$, meaning that $\delta(X) \approx \diam(\Omega)$ on $\Omega \setminus \Omega_b$. It implies that 
\[ T_1 \lesssim \frac{|\Omega|}{\diam(\Omega)} \approx \diam(\Omega)^{1-n} \approx \sigma(\partial \Omega)\]
by the Ahlfors regularity of the boundary. 

We conclude that $\sigma(\partial \Omega) I_\Omega \lesssim (1+K') \sigma(\partial \Omega)$, meaning that $I_\Omega \lesssim (1+K') < +\infty$, as desired for the theorem.
\ep 

\begin{theorem} \label{CorNG>0}
Under the same assumption as Theorem \ref{ThonG}, there exist $c>0$ and $K>0$ depending only on the CAD constants such that
\[ \1_{\{\delta |\nabla G| \leq c G\}} \in CM_2^\infty(K).\]
\end{theorem}

\bp
As in Theorem \eqref{ThonG}, the key computation is done previously in \cite[Lemma 5.8]{FL23}, and states the existence of  $c>0$ and $K'$ such that  
\begin{equation} \label{CMbd3}
r^{1-n} \int_{\Omega \cap B(x,r)}\1_{\{\delta |\nabla G| \leq 2c G\}} \, dX  \leq K'
\end{equation}
whenever $x\in \partial \Omega$ and $r\leq \dist(x,B_\Omega)/2$. So the main difficulty here is to replace $\1_{\{\delta |\nabla G| \leq c G\}}(X)$ by $\osc_{B_X} \1_{\{\delta |\nabla G| \leq c G\}}$. This is done, as one can expect, by using the local H\"older regularity of $G$ and $\nabla G$. Indeed, we have
\[
\osc_{\epsilon B_X} \left| \frac{\delta\nabla G}{G} \right| \leq \frac{\epsilon\delta(X)}{4} \dfrac{\ds \max_{\epsilon B_X}|\nabla G|}{\ds \min_{\epsilon B_X} G} 
+ \big(1+\frac{\epsilon}4\big) \delta(X) \dfrac{\ds \osc_{\epsilon B_X} |\nabla G|}{\ds \min_{\epsilon B_X} G} 
+ \big(1+\frac{\epsilon}4\big) \delta(X) \dfrac{\ds (\max_{\epsilon B_X} |\nabla G|) \osc_{\epsilon B_X} G}{\ds (\min_{\epsilon B_X}G)^2}
\]
As long as $2B_X \cap B_\Omega = \emptyset$ and $\epsilon \in (0,\frac12)$, we have for a $\alpha \in (0,1)$
\begin{enumerate}
\item with the Haranck inequality on $G\geq 0$, 
\[\min_{B_X}G \approx \max_{B_X} G \approx G(X);\] 
\item with the Moser estimate applied to $\nabla G$, the Cacciopoli inequality applied to $G$, and (1)
\[ \max_{\epsilon B_X} |\nabla G| \lesssim \left( \fint_{\frac34 B_X} |\nabla G|^2 \right)^\frac12 \lesssim \frac1{\delta(X)}  \left( \fint_{B_X} |G|^2 \right)^\frac12 \lesssim \frac{G(X)}{\delta(X)};\]
\item with the H\"older continuity and the Moser estimate applied to the solution $G$, and (1),
\[ \osc_{\epsilon B_X} \lesssim \epsilon^\alpha  \left( \fint_{B_X} |G|^2 \right)^\frac12 \lesssim \epsilon^\alpha G(X);\]
\item and similarly, with the  H\"older continuity applied to the solutions $\partial_i G$
\[ \min_{\epsilon B_X}|\nabla G| \lesssim \frac{G(X)}{\delta(X)}.\]
\end{enumerate}
Consequently, the oscillation estimate on $\frac{\delta\nabla G}{G}$ becomes
\[ \osc_{\epsilon B_X} \left| \frac{\delta\nabla G}{G} \right|  \leq C \epsilon^\alpha.\]
We fix $\epsilon$ such that $C\epsilon^\alpha \leq c$, where $c>0$ is the constant in \eqref{CMbd3}, and now we have
\[ \Big\{X\in \Omega, \, 2B_X \cap B_\Omega = \emptyset, \, \inf_{\epsilon B_X}  \left| \frac{\delta\nabla G}{G} \right| < c \Big\} \subset  \Big\{X\in \Omega, \, 2B_X \cap B_\Omega = \emptyset, \,  \left| \frac{\delta(X)\nabla G(X)}{G(X)} \right| < 2c \Big\},\]
hence, together with \eqref{CMbd3}, 
\begin{equation} \label{CMbd3b}
r^{1-n} \int_{\Omega \cap B(x,r)} \left(\sup_{\epsilon B_X} \1_{\{\delta |\nabla G| \leq c G\}}\right) \, dX  \leq K'
\end{equation}
whenever $x\in \partial \Omega$ and $r\leq \dist(x,B_\Omega)/2$.

\medskip

Let $x\in \partial \Omega$ and $r\in (0,\diam\Omega)$. If $r < \dist(x,B_\Omega)/4$, then we use the brutal bound
\[ \sup_{B_X}\1_{\{\delta |\nabla G| \leq c G\}} \lesssim (\epsilon r)^{-n} \int_{B_X} \Big(\sup_{\epsilon B_Y} \1_{\{\delta |\nabla G| \leq c G\}}\Big) dY\]
to get with Fubini's lemma and \eqref{CMbd3b}
\begin{multline} \label{CMbd4}
I(x,r):= r^{1-n} \int_{\Omega \cap B(x,r)} \left(\sup_{B_X} \1_{\{\delta |\nabla G| \leq c G\}}\right) \, dX \\ \lesssim r^{1-n} \int_{\Omega \cap B(x,r)} \left(\frac1{(\epsilon r)^{n}} \int_{B_X} \Big(\sup_{\epsilon B_Y} \1_{\{\delta |\nabla G| \leq c G\}}\Big) dY \right) \, dX \\
\leq C_\epsilon r^{1-n} \int_{\Omega \cap B(x,2r)} \left(\sup_{B_Y} \1_{\{\delta |\nabla G| \leq c G\}}\right) \, dY \lesssim K'.
\end{multline}
If instead $r \geq \dist(x,B_\Omega)/4$, it means that $\Omega$ is bounded and $r\approx \diam (\Omega)$, so 
\[ I(x,r) \lesssim I_\Omega := \frac1{\sigma(\partial \Omega)} \int_{\Omega} \left(\sup_{B_X} \1_{\{\delta |\nabla G| \leq c G\}}\right) \, dX.\]
From there, the computations are similar as the ones found in the proof of Theorem \ref{ThonG}. We cover $\partial \Omega$ by the balls $\{B_i:= B(x,r_x)\}_{x\in \partial \Omega}$, where $r_x = \dist(x,B_\Omega)/40$. By Vitali's lemma, we can find a non-overlapping subcollection $\{B(x_i,r_{x_i})\}_{i\in I}$ such that $\{5B_i\}_{i\in I}$ covers $\partial \Omega$, meaning that $\{10B_i\}_{i\in I}$ covers $\Omega_b:= \{X\in \Omega, \delta(X) \leq \inf_{i\in I} 5r_{x_i}\}$. It yields that 
\[\sigma(\partial \Omega ) I_\Omega \leq \int_{X\in \Omega \setminus \Omega_b}\, \frac{dX}{\delta(X)} \\ + \sum_{i\in I} \int_{X\in 10B_i} \left(\sup_{B_X} \1_{\{\delta |\nabla G| \leq c G\}}\right) \, \frac{dX}{\delta(X)} =: T_1 + T_2.\]
The term $T_2$ is simple. Since $10r_{x_i} \leq \dist(x_i,B_\Omega)/4$, we have by \eqref{CMbd4} that
\[ T_2 = \sum_{i\in I} (10r_{x_i})^{1-n} I(x_i,10r_{x_i}) \lesssim K' \sum_{i\in I} \sigma(B_i) \leq K' \sigma(\Omega).\]
As for the term $T_1$, we use the fact that $r_{x_i} \approx \dist(x_i,B_\Omega) \approx \diam(\Omega)$, meaning that $\delta(X) \approx \diam(\Omega)$ on $\Omega \setminus \Omega_b$ and consequently 
\[ T_2 \lesssim \frac{|\Omega|}{\diam(\Omega)} \approx \diam(\Omega)^{1-n} \approx \sigma(\partial \Omega)\]
by the Ahlfors regularity of the boundary. We conclude that $I_\Omega \lesssim (1+K')$, which ends the proof of the theorem.\ep

\subsection{On the Neumann function}

In this subsection, we aim to give the definition and few basic properties of the Neumann function, which is can be used to represent the weak solution to $Lu=h-\diver \vec F$ with zero Neumann data and zero average. Let us start by recalling the correct spaces for the use of the Lax-Milgram theorem in the Poisson-Neumann problem.

\begin{theorem} \label{ThLM}
Let $\Omega \subset \R^n$ be a bounded 1-CAD, and $L=-\diver A \nabla$ be a uniformly elliptic operator. Let $\widehat{W}^{1,2}(\Omega)$ be the subspace of $W^{1,2}(\Omega)$ of functions with zero average, that is,
\begin{equation} \label{defhatW12}
\widehat W^{1,2}(\om):=\set{f\in W^{1,2}(\om), \, \int_{\Omega} f\, dX=0}.
\end{equation}
Then $\widehat W^{1,2}(\Omega)$ is a Hilbert space equipped with the inner product $\left<u, v\right> = \iint_\Omega \nabla u \cdot \nabla v\, dX$. Moreover, by the Lax-Milgram theorem, for any $h\in L^2(\Omega)$, $\vec F \in L^2(\Omega,\R^n)$, there exists a unique $u_{\vec F,h} \in \widehat{W}^{1,2}(\Omega)$ such that 
    \begin{equation} \label{LMNeumann} 
    \int_{\Omega} A \nabla u_{\vec F,h} \cdot \nabla \varphi \, dX = \int_\Omega \Big( h \varphi  + \vec F \cdot \nabla \varphi\Big) \, dX \qquad \text{ for any } \varphi \in  \widehat{W}^{1,2}(\Omega).
    \end{equation}
    If $\int_\Omega h\, dX = 0$, such $u_{\vec F,h}$ satisfies
    \begin{equation} \label{LMNeumann2} 
    \int_{\Omega} A \nabla u_{\vec F,h} \cdot \nabla \varphi \, dX = \int_\Omega \Big(h \varphi  + \vec F \cdot \nabla \varphi \Big) \, dX \qquad \text{ for any } \varphi \in  W^{1,2}(\Omega),
    \end{equation}
    in particular, $u_{\vec F,h} \in \widehat{W}^{1,2}(\Omega)$ is a weak solution to $Lu=h-\diver \vec F$ with zero Neumann data. 
    
Moreover, the map $\vec F \mapsto u_{\vec F,0}$ is bounded from $L^2(\Omega)$ to $W^{1,2}(\Omega)$.
\end{theorem}

\bp
A simple application of the Lax-Milgram theorem, with the quadratic form
\[ u,v\in \widehat{W}^{1,2}(\Omega) \mapsto \int_{\Omega} A \nabla u \cdot \nabla v \, dX\]
and the linear forms
\[ v \in \widehat{W}^{1,2}(\Omega) \mapsto \int_\Omega \Big( h v  + \vec F \cdot \nabla v\Big) \, dX.\]
Note that $v \in \widehat{W}^{1,2}(\Omega) \mapsto \int_\Omega h v \, dX$ is indeed a continuous form thanks to the Poincar\'e inequality 
\[ \int_{\Omega} v^2 \, dX \leq C\diam(\Omega) \left(\int_{\Omega} |\nabla v|^2\right)^\frac12\]
which holds for functions $v\in \widehat{W}^{1,2}(\Omega)$  with a constant $C$ that depends only on the 1-CAD constants of $\Omega$.
\ep

The next theorem introduces the Neumann function.

\begin{theorem} \label{ThNeumann}
Let $\Omega \subset \R^n$ be bounded 1-CAD domain, and $L=-\diver A \nabla$ be a uniformly elliptic operator. There exists a unique function $N$ --- called Neumann function --- defined on $\Omega \times \Omega$ such that $N(X,\cdot)$ is continuous on $\overline{\Omega}\setminus \{X\}$, locally integrable in $\Omega$, $\int_{\Omega}N(X,Y)dX =0$, and that the 
 following representation formula holds: for any $Y\in\overline{\om}$, 
any $\phi\in W^{1,2}(\om)$, 
\begin{equation}\label{eq.RieszNmn'}
     \int_{\om}A(X)\nabla_X N(X,Y)\cdot \nabla \phi(X)dX=\phi(Y) - \fint_{\Omega}\phi(Z)\, dZ.
\end{equation}
Moreover,
\begin{enumerate}[(i)]
    \item for any $X,Y\in \Omega$, $N^*(Y,X) = N(X,Y)$, where $N^*(X,Y)$ is the Neumann function for the adjoint operator  $L^*=-\diver(A^T\nabla)$;
    \item for any $\vec F \in L^2(\Omega)$ and any $h\in L^2(\Omega)$ with zero average, the function $ u_{\vec F,h}$ defined as
  \begin{equation}  \label{defuFh}
  u_{\vec F,h}(X) := \int_\Omega \Big[\nabla_Y N(X,Y) \cdot \vec F(Y) + N(X,Y)h(Y) \Big]\, dY
  \end{equation}
  is the unique weak solution $u_{\vec F,h} \in \widehat W^{1,2}(\Omega)$ to $L  u_{\vec F,h} = h - \diver \vec F$ with zero Neumann data given  by the Lax-Milgram theorem (Theorem \ref{ThLM}).
\end{enumerate}
\end{theorem}

\bp
When $\widehat W^{1,2}(\Omega)$ is the space of function in $W^{1,2}(\Omega)$ with zero average trace instead, the theorem is Theorem 2.29 and Proposition 2.34 in \cite{FL-PN}, which themselves comes from the study of the Neumann function in \cite{KP93}.

Let us recall the strategy of the proof. The main difficulty stems from the pole $Y$ around which $N(.,Y)$ is not in $W^{1,2}_{loc}(\Omega)$. One strategy would be to construct $N(.,Y)$ by approaching it using solutions to $L u = |B(Y,\epsilon)|^{-1}\1_{B(Y,\epsilon)}$: this argument was applied to construct the Green function $G(.,Y)$ in \cite{GW82,DFMhighcod,DFMmixed}. The strategy of Kenig and Pipher in \cite{KP93} is to instead study the properties of $v(X,Y) = N(X,Y) - G(X,Y)$ --- which is a solution inside $\Omega$ --- and to invoke the already known results on the Green function to deal with the pole. Let us take this route and give a sketch of the proof.

\medskip

\noindent {\bf Step 1:} First, an immediate property of the elliptic measure is the fact that for any $\varphi \in C^\infty_0(\R^n)$, the function $\varphi_{L^*}$ defined on $\Omega$ as
\[ \varphi_{L^*}(Y):= \int_{\partial \Omega} \varphi \, d\omega_{L^*}^Y\]
is the solution to $L^* u = 0$ with boundary data $\varphi$ on $\partial \Omega$. The function $\varphi - \varphi_{L^*}$ lies then in $W^{1,2}(\Omega)$ and has zero trace, meaning that 
\[ \int_\Omega A \nabla G(.,Y) \cdot \nabla \varphi \, dX = \int_\Omega A \nabla G(.,Y) \cdot \nabla [\varphi-\varphi_{L^*}] \, dX + \int_\Omega A \nabla G(.,Y) \cdot \nabla \varphi_{L^*} \, dX =:I_1 + I_2\]
The term $I_1$ is $\varphi(Y) - \varphi_{L^*}(Y)$ by definition of the Green function, and the term $I_2$ is 0, since $\varphi_{L^*}$ is a solution. We recover a classical formula:
\begin{equation} \label{AGphi}
 \int_\Omega A \nabla G(.,Y) \cdot \nabla \varphi \, dX = \varphi(Y) - \int_{\partial \Omega} \varphi \, d\omega_{L^*}^Y.
 \end{equation}

\medskip

For $Y\in \Omega$, the form 
\[\varphi \in \widehat W^{1,2}(\Omega) \mapsto \int_{\partial \Omega} \varphi \, d\omega_{L^*}^Y\]
is continuous, with norm $C\delta(Y)^{2-n}$. The proof is the same as \cite[Lemma 2.3]{KP93} and is left to the reader. By the Lax-Milgram theorem, it means that there is a $v_Y \in \widehat W^{1,2}(\Omega)$ such that 
\[\int_\Omega A\nabla v_Y \cdot \nabla \varphi \, dX =  \int_{\partial \Omega} \varphi \, d\omega_{L^*}^Y \quad \text{ for } \widehat W^{1,2}(\Omega),\]
that is
\begin{equation} \label{Avphi}
\int_\Omega A\nabla v_Y \cdot \nabla \varphi \, dX = - \fint_{\Omega} \varphi \, dX + \int_{\partial \Omega} \varphi \, d\omega_{L^*}^Y \quad \text{ for } \varphi \in W^{1,2}(\Omega).
\end{equation}
We define $v(X,Y) := v_Y(X) - \fint_\Omega G(Z,Y) \, dZ$ and then
\[ N(X,Y) := v(X,Y) + G(X,Y).\]
By construction, we have that $\int_\Omega N(Z,Y)\, dZ = 0$. From the definition of $N(X,Y)$, \eqref{AGphi} and \eqref{Avphi}, we get the desired identity \eqref{eq.RieszNmn'}. The continuity of $N(.,Y)$ outside of the pole comes from the fact that $N(.,Y)$ is a solution to $Lu=-1$ in $\Omega$ with zero Neumann data, meaning that the H\"older continuity of solutions with zero Neumann data applies, see for instance \cite[Corollary 4.49]{HS25}.

\medskip

\noindent {\bf Step 2:} We know that $G^*(X,Y) = G(Y,X)$, where $G^*$ is the Green function for $L^*$, see for instance \cite[Lemma 14.79]{DFMmixed}. In order to prove that $N^*(X,Y) = N(Y,X)$, it is thus enough to prove that $v(X,Y) = v^*(X,Y)$, where $v^*(.,Y) \in \widehat W^{1,2}(\Omega)$ satisfies 
\begin{equation} \label{Av*phi}
\int_\Omega A^*\nabla v_Y \cdot \nabla \varphi \, dX = - \fint_{\Omega} \varphi \, dX + \int_{\partial \Omega} \varphi \, d\omega_{L}^Y \quad \text{ for } \varphi \in W^{1,2}(\Omega).
\end{equation}
By using both \eqref{Avphi} and \eqref{Av*phi}, we obtain that 
\[ \int_\Omega A \nabla v_X \cdot \nabla v^*_Y\, dZ = \int_{\partial \Omega} v^*_Y d\omega_{L^*}^X = \int_{\partial \Omega} v_X d\omega_L^Y,\]
so 
\[\int_{\partial \Omega} v^*(.,Y) d\omega_{L^*}^X + \fint_\Omega G^*(Z,Y) \, dZ = \int_{\partial \Omega} v(.,X) d\omega_L^Y + \fint_\Omega G(Z,X) \, dZ,\]
We use the fact that  $G^*(Z,Y) = G(Y,Z)$ and $G(Z,X) = G^*(X,Z)$, and the above identity becomes
\[\int_{\partial \Omega} v^*(.,Y) d\omega_{L^*}^X - \fint_\Omega G^*(X,Z) \, dZ = \int_{\partial \Omega} v(.,X) d\omega_L - \fint_\Omega G(Y,Z) \, dZ.\]
However, the left hand side is solution to $L^*u=-1$ in $\Omega$ with boundary data $v^*(.,Y)|_{\partial \Omega}$ evaluated at the point $X$, which is $v^*(X,Y)$. Similarly, the right hand side is $v(Y,X)$. It gives the desired equality.

\medskip

\noindent {\bf Step 3:} If $\vec F \in C^\infty_0(\Omega)$ and $h \in C^\infty_0(\Omega)$ with zero average, then Lemma 14.87 in \cite{DFMmixed} shows that the function $u^D_{\vec F,h}$ defined as
\[ u^D_{\vec F,h}(X):= \int_\Omega [\nabla_Z G(X,Z) \cdot \vec F(Z) + G(X,Z) h(Z)]\, dZ \]
lies in $W^{1,2}(\Omega)$ and verifies
\[  \int_\Omega A\nabla u^D_{\vec F,h} \cdot \nabla \varphi \, dX = \int_\Omega [\vec F \cdot \nabla \varphi + h\varphi] \, dX \qquad \text{ for } \varphi \in W_0^{1,2}(\Omega),\]
that is
\[  \int_\Omega A\nabla u^D_{\vec F,h} \cdot \nabla \varphi \, dX = \int_\Omega [\vec F \cdot \nabla (\varphi-\varphi_{L^*}) + h(\varphi - \varphi_{L^*})] \, dX \qquad \text{ for } \varphi \in C^\infty(\R^n),\]
where $\varphi_{L^*}$ is introduced in Step 1. The function $v_{\vec F,h}$ defined by
\[ v_{\vec F,h}(X):= \int_\Omega [\nabla_Z v(X,Z) \cdot \vec F(Z) + v(X,Z) h(Z)]\, dZ \]
lies also in $W^{1,2}(\Omega)$ and verifies 
\[  \int_\Omega A\nabla v_{\vec F,h} \cdot \nabla \varphi \, dX = \int_\Omega [\vec F \cdot \nabla \varphi_{L^*} + h \varphi_{L^*}] \, dX \qquad \text{ for } \varphi \in \widehat W^{1,2}(\Omega)\]
by a simple differentiation under the integral, since $v(.,Y) \in W^{1,2}(\Omega)$ and $-\diver \vec F + h \in C^\infty_0(\Omega)$. It entails that $u_{\vec F, h} = v_{\vec F,h} + u^D_{\vec,h} \in \widehat{W}^{1,2}(\Omega)$ verifies
\begin{equation} \label{LMuFh} \int_\Omega A\nabla u_{\vec F,h} \cdot \nabla \varphi \, dX = \int_\Omega [\vec F \cdot \nabla \varphi + h \varphi] \, dX \qquad \text{ for } \varphi \in \widehat W^{1,2}(\Omega),
\end{equation}
meaning that $u_{\vec F,h}$ is the solution to $Lu = h -\diver \vec F $ with zero Neumann data given by the Lax-Milgram theorem. 

\medskip

When $\vec F \in L^2(\Omega)$, $h \in L^2(\Omega)$ with zero average, we want to make sense to the expression \eqref{defuFh} as the limit of the $u_{\vec F_n,h_n}$, where $\vec F_n \in C_0^\infty(\Omega)$ and $h_n \in C_0^\infty(\Omega)$ with zero average are such that $\vec F_n \to \vec F$ and $h_n \to h$ in $L^2$. Indeed, the linear forms
\[ \varphi \in \widehat W^{1,2}(\Omega) \mapsto \int_\Omega [\vec F_n \cdot \nabla \varphi + h_n \varphi ] dX\]
are a Cauchy sequence in the space of forms on $\widehat W^{1,2}(\Omega)$. The Lax-Milgram implies then that $u_{\vec F_n,h_n}$ is a Cauchy sequence in $\widehat W^{1,2}(\Omega)$  that converges to a function $u_{\vec F,h} \in \widehat W^{1,2}(\Omega)$ that verifies \eqref{LMuFh}. Moreover, the limit $u_{\vec F,h}$ does not depend  on the approximation of $\vec F$ and $h$. The Theorem follows.
\ep

Finally, we conclude the paragraph with a stability result.

\begin{theorem} \label{ThfixedpointW}
Let $\Omega$ be a bounded 1-sided CAD, and let $L_i=-\diver A_i \nabla$, $i\in \{0,1\}$ be two uniformly elliptic operators. Let $N_i(X,Y)$, $i\in \{0,1\}$ be the Neumann functions of the two operators.

For any $h\in L^2(\Omega)$ with zero average in $\Omega$, and any $\vec F \in L^2(\Omega,\R^n)$, we define the map $\mathscr{N}_{\vec F} : \, \widehat W^{1,2}(\Omega) \to \widehat W^{1,2}(\Omega)$ as
\begin{equation} \label{defmapN}
\mathscr{N}_{\vec F,h}(v)(Y) := \int_\Omega \Big[ \nabla_X N_0(X,Y) \cdot [(A_0-A_1)(X)\nabla v (X)+ \vec F(X)]\,  + N_0(X,Y) h(X) \Big] dX.
\end{equation} 

Then there exists $\epsilon>0$ (depending only on the ellipticity and boundedness of $L_0$, but not on $\vec F$ and $h$) such that when $\|A_0 - A_1\|_{L^\infty(\Omega)} < \epsilon$ and $v\in  \widehat W^{1,2}(\Omega)$, the sequence $\{(\mathscr{N}_{\vec F,h})^k(v)\}_{k\in \bN}$ converges in $W^{1,2}(\Omega)$ to the unique fixed point of $\mathscr{N}_{\vec F,h}$. This fixed point $u$ verifies
\begin{equation} \label{fixedpoint1}
u(X) = \int_\Omega \Big[ \nabla_Y N_1(X,Y) \cdot \vec F(Y) + N_1(X,Y) h(Y) \Big] \, dY \quad \text{ for a.e. } X \in \Omega,
\end{equation}
or equivalently
\begin{equation} \label{fixedpoint2}
\int_\Omega A_1\nabla u\cdot \nabla \phi \, dX = \int_\Omega \Big[ \vec F \cdot \nabla  \phi + h\phi\Big]  \, dX \quad \text{ for any } \phi \in \widehat W^{1,2}(\Omega).
\end{equation}
\end{theorem}

\bp
One can immediately recognize Picard's fixed point theorem. Since we already know that $\widehat{W}^{1,2}(\Omega)$ is a Banach space, we need to check whether $\mathscr N_{\vec F,h}$ is a contraction mapping.

Since the $\mathscr N_{\vec F,h}(v)$ is the weak solution in $\widehat{W}^{1,2}(\Omega)$ to 
\[L \mathscr N_{\vec F,h}(v) = h - \diver [(A_0-A_1)\nabla v (X)+ \vec F] \ \text{ in $\Omega$ with zero Neumann data,}\]
 see  (ii) of Theorem \ref{ThNeumann}, then by the Lax-Milgram theorem,
\begin{multline*} 
\|\mathscr{N}_{\vec F,h}(v-w)\|_{\dot W^{1,2}(\Omega)} \leq C \|(A_0-A_1) \nabla (v-w)|\|_{L^2(\Omega)} \\
\leq C\|A_0-A_1\|_\infty \|v-w\|_{\dot W^{1,2}(\Omega)} \leq \frac12 \|v-w\|_{\dot W^{1,2}(\Omega)}.
\end{multline*}
if $\|A_0-A_1\|_\infty$ is small enough. Hence the map is contractive.

As for the fixed point $u$, it is the solution in $\widehat{W}^{1,2}(\Omega)$  to $L_0 u = h - \diver [(A_0-A_1)\nabla u + \vec F]$ with zero Neumann data and zero average, meaning by  (ii) of Theorem \ref{ThNeumann} that $u$ is the solution function in $\widehat{W}^{1,2}(\Omega)$ verifying 
\[  \int_\Omega A_0\nabla u\cdot \nabla \phi \, dX = \int_\Omega \Big( [(A_0-A_1)\nabla u + \vec F] \cdot \nabla  \phi + h\phi\Big)\, dX \quad \text{ for any } \phi \in \widehat W^{1,2}(\Omega).\]
This immediately implies \eqref{fixedpoint2}, which is equivalent to  \eqref{fixedpoint1}, again by (ii) of Theorem \ref{ThNeumann}.
\ep

\begin{theorem} \label{ThwPN} Let $\Omega$ be a bounded 1-CAD, $L:=-\diver A \nabla$ be a uniformly elliptic operator on $\Omega$, and $p\in (1,\infty)$. 
Then the following are equivalent:
\begin{enumerate}[(i)]
\item The Poisson-Neumann problem \eqref{PNp}$_{L,\Omega}$ is solvable in $L^p$, meaning that there exists $C>0$ such that, for any $h \in L^\infty_c(\Omega)$, the function $u_{\vec F,0}$ defined by
\[ u_{\vec F,0}(X) := \int_\Omega \nabla_Y N(X,Y) \cdot \vec F (Y) \, dY \]
satisfies
\[ \|\N^2(u_{0,h})\|_{L^{p'}(\partial \Omega)} \leq C \|\vec F\|_{T^{p,1,2}_{-1}(\Omega)}.\]
\item The Poisson-Neumann-Regularity problem for $L^*$ is solvable in $L^{p'}$, meaning that there exists $C>0$ such that, for any $h \in L^\infty_c(\Omega)$, the function $u^*_{0,h}$ defined by
\[ u_{0,h}^*(X) := \int_\Omega N^*(X,Y)h(Y) \, dY \]
satisfies
\[ \|\N^2(\nabla u_{0,h}^*)\|_{L^{p'}(\partial \Omega)} \leq C \|h\|_{T^{p,1,2}_{-1}(\Omega)}.\]
\end{enumerate}
The tent space $T^{p,1,2}_{-1}(\Omega)$ was defined in Definition \ref{defTp}.
\end{theorem}

\bp
This ``theorem'' is a simple consequence of the duality of tent spaces. The proof can be found as \cite[Lemma 3.1]{FL-PN}. Note that our choice of the Neumann function $N(X,Y)$ --- with zero average on $\Omega$ instead of zero average on the trace --- render useless the need to remove the average of $u_{\vec F,0}$ in the Poisson-Neumann problem.
\ep

To finish the subsection, we give relation between the Poisson-Neumann problem and a localization property. First, we start with the definition of localization.

\begin{definition}
Given $\Omega$ a bounded 1-CAD, $L:=-\diver A \nabla$ a uniformly elliptic operator on $\Omega$, and $p\in (1,\infty)$, we say that \eqref{Locp} --- or \eqref{Locp}$_{L^*}$,  \eqref{Locp}$_{L^*,\Omega}$ --- holds if there exists $C>0$ such that for any ball $B$ centered on $\partial \Omega$ and any weak solution $u\in W^{1,2}(2B\cap \Omega)$ to $L^*u=0$ with zero Neumann data, we have
\begin{equation} \label{Locp} \tag{Loc$_{p'}$}
\|\N^2(\nabla u \1_{B})\|_{L^{p'}(\partial \Omega)} \leq C \left(\fint_{2B\cap \Omega} |\nabla u|^2 dX\right)^\frac12.
\end{equation}
\end{definition}

The following proposition shows that we can change $2B$ by $\kappa B$ in the right-hand side of \eqref{Locp} and limit ourselves to balls of sufficiently small radius.

\begin{proposition} \label{prLocp}
Let $\Omega$ be a bounded 1-CAD, $L:=-\diver A \nabla$ be a uniformly elliptic operator on $\Omega$, $p\in (1,\infty)$, $r_0\in (0,\diam\Omega)$, and $\kappa >1$. 

Then the localization property \eqref{Locp} holds if and only if there is a $C>0$ such that for any ball $B$ centered on $\partial \Omega$ and radius $r\leq r_0$ and any weak solution $u\in W^{1,2}(2B\cap \Omega)$ to $L^*u=0$ with zero Neumann data, we have
\begin{equation}
\|\N^2(\nabla u \1_{B})\|_{L^{p'}(\partial \Omega)} \leq C \left(\fint_{\kappa B\cap \Omega} |\nabla u|^2 dX\right)^\frac12.
\end{equation}
\end{proposition}

\bp
Let $r_0 \in (0,\diam\Omega)$ and $\kappa >1$. Assume that there exists $C>0$ such that for any ball $B$ centered on $\partial \Omega$ and radius $r\leq r_0$ and any weak solution $u\in W^{1,2}(2B\cap \Omega)$ to $L^*u=0$ with zero Neumann data, we have
\begin{equation} \label{LocpA1}
\|\N^2(\nabla u \1_{B})\|_{L^{p'}(\partial \Omega)} \leq C \left(\fint_{\kappa B\cap \Omega} |\nabla u|^2 dX\right)^\frac12.
\end{equation}
To prove the proposition, it suffices to show that for any $\lambda \in (1,\kappa)$, we have $C_\lambda$ --- that depends on ($\lambda$,) $r_0$ and $\kappa$ --- such that that for any ball $B$ centered on $\partial \Omega$ and any weak solution $u\in W^{1,2}(2B\cap \Omega)$ to $L^*u=0$ with zero Neumann data, we have
\begin{equation}
\|\N^2(\nabla u \1_{B})\|_{L^{p'}(\partial \Omega)} \leq C \left(\fint_{\lambda B\cap \Omega} |\nabla u|^2 dX\right)^\frac12.
\end{equation}

\medskip

The proof is a simple consequence of a covering lemma. Let $B=:B(x,r)$ be a ball centered on the boundary and of radius $r\in (0,\diam(\Omega)$. The collection $\{B(z,r r_0(\lambda-1)/10\kappa \diam(\Omega))\}_{z\in B\cap \partial \Omega}$ covers $B\cap \partial \Omega$, so by the Vitali lemma, we can find a non-overlapping subcollection  $\{B_{i}\}_{i\in I_1}$ such that $B\cap \partial \Omega \subset \bigcup_{i\in I_1} 5B_i$. 
Let $S$ be the Whitney set $B \setminus \bigcup_{i\in I_1} 10B_i$. The collection of balls $\{\frac15 \min\{1, \lambda-1\} B_X\}_{X\in S}$ covers $S$, so by the Vitali lemma again,  we can find a non-overlapping subcollection  $\{B_{i}\}_{i\in I_2}$ such that $S \subset \bigcup_{i\in I_2} 5B_i$. Since the two collections are non-overlapping, it forces $\# I_1 + \# I_2$ to be bounded by a constant $N$ that depends only on $r_0(\kappa-1)/10\lambda \diam(\Omega)$ and the dimension $n$. 

For $i\in I_1$, we have that $10 B_i$ has radius less than $r_0$, so \eqref{LocpA1} gives that 
\[ \|\N^2(\nabla u \1_{10B_i})\|_{L^{p'}(\partial \Omega)} \lesssim \left(\fint_{\kappa B_i \cap \Omega} |\nabla u|^2 dX\right)^\frac12 \lesssim \left(\fint_{\lambda B \cap \Omega} |\nabla u|^2 dX\right)^\frac12  . \]
since $10\kappa B_i$ is included in $\lambda B$ and both ball have comparable radius (up to a constant that depends on $\lambda$, $\kappa$ and $r_0$). 

For $i \in I_2$, we simply observe that $5B_i$ is a Whitney ball included in  $\lambda B$, so  
\[  \|\N^2(\nabla u \1_{5B_i})\|_{L^{p'}(\partial \Omega)} \lesssim \left(\fint_{5B_i} |\nabla u|^2 dX\right)^\frac12 \lesssim \left(\fint_{\lambda B \cap \Omega} |\nabla u|^2 dX\right)^\frac12.\]

Altogether, 
\begin{multline*}
     \|\N^2(\nabla u \1_{B})\|_{L^{p'}(\partial \Omega)} \leq \sum_{i \in I_1}\| \N^2(\nabla u \1_{10B_i})\|_{L^{p'}(\partial \Omega)} +  \sum_{i \in I_2} \|\N^2(\nabla u \1_{5B_i})\|_{L^{p'}(\partial \Omega)} \\
     \leq C N  \left(\fint_{\lambda B \cap \Omega} |\nabla u|^2 dX\right)^\frac12,
\end{multline*}

as desired for the lemma.
\ep

\begin{theorem} \label{ThLocp}
Let $\Omega$ be a bounded 1-CAD, $L:=-\diver A \nabla$ be a uniformly elliptic operator on $\Omega$, and $p\in (1,\infty)$.
\begin{enumerate}[(a)]
\item If the Poisson-Neumann problem \eqref{PNp}$_{L}$ is solvable, then \eqref{Locp}$_{L^*}$ holds. Moreover, the constant in \eqref{Locp}$_{L^*}$ depends only on $p$, the constant in \eqref{PNp}$_{L}$, the dimension $n$, and the 1-CAD constants of $\Omega$.

\item If \eqref{Locp}$_{L^*}$ holds, and if there is a $q\in (1,\infty)$ such that the the Poisson-Neumann problem {\rm(PN$_{q}$)}$_{L}$ is solvable, then the Poisson-Neumann problem \eqref{PNp}$_L$ is solvable as well.
\end{enumerate}
\end{theorem} 

\bp
Part (a) is \cite[Lemma 4.1]{FL-PN}, and Part (b) is a combination of \cite[Theorem 1.20, Theorem 1.19, and Proposition 1.17]{FL-PN}. Both implications are more easily seen on the implication charts in \cite[Subsection 1.3]{FL-PN} that summarize the results of the article.
\ep

We will need the theorem only via the following corollary.

\begin{corollary} \label{CorLocp}
Let $\Omega$ be a bounded Lipschitz domain and $p\in (1,\infty)$. Then the following are equivalent:
\begin{enumerate}[(i)]
\item The Poisson-Neumann \eqref{PNp}$_{-\Delta}$ is solvable,
\item The localization property \eqref{Locp}$_{-\Delta}$ holds.
\end{enumerate}
\end{corollary}

\bp
Immediate consequence of Theorems \ref{ThRellich} and \ref{ThLocp}.
\ep

\subsection{Caccioppoli inequalities}

\begin{lemma}[Double Caccioppoli inequality] \label{lemDCaccio}
Let $\Omega$ be a $1$-CAD, and let $L:=-\diver A \nabla$ be a uniformly elliptic operator on $\Omega$. There there exists $C>0$ such that for any ball $B:=B(X_B,r_B)$ satisfying $4B\subset \Omega$, any cutoff function $\varphi \in C^\infty_0(\Omega)$, any positive weak solution $v \in W^{1,2}(3B)$ to $L^*v = 0$ in $3B$, any $\vec F \in L^2(2B)$ and any weak solution $u\in W^{1,2}(2B)$ to $Lu = - \diver \vec F$, we have 
\begin{equation} \label{DCaccio1aa}
\fint_{B}  \frac{|\nabla v|^2}{v^2} u^2 \varphi^4 \, dX \lesssim \frac{1}{r_B^2} \fint_{2B_X}  u^2 \varphi^2 \, dX + \fint_{2B_X} |\vec F|^2 \varphi^4 \, dX.
\end{equation}
\end{lemma}

\bp
As the name suggests, the proof relies on using a Caccioppoli type argument twice. The case $\vec F = 0$ was treated in \cite[(5.20)]{FP23}, and the general case is not much more complicated. 

\medskip

Let $B$, $\varphi$, $v$ and $u$ be as in the Lemma. Let $\psi \in C^\infty_0(2B)$ be such that $0 \leq \psi \leq 1$, $\psi \equiv 1$ on $B$, and $r |\nabla \Psi| \leq 2$.
Since $v$ is positive in $3B$, the Harnack inequality shows that $v \approx v(X_B)$ on $2B$, and thus the bound \eqref{DCaccio1aa} will be proven as long as we show that 
\begin{equation} \label{DCaccio1ab}
 \fint_B |\nabla v|^2 u^2 \varphi^4 \, dX \lesssim v(X_B)^2 \Big(\frac{1}{r_B^2} \fint_{2B_X}  u^2 \varphi^2 \, dX + \fint_{2B_X} |\vec F|^2 \varphi^4 \, dX\Big)
 \end{equation}
We write $\phi$ for $\varphi \psi$. Note that since $4B \subset \Omega$, $r \leq \delta$ on $2B \supset \supp \phi$, and so $|\nabla \phi| \lesssim \frac1{r_B}$. The bound \eqref{DCaccio1ac} is now a consequence of
\begin{equation} \label{DCaccio1ac}
 I_0 := \int_\Omega |\nabla v|^2 u^2 \phi^4 \, dX \lesssim v(X_B)^2 \Big(\int_{\Omega}  u^2 \phi^2 |\nabla \phi|^2 \, dX + \int_{\Omega} |\vec F|^2 \phi^4 \, dX\Big).
 \end{equation}
 
 \medskip
 
 We first show that 
 \begin{equation} \label{DCaccio1ad}
 T_0 := \int_\Omega |\nabla u|^2 \phi^4 \, dX \lesssim \int_{\Omega}  u^2 \phi^2 |\nabla \phi|^2 \, dX + \int_{\Omega} |\vec F|^2 \phi^4 \, dX =: J_1 + J_2.
 \end{equation}
 Indeed, by the ellipticity of $A$, we have
\[
 T_0 \lesssim \int_\Omega A\nabla u\cdot \nabla u \phi^4 \, dX = \int_{\Omega} A\nabla u \cdot \nabla [u\phi^4] \, dX - 4 \int_\Omega A\nabla u \cdot \nabla \phi \, [u\phi^3]\, dX =: T_1 + T_2.
\]
 The term $T_2$ is bounded with the help of the Cauchy-Schwarz inequality and the boundedness of $A$ by
 \[ T_2 \lesssim  \left(  \int_\Omega |\nabla u|^2 \phi^4 \, dX \right)^\frac12 \left( \int_\Omega |\nabla \phi|^2 u^2 \phi^2 \right)^\frac12 = (T_0)^\frac12 (J_1)^\frac12.\]
 As for $T_1$, we use the fact that $u$ is a weak solution to $Lu = -\diver \vec F$ on $2B$ and we get 
\[
 T_1= \int_\Omega \vec F \cdot  \nabla [u\phi^4] \, dX = \int_\Omega \vec F \cdot  \nabla u\, \phi^4 \, dX + 4 \int_\Omega \vec F \, \nabla \phi \, u \phi^3 =: T_3+T_4.
\]
With the Cauchy-Schwarz inequality, $T_3 \leq (T_0)^\frac12 (J_2)^\frac12$ and $T_4 \leq 4 (J_1)^\frac12 (J_2)^\frac12$. Altogether, 
\[ T_0 \lesssim (T_0)^\frac12 (J_1+J_2)^\frac12 + (J_1J_2)^\frac12,\]
which easily self improves into the desired claim \eqref{DCaccio1ad}, since $T_0$ is finite.

\medskip

To prove \eqref{DCaccio1ac}, we use the ellipticity of $A$ to get
\begin{multline*}
I_0 \lesssim \int_\Omega A\nabla v \cdot \nabla v \, [u^2\phi^4]\, dX \\
=  \int_\Omega A\nabla [vu^2\phi^4] \cdot \nabla v \, dX - 2 \int_\Omega A\nabla u \cdot \nabla v \, [vu\phi^4] \, dX - 4 \int_\Omega A\nabla \phi \cdot \nabla v \, [vu^2\phi^3] \, dX\\=: I_1+I_2+I_3.
\end{multline*}
The term $I_1$ is 0, since $v$ is a weak solution to $L^* v = 0$. The term $I_2$ and $I_3$ are bounded with the Cauchy-Schwarz inequality, the boundedness of $A$, and the fact that $v\approx v(X_B)$ on $2B$ as follows:
\[ I_2 \lesssim (I_0)^\frac12 \left(  \int_\Omega |\nabla u|^2 \phi^4 v^2 \, dX \right)^\frac12 \lesssim  v(X_B) (I_0)^\frac12 (T_0)^\frac12\]
and
\[ I_3 \lesssim (I_0)^\frac12 \left(  \int_\Omega |\nabla \phi|^2 \phi^2 u^2 v^2 \, dX \right)^\frac12 \lesssim  v(X_B) (I_0)^\frac12 (J_1)^\frac12.\]
We conclude that 
\[ I_0 \lesssim (I_0)^\frac12 v(X_B) [J_1 + T_0]^\frac12 \lesssim (I_0)^\frac12 v(X_B) [J_1 + J_2]^\frac12\]
with our bound on $T_0$, which again self improves into the claim \eqref{DCaccio1ac}. The lemma follows.
\ep

\begin{lemma} \label{lemCaccioN2u}
Let $\Omega \subset \R^n$, and let $L = - \diver A \nabla$ be a uniformly elliptic operator satisfying $|\delta \nabla A| \leq K$. Then for any $X\in \Omega$ and any $h\in L^2(2B_X)$, any weak solution $u\in W^{1,2}(\Omega)$ to $Lu=h$ in $2B_X$ lies in $W^{2,2}_{loc}(2B_X)$. More precisely, for any cut-off function $\varphi \in C^\infty_0(2B_X)$, we have
\[ \int_{\Omega} |\nabla^2 u|^2 \varphi^4 \, dZ \leq C \left(  \frac1{\delta(X)^2} \int_{\Omega} |\nabla u|^2 \varphi^2 \, dZ + \int_{\Omega} h^2 \varphi^4 \, dZ \right),\]
where $C$ depends only on $K$, the cutoff constants of $\varphi$, and the elliptic and boundedness constants of $L$.
\end{lemma}

\bp
The lemma is a variant of (i) of Lemma \ref{lembdonG}. We prove it for completeness.

\medskip

One difficulty of the proof is the fact that we do not know {\em a priori} that $u\in W^{2,2}_{loc}(2B_X)$. However, this fact was proven in \cite[Theorem 8.8]{GT83}, so we will not show it here.

\medskip

By the ellipticity and boundedness of $A$, we have 
\[ I := \int_{\Omega} |\nabla^2 u|^2 \varphi^4 \, dZ \approx \sum_{i=1}^n \int_{\Omega} A \nabla [\partial_i u] \cdot \nabla [\partial_i u] \varphi^4 \, dZ.\]
Define
 \[J :=  \frac1{\delta(X)^2} \int_{\Omega} |\nabla u|^2 \varphi^2 \, dZ + \int_{\Omega} h^2 \varphi^4 \, dZ.\]
 Since $I$ is finite, to prove the lemma, it suffices to thus show that
 \[ I_0 := \int_\Omega  A \nabla [\partial u] \cdot \nabla [\partial u] \varphi^4 \, dZ \lesssim (IJ)^{1/2} + J,\]
 where $\partial$ stands for any partial derivative of $u$.
 
 \medskip
 
 We first write
 \[ I_0 = \int_\Omega  A \nabla [\partial u] \cdot \nabla [\varphi^4 \partial u] \, dZ  - 4 \int_\Omega A \nabla [\partial u] \cdot \nabla \varphi \, \varphi^3 \partial u \, dZ =: I_1 - 4I_2.\]
We deal with $I_2$ with the boundedness of $A$ and the Cauchy-Schwarz inequality
\[ I_2 \lesssim I^\frac12 \left( \int_\Omega |\nabla u|^2 |\nabla \varphi|^2 \varphi^2 \, dZ \right)^\frac12 \lesssim (IJ)^{1/2}\]
by using one property of cut-off functions: the fact that $|\delta \nabla \varphi| \lesssim 1$, which implies $|\nabla \varphi| \lesssim 1/\delta(X)$ since $\delta \approx \delta(X)$ on $\supp \varphi \subset 2B_X$.

To treat the term $I_1$, we first notice that $\partial u \in W^{1,2}_{loc}(2B_X)$ by \cite[Theorem 8.8]{GT83} and then it is simple to observe that $v = \partial u$ is a weak solution to $Lv = \partial h + \diver (\partial A) \nabla u$. Consequently, 
\begin{multline*} 
I_1 = - \int_\Omega h \partial [\varphi^4 \partial u]\, dZ  - \int_\Omega  (\partial A) \nabla u \cdot \nabla [\varphi^4 \partial u] \, dZ \\
= - 4 \int_\Omega h \partial \varphi \, (\varphi^3 \partial u)\, dZ - \int_\Omega h \partial^2 u\,  \varphi^4 \, dZ  - \int_\Omega  (\partial A) \nabla u \cdot \nabla [\partial u] \, \varphi^4  \, dZ \\
- 4  \int_\Omega  (\partial A) \nabla u \cdot \nabla \varphi \,  [\varphi^3 \partial u]  \, dZ =: I_3 + I_4 + I_5 + I_6.
\end{multline*}
The term $I_3$ is bounded with the help of the Cauchy-Schwarz inequality, and we obtain
\[|I_3| \lesssim \left( \int_\Omega h^2 \varphi^4 \, dZ \right)^\frac12 \left( \frac1{\delta(X)} \int_\Omega |\nabla u|^2 \varphi^2 \, dZ \right)^\frac12 \lesssim J \]
since $|\nabla \varphi| \lesssim 1/\delta(X)$ as before. For the term $I_4$, we easily have
\[|I_4| \leq I^\frac12 \left( \int_\Omega h^2 \varphi^4 \, dZ \right)^\frac12 \leq (IJ)^\frac12. \]
The term $I_5$ is similar to $I_2$, since $|\partial A| \leq |\nabla A| \lesssim 1/\delta(X)$ by assumption on $A$, and we have
\[ I_5 \lesssim I^\frac12 \left( \int_\Omega |\nabla u|^2 |\partial A|^2 \varphi^4 \, dZ \right)^\frac12 \lesssim (IJ)^{1/2}.\]
Finally, for $I_6$, we use both $|\nabla \varphi| \lesssim 1/\delta(X)$ and $|\partial A| \lesssim 1/\delta(X)$ to obtain
\[I_6 \lesssim \int_\Omega |\nabla u|^2 |\partial A| |\nabla \varphi| \varphi^3 \, dZ \lesssim J.\]
 The lemma follows.
\ep

\subsection{An integration by parts in Chord-Arc Domains}

\begin{lemma} \label{IPPwithNG}
Let $\Omega \subset \R^n$ be a 1-sided CAD, and write $B_\Omega$ for the center ball of $\Omega$ if $\Omega$ is bounded, and $B_\Omega = \emptyset$ when $\Omega$ is unbounded. Write $G$ be the Green function (for $-\Delta$) with pole at center/infinity defined in Definition \ref{defomG}. 

Assume that $\vec F \in W^{1,2}_{loc}(\Omega,\R^{n})$ is a vector  such that $\vec F \cdot \nabla G = 0$ for a.e. in $\Omega$. Then for any $v\in W^{1,2}_{loc}(\Omega)$, any $\phi \in C^\infty_c(\Omega \setminus \overline{2B_\Omega}))$ and $c>0$,
\begin{multline} \label{ExtraIPP}
\left|\int_{\Omega}\vec F \cdot \nabla v\, \phi\,  dX\right| \lesssim \int_{\Omega}  |\nabla v| \, \delta \big[ |\vec F| |\nabla \phi| + |\phi| |\nabla \vec F| \big]\, dX  \\
+  \int_{\Omega}  \Big(\frac{\delta^2 \big|\nabla^2 G\big|}{G} + \1_{\set{\delta |\nabla G|<cG}}\Big) |\vec F| |\nabla v| |\phi| \, dX,
\end{multline}
with a constant that depends on the dimension $n$, the Ahlfors regular constant of $\partial \Omega$, and $c$.
\end{lemma}

\bp Since the Green function $G$ is a positive harmonic function on $\Omega \setminus B_\Omega$, we have by Lemma \ref{lembdonG} that
\begin{equation} \label{bdonN^2G} \delta^2 |\nabla^2 G| + \delta |\nabla G| \lesssim G \quad \text{ on } \Omega \setminus 2B_\Omega.
\end{equation}
Let $c$ be as in the Lemma. We construct $\psi : \Omega \to \R$ to be a smooth cut-off function such that $\supp \psi \subset \{\delta|\nabla G| > cG/2\}$ and $\psi \equiv 1$ on $\{\delta|\nabla G| > cG\}$. 
Using \eqref{bdonN^2G}, one can check that 
\[ \left| \nabla \left( \frac{\delta |\nabla G|}{G} \right)\right|  \leq \frac{|\nabla \delta| |\nabla G|}{G} +\frac{\delta |\nabla^2 G|}{G}  + \frac{\delta |\nabla G|^2}{G^2} \lesssim \frac1{\delta},\]
meaning that we can choose $\psi$ such that $|\nabla \psi| \lesssim 1/\delta$.

\medskip 

Finally, we can assume that $v\in W^{2,2}_{loc}(\Omega)$. Indeed, once the inequality \eqref{ExtraIPP} is proven for $v\in W^{2,2}_{loc}(\Omega)$, we can extend it to $v\in W^{1,2}_{loc}(\Omega)$ by density of smooth functions in $W^{1,2}(\supp \phi)$. With this set up done, we introduce the Green function by writing $1=\frac{\nabla G \cdot \nabla G}{|\nabla G|^2}$ and we write
\[
I := \int_\Omega \vec F \cdot \nabla v \, \phi\, dX = \int_\Omega (\vec F \cdot \nabla v) \, \phi (1-\psi)\,  dX +  \int_\Omega (\vec F \cdot \nabla v) \, \phi\psi \, \frac{\nabla G \cdot \nabla G}{|\nabla G|^2}  dX =: I_1 + I_2.
\]
The bound on $I_1$ is straightforward: we have
\[ |I_1| \leq \int_{\Omega} |\vec F| |\nabla v| |\phi| \1_{\{\delta|\nabla G| \leq cG\}} \, dX,\]
which is one of terms in the right-hand side of \eqref{ExtraIPP}. As for $I_2$, we proceed to an integration by parts to move one of the gradient away from (one) $G$. We obtain
\begin{multline*}
I_2 =  \int_\Omega (\vec F \cdot  \nabla v) \, \frac{\psi \phi G}{|\nabla G|^2} \, (-\Delta G)\,  dX - 
\int_\Omega (\vec F \cdot  \nabla v) \, (\nabla [\phi\psi] \cdot \nabla G) \frac{G}{|\nabla G|^2} \,  dX \\ 
+ 2 
\int_\Omega (\vec F \cdot  \nabla v) \, \frac{\psi\phi G}{|\nabla G|^3} \, \nabla G \cdot \nabla |\nabla G|\,  dX - \int_\Omega \nabla(\vec F \cdot  \nabla v) \cdot \nabla G \, \frac{\psi\phi G}{|\nabla G|^2} \, dX \\
=: I_3 + I_4 + I_5 + I_6. 
\end{multline*}
The first 3 terms ($I_3$, $I_4$, and $I_5$) are the easy ones. First, $I_3 = 0$ since $G$ is a solution to $-\Delta G=0$ on $\supp \phi$. For $I_4$, we use the fact that $\delta |\nabla G| \approx G$ on the support of $\psi$  to get 
\[ |I_4| \lesssim \int_{\Omega} |\vec F| |\nabla v| \, \delta (|\phi \nabla \psi| + |\psi \nabla \phi| ) \, dX \lesssim \int_{\Omega} |\vec F | |\nabla v| (|\phi| \1_{\{\delta|\nabla G| \leq cG\}} + \delta|\nabla \phi|) \, dX,\]
since $\delta |\nabla \psi | \lesssim 1_{\supp \psi}$ by construction of $\psi$. For $I_5$, we use again the fact that $\delta |\nabla G| \approx G$ on the support of $\psi$ to obtain
\[ |I_5| \lesssim  \int_{\Omega} |\vec F| |\nabla v||\phi| \,  \frac{\delta ^2\big |\nabla^2 G\big |}{G} \, dX,\]
which is again a term in the right-hand side of \eqref{ExtraIPP}.

Before we move to $I_6$, let us introduce a notation: we write $w_G$ for the unit vector $\frac{\nabla G}{|\nabla G|}$. We let the reader check that on the support of $\psi$,
\begin{equation} \label{bdonNwG}
|\nabla w_G| \lesssim \frac{\delta |\nabla^2 G|}{G}.
\end{equation}
We write $w_G^k$ for the $k^{th}$ coordinates of $w_G$. So $I_6$ becomes
\begin{multline*}
I_6 = \sum_{1\leq j,k\leq n}  \int_\Omega w_G^k \partial_k (\vec F_j\, \partial_j v) \, \frac{\psi\phi G}{|\nabla G|} \, dX  
\\
=  \sum_{1\leq j,k\leq n}  \int_\Omega (w_G^k \partial_k \vec F_j) \, \partial_j v \, \frac{\psi\phi G}{|\nabla G|} \, dX  + \sum_{1\leq j,k\leq n}  \int_\Omega \vec F_j\, (w_G^k \partial_k\partial_j v) \, \frac{\psi\phi G}{|\nabla G|}\, dX \\ =: I_7 + I_8.
\end{multline*}
For $I_7$, we have the bound
\[ |I_7| \lesssim \int_{\Omega}  |\delta \nabla \vec F| |\nabla v||\phi| \, dX,\]
using $G/|\nabla G| \lesssim \delta$ on $\supp \psi$. 

The term $I_8$ is dealt by switching the $k$ and $j$ derivative, and by moving the $j$ derivative on the other terms by a last integration by parts. First, note that 
\[ w_G^k \partial_k \partial _j v = \partial_j ( w_G^k \partial_k v) - (\partial_j w_G^k)(\partial_k v),\]
meaning that 
\begin{multline*}
 \Big| I_8 -  \underbrace{\sum_{1\leq i,j,k\leq n}  \int_\Omega \vec F_j\, \partial_j (w_G^k \partial_k v) \, \frac{\psi\phi G}{|\nabla G|} \, dX}_{I_{9}} \Big| \lesssim \int_\Omega |\vec F| |\nabla w_G| |\nabla v| |\phi| \delta \, dX \\
 \lesssim \int_\Omega \frac{\delta^2|\nabla^2 G|}{G}  |\vec F| |\nabla v| |\phi| \, dX
 \end{multline*}
 thanks to \eqref{bdonNwG}. Then, we write
 \begin{multline*}
 I_{9} = \sum_{1\leq j\leq n} \int_\Omega \partial_j \Big(\vec F_j \frac{\psi\phi G}{|\nabla G|} \Big) \, (\nabla v \cdot w_G)\,  dX \\
  = \sum_{1\leq j\leq n} \int_\Omega  \vec F_j \frac{\psi\phi \partial_j G}{|\nabla G|} \, (\nabla v \cdot w_G)\,  dX + \sum_{1\leq j\leq n} \int_\Omega  (\partial_j \vec F_j) \frac{\psi\phi G}{|\nabla G|} \, (\nabla v \cdot w_G)\,  dX \\
 + \sum_{1\leq j\leq n} \int_\Omega  \vec F_j \frac{G \partial_j[\psi\phi]}{|\nabla G|}  \, (\nabla v \cdot w_G)\,  dX - \sum_{1\leq j\leq n} \int_\Omega  \vec F_j \frac{\psi\phi G \partial_j |\nabla G|}{|\nabla G|^2} \, (\nabla v \cdot w_G)\,  dX \\
 =: I_{10} + I_{11} + I_{12} + I_{13}.
 \end{multline*}
 The term $I_{10}$ is the only new one, as the other ones can be bounded as earlier terms. But by assumption, we have $\vec F \cdot \nabla G = 0$, and $I_{10}$ can be rewritten as
 \[ I_{10} = \int_\Omega \vec F\cdot \frac{\nabla G}{|\nabla G|} \, \psi\phi \, (\nabla v \cdot w_G)\, dX = 0.\]
 Finally, the term $I_{11}$ is bounded as $I_7$ by
 \[ |I_{11}| \lesssim \int_{\Omega}  |\delta \nabla \vec F||\nabla v||\phi| \, dX,\]
 the term $I_{12}$ is bounded as $I_4$ by
 \[ |I_{12}| \lesssim \int_{\Omega} |\vec F| |\nabla v| (|\phi| \1_{\{\delta|\nabla G| \leq cG\}} + \delta |\nabla \phi|)\, dX,\] 
and the term $I_{13}$ is bounded as $I_5$ as
 \[ |I_{13}| \lesssim \int_{\Omega} |\vec F| |\nabla v| |\phi| \, \frac{\delta ^2\big |\nabla^2 G\big |}{G} \, dX.\] 
 The lemma follows.\ep

\section{$S < N$ estimates}

\subsection{$S<N$ for the Poisson problem}
In this section, we aim to establish results analogous to those presented in \cite[Lemma 2.11]{KP95} and \cite[Lemma 4.7]{DFMpert}, but with the distinction that our results pertain to elliptic measures, whereas the preceding results were with respect to surface measures.
Contrary to \cite{KP95} and \cite{DFMpert}, we will separate the result in 2 parts: an integration by parts (Theorem \ref{ThS<N1}) and a good lambda argument (Theorem \ref{Thgoodlambda}), to make it easier to follow. 

\begin{theorem} \label{ThS<N1} 
Let $\Omega \subset \R^n$ be a 1-CAD, and let $L:=-\diver A \nabla$ be a uniformly elliptic operator. We write $\omega_L$ for the elliptic measure (associated to $L$) with pole at center/infinity.

For $\vec F \in L^2_{loc}(\Omega)$, a weak solution $u$ to $Lu=-\diver \vec F$ in $\Omega$, and a cut-off function $\varphi \in C^\infty_c(\Omega)$ (see Definition \ref{defcutoff}), we have
\[ \|\mathcal A(\varphi^2\delta\nabla u)\|_{L^{2}(\partial \Omega,\omega_L)} \leq C_\varphi \Big( \|\N^2(\varphi u)\|_{L^{2}(\partial \Omega,\omega_L)} + \|\mathcal A_1^2(\varphi\delta \vec F) \|_{L^{2}(\partial \Omega,\omega_L)} \Big),\]
where $C_\varphi$ depends only on the cut-off constants of $\varphi$, the ellipticity and boundedness of constants of $L$, the 1-sided CAD constants of $\Omega$.
\end{theorem}

\bp
Keep in mind that the presence of $\varphi$, which is smooth and compactly supported, guarantees the finiteness of all the quantities that will be discussed in the proof.

\smallskip

Assume first that $\varphi \in C^\infty_c(\Omega \setminus \overline{2B_\Omega})$, where $B_\Omega$ is the center of $\Omega$ when $\Omega$ is bounded, and $B_\Omega = \emptyset$ otherwise. 
We write $\omega$ for $\omega_L$ to lighten the notation, and we will also need $G$, the Green function with pole at infinity/center of $L^* = - \diver A^T \nabla$, which is now a solution on $\supp \varphi$ with our current assumption on $\varphi$.  The combination of the equivalence \eqref{FubiniAtoG} and the ellipticity and boundedness of $A$ entails that
\[\|\mathcal A(\varphi^2\delta\nabla u)\|_{L^{2}(\partial \Omega,\omega)}^2  \approx  I:= \int_\Omega (A\nabla u\cdot \nabla u)\, \varphi^4 G\, dX.\]
We write
\[ J:=  \underbrace{\|\N^2(\varphi u)\|_{L^{2}(\partial \Omega,\omega)}^2}_{J_1} + \underbrace{\|\mathcal A_1^2(\varphi\delta \vec F) \|_{L^{2}(\partial \Omega,\omega)}^2}_{J_2} \]
and we will prove that $I \lesssim J + (IJ)^{1/2}$, which will self improve to $I \lesssim J$ since $I$ is finite.

\medskip

From there, we can start the integration by parts, and we write
\begin{multline*}
I = \int_\Omega A \nabla u \cdot \nabla (u \, \varphi^4 G) \, dX - 4 \int_\Omega (A \nabla u \cdot \nabla \varphi )\, \varphi^3 G u \, dX - \int_\Omega (A\nabla u\cdot \nabla G) \, \varphi^4 u \, dX \\ =: I_1 + I_2 + I_3,
\end{multline*}
and we further rewrite $I_3$ as
\[
I_3 = -\frac12 \int_\Omega A\nabla (\varphi^4 u^2) \cdot \nabla G \, dX + 2 \int_\Omega (A \nabla \varphi \cdot \nabla G)\, \varphi^3 u^2\, dX =: I_4 + I_5.
\]
The term $I_4$ is 0, since $G$ is a weak solution to $L^* G=0$ on $\supp \varphi$, and $\varphi^4 u^2$ is a $W^{1,2}_{loc}(\Omega)$ function that is compactly supported in $\Omega$. For the term $I_5$, we use the fact that $\varphi$ is a cut-off function and Theorem \ref{ThCarlesonCO} to get
\begin{multline*}
|I_5| \lesssim \int_{\Omega} \1_{\supp\nabla\varphi}  \frac{\delta|\nabla G|}{G} \varphi^3 u^2 \, \frac{G}{\delta^2} dX \lesssim \Big\|\N^1\Big( \frac{\delta|\nabla G|}G u^2 \varphi^3 \Big) \Big\|_{L^1(\partial \Omega,\omega)} \\
\lesssim  \Big\|\N^1\Big( \frac{\delta^2|\nabla G|^2}{G^2} u^2 \varphi^4 \Big) \Big\|_{L^1(\partial \Omega,\omega)} + \|\N^1(u^2 \varphi^2)\|_{L^1(\partial \Omega,\omega)} \\
\lesssim  \Big\|\N^1_{(2,\kappa)}\Big( \frac{\delta^2|\nabla G|^2}{G^2} u^2 \varphi^4 \Big) \Big\|_{L^1(\partial \Omega,\omega)} + \|\N^1(u^2 \varphi^2)\|_{L^1(\partial \Omega,\omega)},
\end{multline*}
where $\kappa \in (0,1]$ is chosen small enough (depending only on the Corkscrew condition) such that $3\kappa B_X \cap B_\Omega = \emptyset$ when $X\in \Omega \setminus 2B_\Omega$. Since the radius of $\kappa B$ is comparable to $\delta$ in $2\kappa B$, and since $G$ is a positive solution to $L^*G=0$ in $2B$, the double Caccioppoli inequality (Lemma \ref{lemDCaccio}) gives then that, whenever  $4B_X \subset \Omega$
\begin{equation} \label{DCaccio1}
\fint_{\kappa B_X}  \frac{\delta^2|\nabla G|^2}{G^2} u^2 \varphi^4 \, dX \lesssim \fint_{2\kappa B_X}  u^2 \varphi^2 \, dX + \fint_{2\kappa B_X}  \delta^2|\vec F|^2 \varphi^4 \, dX.
\end{equation}
Altogether, 
\[|I_5| \lesssim \|\N^1_{(2,2\kappa)}(u^2 \varphi^2)\|_{L^1(\partial \Omega,\omega)}  + \|\N^1_{(2,2\kappa)}(\delta^2\varphi^2|\vec F|^2)\|_{L^1(\partial \Omega,\omega)} \lesssim J_1+J_2\]
since $\|\N^2_{(2,2\kappa)}(\delta \varphi |\vec F|)\|_{L^2(\omega)} \lesssim \|\A_1^2(\delta \varphi |\vec F|)\|_{L^2(\omega)}$ by Lemma \ref{lemcompNA}. We move to $I_2$, that we bound as follows:
\begin{multline*}
|I_2| \lesssim \int_\Omega |\nabla u| u \1_{\supp \nabla \varphi} \varphi^3 \frac{G}{\delta} \, dX  \leq \left(  \int_\Omega |\nabla u|^2 \varphi^4 G\, dX \right)^\frac12 \left( \int_\Omega |\varphi  u|^2 \1_{\supp \varphi} \frac{G}{\delta^2}\, dX \right)^\frac12 \\
\lesssim I^{1/2} \Big( \|\N^1(|\varphi u|^2)\|_{L^1(\partial \Omega,\omega)}\Big)^{\frac12} \leq (IJ_1)^{1/2},
\end{multline*}
by using, again, the fact that $\varphi$ is a cut-off function and the Carleson inequality found in Theorem \ref{ThCarlesonCO}.

It remains to estimate $I_1$. Since $u$ is a weak solution to $Lu=-\diver \vec F$, one has
\begin{multline*}
I_1= \int_\Omega \vec F \cdot \nabla (u\varphi^4 G)\, dX \\
=  \int_\Omega (\vec F \cdot \nabla u)\, \varphi^4 G\, dX + 4 \int_\Omega (\vec F \cdot \nabla \varphi) \, \varphi^3 uG \, dX + \int_\Omega (\vec F \cdot \nabla G)\, u\varphi^4\, dX \\
=:I_6+I_7+I_8.
\end{multline*}
By using \eqref{FubiniAtoG} and Lemma \ref{lemFubini},
\[|I_6| \leq I^\frac12 \left( \int_\Omega  |\vec F|^2 \varphi^4 G\, dX\right)^\frac12 \lesssim I^\frac12 \|\mathcal A_2(\delta\varphi^2\vec F)\|_{L^2(\partial \Omega,\omega)} \leq (IJ_2)^{1/2}.\]
For $I_7$, we invoke the $L^1$-variant of \eqref{FubiniAtoG} and the fact that $|\delta \nabla \varphi| \lesssim 1$ to get
\begin{multline*}
|I_7| \approx \int_\Omega  |\vec F| \varphi^2 u \frac{G}{\delta}\, dX  \lesssim  \int_{\pom} \mathcal A_1(\delta |\vec F| \varphi^2 u)\, d\omega
\leq \|\mathcal A_1^2(\delta\varphi\vec F) \|_{L^2(\partial \Omega,\omega)} \|\N^2(\varphi u)\|_{L^2(\partial \Omega,\omega)} \\ \leq (J_2J_1)^{\frac12} \leq J.
\end{multline*}
Finally, for $|I_8|$, by using  again the $L^1$-version of \eqref{FubiniAtoG},
\begin{multline*}
|I_8| \leq \int_\Omega  |\vec F||\nabla G| \varphi^4 u \, dX  \lesssim  \int_{\partial \Omega} \mathcal A_1\Big(\delta |\vec F| \frac{\delta|\nabla G|}{G} \varphi^4 u\Big)\, d\omega \\ \leq  \int_{\partial \Omega} \mathcal A_1^2(\delta |\vec F| \varphi) \N^2\Big(\varphi^2 u \frac{\delta |\nabla G|}{G}\Big) \, d\omega \leq \|\mathcal A_1^2(\delta |\vec F| \varphi)\|_{L^2(\partial \Omega, \omega)}  \Big\|\N^2\Big(\varphi^2 u \frac{\delta |\nabla G|}{G}\Big)\Big\|_{L^2(\partial \Omega,\omega)}  \\
\lesssim \|\mathcal A_1^2(\delta |\vec F| \varphi)\|_{L^2(\partial \Omega, \omega)}  \left( \|\N^2(\varphi u)\|_{L^2(\partial \Omega,\omega)} + \|\N^2(\varphi \delta |\vec F|)\|_{L^2(\partial \Omega,\omega)}\right) \\
\lesssim (J_2)^{1/2} (J_1+J_2)^{1/2} \leq J
\end{multline*}
where we used the double Caccioppoli inequality \eqref{DCaccio1} for the last but one line. We finished to prove that $I \lesssim J + (IJ)^{1/2}$, which concludes the theorem when $\Omega$ is unbounded or the support $\varphi$ avoids the center ball.

\medskip

To end the proof, consider that $\varphi$ is only in $C^\infty_c(\Omega)$. Then we take a function $\varphi_\Omega \in C^\infty(\Omega,[0,1])$ such that $\varphi_\Omega \equiv 0$ on $\Omega \setminus 4B_\Omega$, $\varphi_\Omega \equiv 1$ on $3B_\Omega$ and $|\nabla \varphi_\Omega| \leq 20/\delta(X_\Omega)$. Such function $\varphi_\Omega$ is a cut-off function with cutoff constants that depend only on the Corkscrew point condition. We have 
\[ \|\mathcal A(\varphi^2\delta\nabla u)\|_{L^{2}(\partial \Omega,\omega_L)} \lesssim  \|\mathcal A(\varphi^2\varphi^2_\Omega \delta\nabla u)\|_{L^{2}(\partial \Omega,\omega_L)} + \|\mathcal A(\varphi^2(1-\varphi_\Omega)^2 \delta\nabla u)\|_{L^{2}(\partial \Omega,\omega_L)} =: T_1 + T_2.\]
The term $T_2$ can be treated using the above argument, since $\supp [\varphi (1-\varphi_\Omega)] \cap 2B_\Omega = \emptyset$. As for the term $T_1$, we simply observe that $\varphi^2\varphi^2_\Omega \delta\nabla u$ is supported away from the boundary, and so
\begin{multline*}
    T_1 \approx \delta(X_\Omega) \left( \frac1{|4B_\Omega|} \int_{\Omega} \varphi^4 \varphi_\Omega^4 |\nabla u|^2 \, dX \right)^\frac12
\\ \lesssim \left( \fint_{4B_\Omega} \varphi^2 \varphi_\Omega^2  u^2 \, dX + \delta(X_\Omega) \fint_{4B_\Omega} \varphi^4 \varphi_\Omega^4 |\vec F|^2\, dX  \right)^\frac12 \\ 
\approx \|\N^2(\varphi \varphi_\Omega u)\|_{L^2(\partial \Omega,\omega)} + \|\A^2_1(\varphi^2 \varphi_\Omega^2 \delta \vec F)\|_{L^2(\partial \Omega,\omega)} \\ \lesssim \|\N^2(\varphi u)\|_{L^2(\partial \Omega,\omega)} + \|\A^2_1(\varphi \delta \vec F)\|_{L^2(\partial \Omega,\omega)},
\end{multline*} 
where the first inequality comes from the Caccioppoli inequality, and because
\[\dist(\supp \varphi_\Omega,\partial \omega) \approx \diam(B_\Omega) \approx \diam(\Omega).\] The theorem follows.
\ep

\subsection{$S<N$ for the Poisson-Regularity problem}

Before we move on, let us present first a nice theorem from \cite{HMMTZ21}, extending the works \cite{HM14,HMU14} from the Laplacian to DKP operators. 

\begin{theorem}[\cite{HMMTZ21}] \label{ThUR=>Ainfty}
Let $\Omega \subset \R^n$ be a 1-CAD. Let $L$ be a uniformly elliptic DKP operator. Then the following are equivalent:
\begin{enumerate}[(i)]
\item $\partial \Omega$ is uniformly rectifiable or equivalently that $\Omega$ is a Chord Arc Domain (by Proposition \ref{prUR<=>CAD}),
\item the elliptic measure with pole at center/infinity (on $\partial \Omega$, for the operator $L$) is $A_\infty$-absolutely continuous with respect to the Ahlfors regular measure of $\partial \Omega$.
\end{enumerate}
\end{theorem}

The main result that we prove in this section is the following.

\begin{theorem} \label{ThS<N2}
Let $\Omega \subset \R^n$ be a Chord Arc Domain, and let $B_\Omega$ be the center ball of $\Omega$ (when $\Omega$ is bounded) or the empty set (when $\Omega$ is unbounded). Let $\omega$ and $G$ be the harmonic measure and the Green function with pole at center/infinity. Let $L=-\diver A \nabla$ be a uniformly elliptic DKP operator that satisfies $(A^T-I) \nabla G = 0$ on $\Omega \setminus 2B_\Omega$.  

For $h\in L^{2}_{loc}(\Omega)$, $\vec F \in W^{1,2}_{loc}(\Omega,\R^n)$ satisfying $\vec F \cdot \nabla G = 0$ on $\Omega \setminus 2B_\Omega$, a weak solution $u$ to $Lu= h -\diver \vec F$ in $\Omega$, and a cut-off function $\varphi \in C^\infty_c(\Omega)$ (see Definition \ref{defcutoff}), we have
\begin{multline} \label{ThS<N2G}
 \|\mathcal A(\varphi^2\delta\nabla^2 u)\|_{L^{2}(\partial \Omega,\omega)} \leq C_\varphi \Big( \|\N^2(\varphi \nabla u)\|_{L^{2}(\partial \Omega,\omega)} + \|\mathcal A_1^2(\varphi\delta h) \|_{L^{2}(\partial \Omega,\omega)} \\ + \|\mathcal A(\varphi\delta \nabla \vec F) \|_{L^{2}(\partial \Omega,\omega)} + \|\N^2(\varphi\vec F) \|_{L^{2}(\partial \Omega,\omega)}  \Big).
 \end{multline}
The constant $C_\varphi$ depends only on the cut-off constants of $\varphi$, the ellipticity, boundedness, and DKP constant of $L$, and the CAD constants of $\Omega$.
\end{theorem}

\begin{remark}
Note that we treat $h$ and $\vec F$ differently: the bound on $\vec F$ in the right-hand side of \eqref{ThS<N2G} uses the functional $\A$, while the bound on $h$ use the functional $\A^2_1$ (which is larger than $\A$, see (ii) of Lemma \ref{lemcompNA}). 

We claim that if we were to replace $ \|\mathcal A(\varphi\delta \nabla \vec F) \|_{L^{2}(\partial \Omega,\omega)}$ by $ \|\mathcal A^2_1(\varphi\delta \nabla \vec F) \|_{L^{2}(\partial \Omega,\omega)}$, the theorem would be true without assuming $A^T\nabla G = \nabla G$ nor $\vec F \cdot \nabla G = 0$. 
This conjecture would give a new and potentially very interesting bound. However, the proof would be fairly technical and yet not useful for our present article, so we leave it for another project. 
\end{remark}

\bp The beginning of the proof is similar to the one of Theorem \ref{ThS<N1}, and the end will use the key integration by part established in Lemma \ref{IPPwithNG}. 

First, as presented at the end of the proof of Theorem \ref{ThS<N1}, we only need to consider cut-off functions that avoids the center ball, i.e. $\varphi \in C^\infty_0(\Omega \setminus 2B_\Omega)$. If $\varphi \in C^\infty_0(\Omega)$ , then we construct $\varphi_\Omega$ as in the proof of Theorem \ref{ThS<N1} and we have
\[  \|\mathcal A(\varphi^2\delta\nabla^2 u)\|_{L^{2}(\partial \Omega,\omega)} \lesssim  \|\mathcal A(\varphi^2\varphi_\Omega^2 \delta\nabla^2 u)\|_{L^{2}(\partial \Omega,\omega)} +  \|\mathcal A(\varphi^2(1-\varphi_\Omega)^2\delta\nabla^2 u)\|_{L^{2}(\partial \Omega,\omega)} =: T_1 + T_2.\]
The function $\varphi(1-\varphi_\Omega) \in C^\infty_0(\Omega \setminus 2B_\Omega)$ is a cut-off function that avoids the center ball, so $T_2$ can be treated with the main argument below. As for $T_1$, we simply observe again that 
\begin{multline*}
    T_1 \approx \delta(X_\Omega) \left( \frac1{|4B_\Omega|} \int_{\Omega} \varphi^4 \varphi_\Omega^4 |\nabla^2 u|^2 \, dX \right)^\frac12 \\
    \lesssim \left( \fint_{4B_\Omega} \varphi^2 \varphi_\Omega^2  |\nabla u|^2 \, dX + \delta(X_\Omega) \fint_{4B_\Omega} \varphi^4\varphi_\Omega^4 [h^2 + |\nabla \vec F|^2] \, dX  \right)^\frac12 \\ \approx \|\N^2(\varphi \varphi_\Omega \nabla u)\|_{L^2(\partial \Omega,\omega)} + \|\A^2_1(\delta \varphi^2 \varphi_\Omega^2 h )\|_{L^2(\partial \Omega,\omega)} + \|\A(\varphi^2 \varphi_\Omega^2 \delta \nabla \vec F)\|_{L^2(\partial \Omega,\omega)} \\
    \lesssim \|\N^2(\varphi \nabla u)\|_{L^2(\partial \Omega,\omega)} + \|\A^2_1(\varphi \delta h)\|_{L^2(\partial \Omega,\omega)} + \|\A(\varphi \delta \vec F)\|_{L^2(\partial \Omega,\omega)},
\end{multline*} 
where the second inequality comes from Lemma \ref{lemCaccioN2u}, since $\varphi \varphi_\Omega \in C^\infty_0(2B_{X_\Omega})$\footnote{We recall that $B_{X_\Omega} = 4B_\Omega$} is a cut-off function. 


\medskip

We move to the main argument of the theorem, and we assume now that $\varphi \in C^\infty_0(\Omega \setminus 2B_\Omega)$. Let $\partial_i$ be the directional derivatives in $\R^n$. We have, with \eqref{FubiniAtoG} and the ellipticity of $A$, that 
\[ \|\mathcal A(\varphi^2\delta\nabla^2 u)\|_{L^{2}(\partial \Omega,\omega_L)}^2  \lesssim I:= \sum_{i=1}^n \int_\Omega (A\nabla \partial_i u \cdot \nabla \partial_i u) \, \varphi^4 G\, dX.\]
To lighten the notation, we write $\partial$ for any of the $\partial_i$, $i\in \{1,\dots,n\}$. Similarly to the proof of Theorem \ref{ThS<N1}, we want to prove that for any partial derivative $\partial$,
\[ I':= \int_\Omega (A\nabla \partial u \cdot \nabla \partial u) \, \varphi^4 G\, dX \lesssim J + (IJ)^\frac12,\]
where
\[J:= \underbrace{\|\N^2(\varphi \nabla u)\|_{L^{2}(\partial \Omega,\omega)}^2}_{J_1} + \underbrace{\|\mathcal A_1^2(\varphi\delta h) \|_{L^{2}(\partial \Omega,\omega)}^2}_{J_2} + \underbrace{\|\mathcal A(\varphi\delta \nabla \vec F) \|_{L^{2}(\partial \Omega,\omega)}^2}_{J_3} + \underbrace{\|\N^2(\varphi\vec F) \|_{L^{2}(\partial \Omega,\omega)}^2}_{J_4}.\]
From the use of a compactly supported function $\varphi$ and from Lemma \ref{lembdonG}, the quantity $I$ is finite and so the above estimate self-improves to $I\lesssim J$, as desired for the Theorem.

\medskip

The proof is similar to the one of Theorem \ref{ThS<N1}, at least at the beginning. We will make a parallel between the terms $I_k$ from this proof and the ones of Theorem \ref{ThS<N1}, in the sense that if the terms have the same notation in both proofs, then they will be treated with a similar argument.
\begin{multline*}
I' = \int_\Omega A \nabla \partial u \cdot \nabla (\partial u \, \varphi^4 G) \, dX - 4 \int_\Omega (A \nabla \partial u \cdot \nabla \varphi )\, \varphi^3 G \partial u \, dX - \int_\Omega (A\nabla \partial u\cdot \nabla G) \, \varphi^4 \partial u \, dX \\ =: I_1 + I_2 + I_3,
\end{multline*}
and we further rewrite $I_3$ as
\[
I_3 = -\frac12 \int_\Omega A\nabla (\varphi^4 [\partial u]^2) \cdot \nabla G \, dX + 2 \int_\Omega (A \nabla \varphi \cdot \nabla G)\, \varphi^3 [\partial u]^2\, dX =: I_4 + I_5.
\]
The term $I_4$ is 0, since $A^T \nabla G = \nabla G$ and $G$ is a weak solution to $-\Delta G=0$ on $\supp \varphi$, and $\varphi^4 [\partial u]^2$ is a $W^{1,2}_{loc}(\Omega)$ function\footnote{This a a simple consequence of \eqref{uW22}} that is compactly supported in $\Omega$. For the term $I_5$, we  use the bounds $|\nabla G| \lesssim G/\delta$ (lemma \ref{lembdonG}), and $\delta |\nabla \varphi|\lesssim 1$ (recall that $\varphi$ is a cut-off function), and we have
\[
|I_5| \lesssim \int_\Omega (\1_{\supp \nabla \varphi}) \varphi^2 |\nabla u|^2 \frac{G}{\delta^2} \, dX \leq  C_\varphi \|\N^1(|\varphi^2 \nabla u|^2)\|_{L^1(\partial \Omega,\omega)} \lesssim J_1
\]
Where we use the fact that $\varphi$ is a cut-off function and the Carleson inequality (Theorem \ref{ThCarlesonCO}) for the second inequality. We move to $I_2$, that we bound as follows:
\[
|I_2| \lesssim \int_\Omega |\nabla \partial u| |\partial u| (\1_{\supp \nabla \varphi}) \varphi^3 \frac{G}{\delta} \, dX 
\lesssim I^{1/2} \|\N^1(|\varphi^2 \nabla u|^2)\|_{L^1(\partial \Omega,\omega)}^{\frac12} \lesssim (IJ_1)^{1/2},\]
by using, again, the Carleson inequality found in Theorem \ref{ThCarlesonCO}.

It remains to bound $I_1$, where all the difficulty lies. To treat this term, we shall first use the fact that $u$ is a weak solution to $Lu= h -\diver \vec F$ in $\Omega$, and thus that $v=\partial u$ is a weak solution to $Lv= \partial h + \diver [(\partial A) \nabla u - \partial \vec F]$. The latter equation is a weak formulation, and in particular $h$ does not have to be differentiable: once we know that $u \in W^{2,2}_{loc}(\Omega)$ --- see Lemma \ref{lemCaccioN2u} --- it is easy to check with few integrations by parts that the relation 
\[ \int_\Omega A \nabla \partial u \cdot \nabla w \, dX =  \int_\Omega\Big( [\partial \vec F - (\partial A) \nabla u]\cdot \nabla w - h\partial w \Big)\, dX  \]
holds whenever $w$ is in $C^\infty_0(\Omega)$. And since $\partial u \in W^{1,2}_{loc}(\Omega)$, the identity extends to all function $w \in W^{1,2}_{loc}(\Omega)$ that are compactly supported in $\Omega$. Using this fact, we get 
\begin{multline*}
I_1 
= - \int_\Omega h \partial (\varphi^4 G \partial u )\, dX + \int_\Omega [\partial \vec F - (\partial A) \nabla u] \cdot \nabla (\partial u \, \varphi^4 G) \, dX \\
= -\int_\Omega h (\partial^2 u) \, \varphi^4 G \, dX - 4 \int_\Omega h \, (\partial u) \, (\partial \varphi) \, \varphi^3 G \, dX -  \int_\Omega h \varphi^4 (\partial G) \, (\partial u) \, dX \\
+ \int_\Omega [\partial \vec F - (\partial A) \nabla u] \cdot [\varphi \nabla[\partial u] + 4(\partial u)\nabla \varphi]   \, (\varphi^3 G) \, dX \\
+  \int_\Omega [\partial \vec F - (\partial A) \nabla u] \cdot \nabla G \, (\varphi^4 \partial u) \, dX
=: I_6+I_7+I_8+I_9+I_{10}.
\end{multline*}
The terms $I_6$, $I_7$, and $I_8$ are simpler analogues to the ones in the proof of Theorem \ref{ThS<N1}. They satisfy
\[ |I_6| \lesssim \|\A(\delta\varphi h)\|_{L^2(\partial\Omega,\omega)} \| \A(\delta \varphi \nabla^2 u)\|_{L^2(\partial\Omega,\omega)} \lesssim (IJ_2)^\frac12\]
and, since $\varphi$ is a cut-off function and \eqref{uW1infty} holds, we have $\delta G |\nabla \varphi| + \delta \varphi |\partial G| \lesssim G$, that is
\[ |I_7|+|I_8| \lesssim \int_\Omega |\delta h| |\nabla u| \varphi^3 \frac{G}{\delta^2}\, dX \lesssim \|\A_1^2(\delta\varphi h)\|_{L^2(\partial\Omega,\omega)} \| \N^2(\varphi \nabla u)\|_{L^2(\partial\Omega,\omega)} \lesssim (J_2J_1)^\frac12 \lesssim J.\]
Theorem \ref{ThUR=>Ainfty} gives that $\omega \in A_\infty(\sigma)$, and so we can use Theorem \ref{ThCarleson1}. From there, we treat $I_9$ without new subtilities, using the Carleson inequalities (Theorems \ref{ThCarleson1} and \ref{ThCarlesonCO}) when needed:
\begin{multline*}
|I_9| \lesssim \int_\Omega \delta\big[|\nabla \vec F| + |\nabla A||\nabla u|\big] \, \big[ \delta |\nabla^2 u| + \1_{\supp\varphi}  |\nabla u|\big] \varphi^3 G\, dX \\
\lesssim \Big\|\A\big(\delta\varphi\big[|\nabla \vec F| + |\nabla A||\nabla u|\big]\Big)\Big\|_{L^2(\partial\Omega,\omega)} \Big\|\A\Big( \varphi\big[ \delta|\nabla^2 u| + \1_{\supp\nabla \varphi}  |\nabla u|\big]\Big)\Big\|_{L^2(\partial\Omega,\omega)} \\
\lesssim \left( \big\|\A\big(\delta \varphi \nabla \vec F\big)\big\|_{L^2(\partial\Omega,\omega)} + \|\delta \nabla A\|_{CM_2^\infty} \big\|\N^2\big(\varphi\nabla u\big)\big\|_{L^2(\partial\Omega,\omega)} \right) \\ \hspace{2cm}\left( \big\|\A\big(\delta \varphi\nabla^2 u\big)\big\|_{L^2(\partial\Omega,\omega)} + 
 \big\|\N^2\big(\varphi\nabla u\big)\big\|_{L^2(\partial\Omega,\omega)} \right)   \\
\lesssim (J_3 + J_1)^{1/2}(I+J_1)^{1/2} \lesssim J + (IJ)^{1/2}.
\end{multline*}
 
 \medskip
 
 The real novelty and difficulty of the proof is the treatment  $I_{10}$. We use the fact that $A^T \nabla G = \nabla G$ and $\vec F \cdot \nabla G =0$ to write
 \[ [\partial \vec F - (\partial A)\nabla u ] \cdot \nabla G = - \vec F \cdot \nabla \partial G - \nabla u \cdot [\partial (A^T \nabla G) - A^T \nabla \partial G] = [- \vec F + (A-I)\nabla u] \cdot \nabla \partial G,\]
 and so 
 \[I_{10} = \int_\Omega [ - \vec F + (A-I)\nabla u] \cdot \nabla[\partial  G] \, \varphi^4 \, \partial u \, dX \]
Note that with our assumption, $(\partial u) [\vec F - (A-I)\nabla u]\cdot \nabla G = 0$.  So we can apply Lemma \ref{IPPwithNG} with $v=\partial G$, and $\phi = \varphi^4$. We obtain then that 
 \begin{multline*}
 |I_{10}| \lesssim \int_{\Omega} \delta|\nabla^2 G| \varphi^3 \Big( |\vec F| |\nabla \varphi| + |\nabla^2 u| \big[ |\vec F| + |A-I| |\nabla u| \big] +  |\nabla u| \big[ |\nabla \vec F| + |\nabla A| |\nabla u| + |A-I| |\nabla u| \Big) \, dX  \\
 + \int_\Omega \Big(\frac{\delta^2 |\nabla^2 G|}{G} + \1_{\{\delta|\nabla G|< cG\}} \Big) |\nabla u| \big[ |\vec F| + |A-I| |\nabla u|\big] |\nabla^2 G| \varphi^4 \, dX. 
 \end{multline*}
 We can rearrange (gather the same terms), remove $|A-I|$ due to the boundedness of $A$, add $\1_{\Omega \setminus 2B_\Omega}$ for free because $\varphi \in C^\infty_0(\Omega \setminus 2B_\Omega)$, and we get
 \begin{multline*}
   |I_{10}| \lesssim \int_{\Omega} \delta \Big[ |\nabla^2 u| + |\nabla \vec F|\Big] \,  \Big[ |\nabla u| + |\vec F|\Big]  \frac{\delta^2 |\nabla^2 G|}{G} \varphi^3  \frac{G}{\delta^2}\, dX \\
 + \int_\Omega \Big(|\delta \nabla A| + \frac{\delta^2 |\nabla^2 G|}{G}\1_{\Omega \setminus 2B_\Omega} + \1_{\{\delta|\nabla G|< cG\}} \Big)^2  \Big[ |\nabla u| + |\vec F|\Big] ^2\varphi^3\, \frac{G}{\delta^2} \, dX \\
 =:I_{11} + I_{12}. 
 \end{multline*}
 
The term $I_{11}$ is bounded, using the Cauchy-Schwarz inequality and then the Carleson inequality provided by Theorem \ref{ThCarleson1} (which we can since $\omega \in A_\infty(\sigma)$ by Theorem \ref{ThUR=>Ainfty}), as follows
\begin{multline*}
I_{11} \lesssim \Big\| \mathcal A\Big( \delta\varphi^2  \Big[ |\nabla^2 u| + |\nabla \vec F|\Big]  \Big)\Big\|_{L^2(\partial \Omega,\omega)} \Big\| \mathcal A\Big( \big[ |\nabla u| + |\vec F|\big]  \frac{\delta^2 |\nabla^2 G|}{G} \varphi  \Big)\Big\|_{L^2(\partial \Omega,\omega)} \\
\lesssim (I+J_3)^\frac12 \Big\|\frac{\delta^2 |\nabla^2 G|}{G} \Big\|_{CM_2^\infty} \Big\| \mathcal N^2\Big( \big[ |\nabla u| + |\vec F|\big] \varphi  \Big)\Big\|_{L^2(\partial \Omega,\omega)} \lesssim (I+J_3)^\frac12 (J_1+J_4)^\frac12 \lesssim J + (IJ)^{1/2}.
\end{multline*}
thanks to Theorem \ref{ThonG}, which guarantees that $\delta^2\nabla^2 G/G \in CM_2^\infty$ is finite. Finally, the fact that $\delta |\nabla A| + \frac{\delta^2  |\nabla^2 G| }{G} \1_{\Omega \setminus 2B_\Omega}  + \1_{\{\delta|\nabla G|< cG\}} \in CM_2^\infty$ --- which is a consequence of the fact that $L$ is a DKP operator, Theorems \ref{ThonG} and \ref{CorNG>0} ---  is used to get that  
\[\Big(|\delta \nabla A| + \frac{\delta^2 |\nabla^2 G|}{G}\1_{\Omega \setminus 2B_\Omega} + \1_{\{\delta|\nabla G|< cG\}} \Big)^2 \in CM_1^\infty,\]
hence that 
\[I_{12} \lesssim\Big\| \mathcal N^1\Big( \big[ |\nabla u| + |\vec F|\big]^2 \varphi^2  \Big)\Big\|_{L^1(\partial \Omega,\omega)} \lesssim  J_1 + J_4 \lesssim J\]
by using a last time the Carleson inequality (Theorem \ref{ThCarleson1}). The Theorem follows.
\ep

\subsection{Good-$\lambda$ argument}

We want to put Theorems \ref{ThS<N1} and \ref{ThS<N2} through the good-$\lambda$ machinery. This allows us to change the $L^2(\partial \Omega,\omega)$ norm in the conclusion of Theorems \ref{ThS<N1} and \ref{ThS<N2} by a $L^p(\partial \Omega,\sigma)$ norm, as long as $p\in (0,\infty)$ and $\omega \in A_\infty(\sigma)$. The good-lambda argument is pretty classical in the area, starting with \cite[Lemma 8]{Dahlberg80}. 
we contemplated to not write it in the present article. However, since it is not classical when the boundary is neither flat nor Lipschitz, or when cut-off functions are used instead of saw-tooth domains\footnote{an alternative strategy preferred by the authors, as in \cite{FMZ21} and \cite{DFMreg}, which do not need to rely on Lipschitz sub-domains.}, we decided to provide a detailed proof.

\begin{theorem} \label{Thgoodlambda}
Let $\Omega$ be a 1-sided CAD. Let $\omega$ and $\sigma$ be two measures on $\partial \Omega$ such that $\omega \in A_\infty(\sigma)$. Assume moreover that, for $N\in \bN^*$, we have some function $f$, $(g_1,\dots,g_n)$ on $\Omega$ and some parameters $(\alpha_1,\dots,\alpha_N)\in [1,\infty]^N$, $(\beta_1,\dots,\beta_N) \in [1,\infty]$ such that 
\begin{enumerate}[(i)]
\item for any $p\in (0,+\infty)$,
\begin{equation} \label{S<Nloc6}
\| \A(f)\|_{L^p(\partial \Omega, \sigma)} + \sum_{k=1}^N \|\A_{\alpha_k}^{\beta_k}(g_k)\|_{L^p(\partial \Omega, \sigma)}< +\infty,
\end{equation}
and
\item
\begin{equation} \label{S<Nloc5}
\| \A(f \phi^2)\|_{L^2(\partial \Omega, \omega)} \leq C_\phi \sum_{k=1}^N \|\A_{\alpha_k}^{\beta_k}(g_k\phi)\|_{L^2(\partial \Omega, \omega)} ,
\end{equation}
whenever $\phi$ is a cut-off function and where $C_\phi$ depends only on $\phi$ via its cut-off constants. In particular we do not need $\phi$ to be compactly supported.
\end{enumerate} 
Then for any $p\in (0,\infty)$, there exists $C_p>0$ such that
\begin{equation} \label{S<Nloc7}
\| \A(f)\|_{L^p(\partial \Omega, \sigma)} \leq C_p \sum_{k=1}^N \|\A_{\alpha_k}^{\beta_k}(g_k)\|_{L^p(\partial \Omega, \sigma)} < + \infty,
\end{equation}
where $C_p$ depends only on the constants in the definition of $\omega \in A_\infty(\sigma)$, the 1-sided CAD constant of $\Omega$, and an upper bound on $C_\phi$ for a class cut-off functions with a uniformly bound on the cut-off constants that depends (i.e. the uniform bound depends) on the 1-sided CAD constant of $\Omega$.
\end{theorem}

\begin{remark}
Note that $\A_\infty^\beta = \N^\beta$, so the theorem also apply when we have non-tangential maximal function on the right-hand side.
\end{remark}

\bp The proof of \eqref{S<Nloc6} $+$ \eqref{S<Nloc5} $\implies $ \eqref{S<Nloc7} will be decomposed in 2 parts: The first and main part is \eqref{S<Nloc5} implies that for any $q>0$, there exists $\eta \in (0,1)$, $\alpha>0$ and $C>0$ such that, for any $\gamma\in (0,1)$ and $\lambda >0$, we have 
\begin{multline}
 \label{S<Nloc8}
\sigma\Big(\Big\{x\in \partial \Omega, \,  \mathcal A(f)(x) >\lambda,\, \mathcal M_q\Big[\sum_{k=1}^N \mathcal A_{\alpha_k}^{\beta_k}(g_k)\Big](x) \leq \gamma\lambda\Big\}\Big) \\
\leq C \gamma^\alpha \sigma(\{x\in \partial \Omega, \mathcal M_q(\mathcal A(f))(x) > \eta \lambda\}),
\end{multline}
where $\mathcal M_q(f)(x) = \sup_{B \ni x}\left( \int_{B} |f|^q d\sigma \right)^\frac1q$ is the $L^q$-version of the Hardy-Littlewood maximal function (i.e. the supremum is taken over all the boundary balls $B\subset \partial \Omega$ that contains $x$). The second part is that \eqref{S<Nloc6} $+$ \eqref{S<Nloc8} $\implies $ \eqref{S<Nloc7}, which is actually a simple application of Cavalieri's theorem. Indeed,  take $p\in (0,\infty)$ and then $q\in (0,p)$. By \eqref{S<Nloc8}, 
\begin{multline*}
\| \A(f)\|_{L^p(\partial \Omega, \sigma)}^p = p \int_0^\infty \lambda^{p-1} \sigma(\{\mathcal A(f)>\lambda\}) \, d\lambda \\
 \leq p \int_0^\infty \lambda^{p-1} \sigma\Big(\Big\{ \mathcal A(f) >\lambda,\, \mathcal M_q\Big[\sum_{k=1}^N \mathcal A_{\alpha_k}^{\beta_k}(g_k)\Big] \leq \gamma\lambda\Big\}\Big) \, d\lambda \\ + p \int_0^\infty \lambda^{p-1} \sigma\Big(\mathcal M_q\Big[\sum_{k=1}^N \mathcal A_{\alpha_k}^{\beta_k}(g_k)\Big] > \gamma\lambda\Big\}\Big) \, d\lambda \\
 \leq C \gamma^\alpha  \int_0^\infty \lambda^{p-1} \sigma(\{\mathcal M_q(\mathcal A(f))>\eta \lambda \}) \, d\lambda + p \int_0^\infty \lambda^{p-1} \sigma\Big(\mathcal M_q\Big[\sum_{k=1}^N \mathcal A_{\alpha_k}^{\beta_k}(g_k)\Big] > \gamma\lambda\Big\}\Big) \, d\lambda \\
 \leq C_{\eta} \gamma^\alpha \| \mathcal M_q(\A(f))\|_{L^p(\partial \Omega, \sigma)}^p + C_{\gamma} \Big\| \mathcal M_q\Big[ \sum_{k=1}^N \mathcal A_{\alpha_k}^{\beta_k}(g_k)\Big]\Big\|_{L^p(\partial \Omega, \sigma)}^p \\
 \leq C_{\eta,p-q} \gamma^\alpha \|\A(f)\|_{L^p(\partial \Omega, \sigma)}^p + C_{\gamma,p-q} \sum_{k=1}^N \|\A_{\alpha_k}^{\beta_k}(g_k)\|_{L^p(\partial \Omega, \sigma)}^p,
\end{multline*}
where the last line is due to the $L^{p/q}$-boundedness of the Hardy-Littlewood maximal function (or the $L^p$-boundedness of $\mathcal M_q$). We choose $\gamma$ small enough such that $C_{\eta,p-q} \gamma^\alpha \leq \frac12$, and together the finiteness of $\|\A(f)\|_{L^p(\partial \Omega, \sigma)}$ given by \eqref{S<Nloc6}, we obtain \eqref{S<Nloc7} by hidding $C_{\eta,p-q} \gamma^\alpha \|\A(f)\|_{L^p(\partial \Omega, \sigma)}^p$ to the left-hand side.

\medskip

It remains to prove that \eqref{S<Nloc5} $\implies$ \eqref{S<Nloc8}. We write 
\[ F_\gamma:=\Big\{x\in \partial \Omega, \,  \mathcal A(f)(x) >\lambda,\, \mathcal M_q\Big[\sum_{k=1}^N \mathcal A_{\alpha_k}^{\beta_k}(g_k)\Big](x) \leq \gamma\lambda\Big\}\]
and
\[ E:= \{x\in \partial \Omega, \mathcal M_q(\mathcal A(f))(x) > \eta \lambda\}.\]
Let $\{B_i\}_{i\in I}$ be a finitely overlapping Whitney decomposition of $E$ by boundary balls. We can obtain it easily from the Vitali covering lemma: the collection of boundary balls $\{B(x,\dist(x,\partial \Omega \setminus E)/10) \cap E\}_{x\in E}$ is a covering of $E$, so the the Vitali covering lemma gives a maximal subcollection of non-overlapping balls $\{B(x_i,\dist(x_i,\partial \Omega \setminus E)/10) \cap E\}_{i\in I}$. The maximality means that $\{B_i\}_{i\in E}:= \{B(x_i,\dist(x_i,\partial \Omega \setminus E)/2) \cap E\}_{i\in I}$ covers $E$. The non-overlapping of $B_i/5$ implies the finitely overlapping of the $B_i$. 

In order to prove \eqref{S<Nloc8}, which is $\sigma(F_\gamma) \leq C\gamma^\alpha \sigma(E)$, it is sufficient to prove that
\begin{equation} \label{S<Nloc9}
\omega( F_\gamma \cap B_i) \leq C\gamma^2 \omega(B_i),
\end{equation}
where $\omega$ is the harmonic measure with pole at center/infinity. Indeed, from the fact that $\omega \in A_\infty(\sigma)$ --- see Definition \ref{defAinfty} --- we have that \eqref{S<Nloc9} gives
\begin{equation} \label{S<Nloc9b}
\sigma( F_\gamma \cap B_i) \leq C'\gamma^\alpha \sigma(B_i),
\end{equation}
for an $\alpha \in (0,2)$. Then we conclude by saying that 
\[\sigma (F_\gamma) \leq \sum_{i\in I} \sigma( F_\gamma \cap B_i) \leq C'\gamma^\alpha \sum_{i\in I} \sigma(B_i) \leq C'' \gamma^\alpha \sigma(E)\]
since $B_i$ is a finitely overlapping covering of $E$ and $F_\gamma$.

So we need to prove \eqref{S<Nloc9} from \eqref{S<Nloc5}. We notice first that there is nothing to prove if $ F_\gamma \cap B_i = \emptyset$. Moreover, since $B_i$ is a Whitney ball, there is a $z_i\in \partial \Omega \setminus E$ such that $\dist(z_i,B_i)= r_i$, the radius of $B_i$. 
When $F_\gamma \cap B_i$ is non-empty, we use ``saw-tooth domains'', which in our case will be done via cut-off functions. We define $\zeta: \, \R\to [0,1]$ such that $\zeta(t) = 0$ if $t\leq 1$, $\zeta(t) = 1$ if $t\geq 2$, and $|\zeta'|\leq 2$. Then we construct $\phi$ on $\Omega$ as
\[ \phi(X) = \zeta\Big( \frac{4\delta(X)}{\dist(X,F_\gamma \cap B_i)} \Big) \Big[ 1 - \zeta\Big( \frac{\delta(X)}{r_i}\Big)\Big] ,\]
where $r_i$ is the radius of $B_i$. The function $\phi$ is really a cut-off function, that is $\delta \nabla \phi \in L^\infty \cap CM^\infty_1$ (and the constants are independent of $\gamma$ or $i$) thanks to Lemma \ref{lemcutoff}.

We want to say that we do not lose much by taking $\mathcal A(f\phi^2)(x)$ instead of $\mathcal A(f)(x)$ when $x\in F_\gamma \cap B_i$. Indeed, if $X\in \Gamma(x)$ and $x\in F_\gamma \cap B_i$, then $\dist(X,F_\gamma\cap B_i) \leq \dist(X,x) \leq 2 \dist(X,\partial \Omega)$, so 
\[\zeta\Big( \frac{4\dist(X,\partial \Omega)}{\dist(X,F_\gamma \cap B_i)} \Big)  =1.\]
So when $x\in F_\gamma \cap B_i$, the difference between $\mathcal A(f\phi^2)(x)$ and $\mathcal A(f)(x)$ comes from the second truncation, i.e. 
\[\mathcal A(f)(x) - \mathcal A(f\phi^2)(x) \leq \mathcal A(f \1_{\delta\geq r_i})(x).\]
Now, we claim that 
\begin{equation} \label{S<Nloc10}
\mathcal A(f \1_{\{\delta\geq r_i\}})(x) \leq C \mathcal M_q(\mathcal A(f))(z),
\end{equation}
which, since $z\in \partial \Omega \setminus E$, gives that $A(f \1_{\delta\geq r_i})(x)\leq C\eta\lambda$ and then 
\begin{equation} \label{S<Nloc11}
\mathcal A(f\phi^2)(x) \geq \lambda/2 \qquad \text{ for } x\in F_\gamma \cap B_i
\end{equation}
 if $\eta$ is chosen (once for all the proof) small enough. To show the claim \eqref{S<Nloc10}, 
we introduce $\Gamma(a,r):= \{X\in \Omega, |X-a| \leq 2 \delta(X) - r\}$, we observe that 
\begin{equation} \label{randominclA} \Gamma(a,r_i) \subset \Gamma (b) \qquad \text{ whenever $a,b\in \partial \Omega$ are such that $|a-b|\leq r_i$}.\end{equation}
If $\{B(x_j,r_i/2) \cap \partial \Omega\}_{j\in J_{x,i}}$ is a covering of $B(x,2r_i) \cap \partial \Omega$ with $x_j\in B(x,2r_i) \cap \partial \Omega$ and $\# J_{x,i}$ to be uniformly finite, which is possible thanks to Vitali's lemma, we have then that 
\begin{multline} \label{inclusioncone}
\Gamma(x)\cap \{\delta \geq r_i\} \subset \bigcup_{a\in B(x,2r_i) \cap \partial \Omega} \Gamma(a,r_i) \subset \bigcup_{j\in J_{x,i}}  \bigcup_{a\in B(x_j,r_i/2) \cap \partial \Omega} \Gamma(a,r_i) \\
\subset \bigcup_{j\in J_{x,i}}  \bigcap_{b\in B(x_j,r_i/2) \cap \partial \Omega} \Gamma(b),
\end{multline}
by \eqref{randominclA}, hence
\[ \mathcal A(f\1_{\{\delta\geq r_i\}}) \lesssim \sum_{j\in J_{x,i}} \fint_{B(x_j,r_i/2)} |\mathcal A(f)|^q \, d\sigma \lesssim \fint_{B(x_j,5r_i)} |\mathcal A(f)|^q \, d\sigma \leq \mathcal M_q(\mathcal A(f))(z)\]
since $|z-x_i| < 5r_i$ and $\# J_{x,i}$ is uniformly finite.

With \eqref{S<Nloc11} in hand, we can estimate $F_\gamma \cap B_i$ by
\begin{multline} \label{bdonomegaF}
\omega(F_\gamma \cap B_i) \leq \frac4 {\lambda^2} \int_{F_\gamma \cap B_i} |\mathcal A(f\phi^2)(x)|^2\, d\omega  \leq  \frac4 {\lambda^2} \|\mathcal A(f\phi^2)\|_{L^2(\partial \Omega,\omega)}^2 \\ \lesssim  \frac4 {\lambda^2} \Big\|\sum_{k=1}^N \mathcal A_{\alpha_k}^{\beta_k}(\phi g_k)\Big\|_{L^2(\partial \Omega,\omega)}^2
\end{multline}
by \eqref{S<Nloc5}. 

In one hand, for $\sigma$-a.e. $x\in F_\gamma \cap B_i$, 
\[ \sum_{k=1}^N \A_{\alpha_k}^{\beta_k}(\phi g_k)(x) = \sum_{k=1}^N \A_{\alpha_k}^{\beta_k}(g_k) \leq \mathcal M_q\Big[ \sum_{k=1}^N \A_{\alpha_k}^{\beta_k}(g_k)\Big](x).\]
In the other hand if $F_\gamma \cap B_i$ is not empty, for $x\in \partial \Omega \cap B_i$, we have $y\in \overline{F_{\gamma} \cap B_i}$ such that $|x-y| = \dist(x,F_{\gamma} \cap B_i)$, and then $\phi = 0$ on $B_X$ for any $X\in \Gamma(x) \cap \{\delta \leq |x-y|/7\}$. But similarly to \eqref{inclusioncone}, we have
\[ \Gamma(x) \cap \{\delta \geq |x-y|/7\}  \subset \bigcup_{j\in J_{x,i,\gamma}} \bigcap_{b\in B(x_j,|x-y|/12) \cap \partial \Omega} \Gamma(b)\]
with $\# J_{x,i,\gamma}$ uniformly finite and $x_j \in B(x, |x-y|/3) \partial \Omega$. So in particular
\begin{multline*}
\sum_{k=1}^N \A_{\alpha_k}^{\beta_k}(\phi g_k)(x)  \leq \sum_{j\in J_{x,i,\gamma}} \left( \fint_{B(x_j,|x-y|/14) \cap \partial \Omega} \left|\sum_{k=1}^N \A_{\alpha_k}^{\beta_k}(\phi g_k)(x) \right|^q \, d\sigma   \right)^\frac1q \\ \lesssim \mathcal M_q\Big[ \sum_{k=1}^N \A_{\alpha_k}^{\beta_k}(g_k)\Big](y) \leq \gamma \lambda.\
\end{multline*}
Altogether, we have
\[ \Big\|\sum_{k=1}^N \mathcal A_{\alpha_k}^{\beta_k}(\phi g_k)\Big\|_{L^2(\partial \Omega,\omega)}^2 \lesssim \omega\left(\supp \Big[\sum_{k=1}^N \mathcal A_{\alpha_k}^{\beta_k}(\phi g_k)\Big]\right) (\gamma \lambda)^2 \lesssim \omega(B_i) (\gamma \lambda)^2  \]
since, by construction of $\phi$, we have $\supp \mathcal A_{\alpha_k}^{\beta_k}(\phi g_k) \subset 10B_i$ and $\omega$ is doubling. Reinserting this bound in \eqref{bdonomegaF}, we get
\[\omega(F_\gamma\cap B_i) \lesssim \gamma^2 \omega(B_i).\]
The claims \eqref{S<Nloc9}, \eqref{S<Nloc8} and \eqref{S<Nloc7} follow, so is the Theorem.
\ep

\subsection{Conclusion on $S<N$ for the Poisson problem}

Theorem \ref{Thgoodlambda} improves Theorems \ref{ThS<N1} and \ref{ThS<N2} in the following manner.

\begin{theorem} \label{ThS<N11} 
Let $\Omega \subset \R^n$ be a 1-sided Chord Arc Domain, where $\sigma$ is the Ahlfors-regular measure on $\partial \Omega$. Let $L:=-\diver A \nabla$ be a uniformly elliptic operator, and assume that  $\omega_L \in A_\infty(\sigma)$ where $\omega_L$ is the elliptic measure for $L$  on $\partial \Omega$ with pole at center/infinity.

For any $p\in (0,\infty)$, there exists $C_p>0$ such that for any $\vec F \in L^2_{loc}(\Omega,\R^n)$, any weak solution $u$ to $Lu=-\diver \vec F$ in $\Omega$, we have
\[
\|\A(\delta \nabla u)\|_{L^{p}(\partial \Omega,\sigma)} \leq C_p \Big( \|\N^2(u)\|_{L^{p}(\partial \Omega,\sigma)}  + \|\A^2_1(\delta \vec F)\|_{L^{p}(\partial \Omega,\sigma)} \Big).\]
The constant $C_p$ depends only on $p$, the ellipticity and boundedness constant of $L$ and the 1-sided CAD constants of $\Omega$. In any case, the finiteness of the quantities involved are not guaranteed. However, the left-hand side will be finite as long as the right-hand side is.
\end{theorem}

\begin{theorem} \label{ThS<N21}
Let $\Omega \subset \R^n$ be a Chord Arc Domain. Let $G$ be the Green function with pole at center/infinity associated to the Laplacian. Take $L=-\diver A \nabla$ be a DKP operator that satisfies $A^T\nabla G = \nabla G$ a.e. on $\Omega$. 

For any $p\in (0,\infty)$, there exists $C_p>0$ such that for any $h\in L^2_{loc}(\Omega)$, any $\vec F \in W^{1,2}_{loc}(\Omega,\R^n)$ satisfying $\vec F \cdot \nabla G = 0$ a.e. on $\Omega$, and any weak solution $u$ to $Lu=-\diver \vec F$ in $\Omega$, we have 
\begin{multline*}
\|\A(\delta \nabla^2 u)\|_{L^{p}(\partial \Omega,\sigma)} \leq C_p \Big( \|\N^2(\nabla u)\|_{L^{p}(\partial \Omega,\sigma)}  + \|\A^2_1(\delta h)\|_{L^{p}(\partial \Omega,\sigma)}  \\
+  \|\A(\delta \nabla \vec F)\|_{L^{p}(\partial \Omega,\sigma)}  + \|\N^2(\vec F) \|_{L^{p}(\partial \Omega,\sigma)} \Big).\end{multline*}
The constant $C_p$ depends only on $p$, the ellipticity, boundedness and DKP constant of $L$, and the 2-sided CAD constants of $\Omega$. In any case, the finiteness of the quantities involved are not guaranteed. However, the left-hand side will be finite as long as the right-hand side is.
\end{theorem}

\bp
The proof of the two above theorems are similar. Let us do only Theorem \ref{ThS<N11}. We take a family of smooth compactly supported cut-off functions $\varphi_k \uparrow \1_\Omega$ that has a uniform bound on their cut-off constants, see Corollary \ref{corcutoff}. Theorem \ref{ThS<N1} gives
\[ \|\mathcal A(\varphi_k^2\phi^2\delta\nabla u)\|_{L^{2}(\partial \Omega,\omega_L)} \leq C \Big( \|\N^2(\varphi_k\phi u)\|_{L^{2}(\partial \Omega,\omega_L)} + \|\mathcal A_1^2(\varphi_k\phi\delta \vec F) \|_{L^{2}(\partial \Omega,\omega_L)} \Big),\]
with a constant $C$ dependent of the cut-off constant of $\phi$, but independent of $k$. The use of the smooth compactly supported function guarantees that
\[ \|\mathcal A(\delta\varphi_k^2\nabla u)\|_{L^{p}(\partial \Omega,\sigma)} +  \|\N^2(\varphi_k u)\|_{L^{p}(\partial \Omega,\sigma)} + \|\mathcal A_1^2(\varphi_k\delta \vec F) \|_{L^{p}(\partial \Omega,\sigma)} <+\infty \]
So Theorem \ref{Thgoodlambda} gives that 
\[ \|\mathcal A(\delta\varphi_k^2 \nabla u)\|_{L^{p}(\partial \Omega,\sigma)} \leq C \Big( \|\N^2(\varphi_ku)\|_{L^{p}(\partial \Omega,\sigma)} + \|\mathcal A_1^2(\delta \varphi_k \vec F) \|_{L^{p}(\partial \Omega,\sigma)} \Big).\]
Taking $k\to \infty$ and thus $\varphi_k \uparrow \1_\Omega$ gives the desired result.
\ep

\section{Proof of the main theorems}

\subsection{Proof of Theorem \ref{ThCarlesonPert}: The stability of the Poisson-Neumann problem under small Carleson perturbations of operators.}

We shall reuse the fixed point argument from Theorem \ref{ThfixedpointW}, and so the map $\mathscr N_{\vec F}: \, \widehat W^{1,2}(\Omega) \to \widehat W^{1,2}(\Omega)$ defined for $\vec F \in L^\infty_c(\Omega,\R^n)$ similarly to \eqref{defmapN} as
\begin{equation}
\mathscr{N}_{\vec F}(v)(X) := \int_\Omega \nabla_Y N_0(X,Y) \cdot [(A_0-A_1)(Y)\nabla v (Y)+ \vec F(Y)]\, dY,
\end{equation} 
where $N_0$ is the Neumann function for $L_0$. 

Provided that $\epsilon$ is small enough, Theorem \ref{ThfixedpointW} gives that for any $\vec F \in L^\infty_c(\Omega)$ and any initial data $v_0 \in \widehat W^{1,2}(\Omega)$ - for instance $v_0 \equiv 0$ - the sequence $(\mathscr{N}_{\vec F})^k(v_0)$ will converge in $\widehat W^{1,2}(\Omega)$ to the fixed point $u$ which is equal to
\begin{equation} \label{formulaforu}
u(Y) = \int_\Omega \nabla_Y N_1(X,Y) \cdot \vec F(Y)\, dY \quad \text{ for a.e. } Y \in \Omega.
\end{equation}

\medskip

For $p\in (1,\infty)$, we introduce the norms 
\[ \|v\|_{Y_p} := \|\N^2(v)\|_{L^p(\partial \Omega,\sigma)}, \  \|v\|_{Z_p} :=  \|\A(\delta \nabla v)\|_{L^p(\partial \Omega,\sigma)}, \ \text{ and } \ \|v\|_{X_p}:= \|v\|_{Y_p} + \|v\|_{Z_p} .\]
Let $X_p:=\{ v \in W^{1,2}_{loc}(\Omega), \, \|v\|_{X_p} < +\infty\}$, and similarly $Y_p :=\{ v \in L^{2}_{loc}(\Omega), \, \|v\|_{Y_p} < +\infty\}$ and $Z_p :=\{ v \in W^{1,2}_{loc}(\Omega), \, \|v\|_{Z_p} < +\infty\}$. Both $X_p$ and $Y_p$ are clearly Banach spaces. 

Since the Poisson-Neumann problem \eqref{PNp}$_L$ is solvable, we have
\[ \|\mathscr{N}_{\vec F}(v)\|_{Y_p} \lesssim \|\mathcal A_1^2(\delta (|\vec F + [A_0-A_1]\nabla v)\|_{L^p(\partial \Omega,\sigma)} \]
whenever $\|\mathcal A_1^2(\delta (|\vec F + [A_0-A_1]\nabla v)\|_{L^p(\partial \Omega,\sigma)} < +\infty$. But, with the Carleson inequality (Theorem \ref{ThCarleson2}) we have that
\[
\|\mathcal A_1^2(\delta (|\vec F + [A_0-A_1]\nabla v)\|_{L^p(\partial \Omega,\sigma)} \lesssim  \|\mathcal A_1^2(\delta \vec F)\|_{L^p(\partial \Omega,\sigma)}  + \|A_0-A_1\|_{CM_2^\infty} \underbrace{\|\A(\delta \nabla v)\|_{L^p(\partial \Omega,\sigma)}}_{= \|v\|_{Z_p} \leq \|v\|_{X_p}} ,
\]
We also know thanks to Theorem \ref{ThS<N11} - which applies since $\omega_L \in A_\infty(\sigma)$ - that 
\begin{multline*}
\|\mathscr{N}_{\vec F}(v)\|_{Z_p} \lesssim \|\mathscr{N}_{\vec F}(v)\|_{Y_p} + \|\mathcal A_1^2(\delta (|\vec F + [A_0-A_1]\nabla v)\|_{L^p(\partial \Omega,\sigma)} \\
\lesssim \|\mathcal A_1^2(\delta \vec F)\|_{L^p(\partial \Omega,\sigma)} + \|A_0-A_1\|_{CM_2^\infty} \|v\|_{X_p}
\end{multline*}
with the previous computations. So as long as $\vec F$ is - for instance - in $L^\infty_c(\Omega)$, the map $\mathscr N_{\vec F}$ sends functions from $X_p$ into itself. 

Let us show that $\mathscr N_{\vec F}$ is also a contraction map. We have with similar computations, using  \eqref{PNp}$_{L_0}$ and Theorem \ref{ThS<N11}, that
\[\|\mathscr{N}_{\vec F}(v-w)\|_{X_p} \lesssim \|A_0-A_1\|_{CM_2^\infty} \|v-w\|_{X_p}.\]
So for $ \|A_0-A_1\|_{CM_2^\infty}$ small enough, depending only on the allowed constants, $\|\mathscr{N}_{\vec F}(v-w)\|_{X_p} \leq \frac12 \|v-w\|_{X_p}$. By Picard's fixed point theorem, the sequence $(\mathscr{N}_{\vec F})^k(v_0)$ introduced earlier also converge to the fixed point $u$ in $X_p$. In particular, $u\in X_p$, and 
\[\|u\|_{X_p} \leq C \|\mathcal A_1^2(\delta \vec F)\|_{L^p(\partial \Omega,\sigma)} + \frac12 \|u\|_{X_p}, \]
meaning that, for any $\vec F \in L^\infty_c(\Omega)$, the function $u$ defined in \eqref{formulaforu} verifies
\[ \|u\|_{X_p} \lesssim \|\mathcal A_1^2(\delta \vec F)\|_{L^p(\partial \Omega,\sigma)}\]
with a constant independent of $\vec F$. It implies that \eqref{PNp}$_{L_1}$, so Theorem \ref{ThCarlesonPert} follows.

\subsection{Proof of Theorem \ref{MainTh}: The stability of the Poisson-Neumann problem for special small DKP operators.} As the proof of this theorem is more technical than the previous one, we divide it into steps. The proof is (very) loosely based on the proof of \cite[Theorem 1.12]{Fenreg}, that established the duality between the $L^p$ Dirichlet problem and the $L^{p'}$ Dirichlet-regularity problem for DKP operators in the half space.

\medskip

First note that $\omega\in A_\infty(\sigma)$ because of Theorem \ref{ThUR=>Ainfty} and the fact that $\Omega$ is a Chord Arc Domain.

\medskip

\noindent{\bf Step 1: The introduction of $\mathscr N_h$:} The setting for the proof of this theorem is the same as for the proof of Theorem \ref{ThCarlesonPert}. For $h \in L^2(\Omega)$ with zero average, we will use the map $\mathscr N_{h}: \, \widehat W^{1,2}(\Omega) \to \widehat W^{1,2}(\Omega)$ defined similarly to \eqref{defmapN} as
\begin{equation}
\mathscr{N}_h(v)(X) := \int_\Omega \Big[ \nabla_Y N_0(X,Y) \cdot (A^T-I)(Y)\nabla v (Y) + N_0(X,Y) h(Y) \Big] \, dY,
\end{equation} 
where $N_0$ is the Neumann function for the Laplacian. 

Provided that $\|A-I\|_{L^\infty(\Omega)} \leq \epsilon$ is small enough, Theorem \ref{ThfixedpointW} gives that for any $h \in L^2(\Omega)$ zero average, the map $\mathscr N_h$ is a contraction map on $\widehat W^{1,2}(\Omega)$, and that sequence the $(\mathscr{N}_{h})^k(0)$ will converge in $\widehat W^{1,2}(\Omega)$ to the fixed point $u^*_h$ which is equal to
\begin{equation} \label{formulaforub}
u^*_h(X) = \int_\Omega N_1^*(X,Y) h(Y)\, dY \quad \text{ for a.e. } Y \in \Omega,
\end{equation}
where $N_1^*$ is the Neumann function for $L^*=-\diver A^T \nabla$. 

For $p'$, the H\"older conjugate of $p$, we introduce the norms 
\[ \|v\|_{X_{p'}}:= \|\N^2(v)\|_{L^{p'}(\partial \Omega,\sigma)} + \|\A(\delta \nabla v)\|_{L^{p'}(\partial \Omega,\sigma)},\]
and 
\[ \|v\|_{\bar X_{p',1}}:= \|\nabla v\|_{X_{p'}} + \|v\|_{W^{1,2}(\Omega)}.\]
Associated to those norms are the Banach space $X_{p'}:= \{ v \in W^{1,2}_{loc}(\Omega), \, \|v\|_{X_{p'}} < +\infty\}$ and $\bar X_{p',1} := \{v \in \widehat W^{1,2}(\Omega) \cap W^{2,2}_{loc}(\omega),\, \|\bar v\|_{X_{p'},1} < +\infty \}$.

\medskip

Recall that $T^{p',1,2}_{-1}(\Omega)$ is the space of functions $f\in W^{1,2}_{loc}(\Omega)$ satisfying $\| \mathcal A_1^2(\delta )\|_{L^{p'}(\partial \Omega,\sigma)}$ is finite. We want to show, for $h\in T^{p',1,2}_{-1}(\Omega)  \cap L^2(\Omega)$ with zero average,
\begin{enumerate}[(a)]
\item $\nabla \mathscr N_h(0) \in X_{p'}$, and more precisely
\begin{equation} \label{claimFPDKPa} 
\|\nabla \mathscr N_h(0)\|_{X_{p'}} \lesssim  \Big\| \mathcal A_1^2(\delta h) \Big\|_{L^{p'}(\partial \Omega,\sigma)}
\end{equation}
with a constant independent of $h$.
\item The map $\mathscr N_h$ verifies the bound
\begin{equation} \label{claimFPDKPb} 
\big\| \nabla [\mathscr N_h(v) - \mathscr N_h(w)] \big\|_{X_{p'}} \lesssim C\epsilon \big\| \nabla [v-w]\big\|_{X_{p'}},
\end{equation}
with $C$ that is independent of $h$ and depend on $A$ only via its ellipticity and boundedness constants.
\end{enumerate}
We also know by Theorem \ref{ThfixedpointW} that
\begin{enumerate}[(a)] \addtocounter{enumi}{2}
\item  $\mathscr N_h$ maps $\widehat W^{1,2}(\Omega)$ to $\widehat W^{1,2}(\Omega)$, and that $\mathscr N_h$ is a contraction map on $\widehat W^{1,2}(\Omega)$ when $\epsilon$ is small enough.
\end{enumerate}

\medskip

The combination of (a), (b) and (c) shows that
\begin{itemize}
\item $\mathscr N_h$ maps $\bar X_{p',1}$ into itself.
\end{itemize}
And then if $\epsilon$ is small enough,
\begin{itemize}
\item $\mathscr N_h$ is a contraction map on $\bar X_{p',1}$.
\item The fixed point $u^*_h$ given in \eqref{formulaforub} belongs to $\bar X_{p',1}$.
\item The map $\mathscr N_h$ is being affine, the fixed point verifies  
\[ \|\nabla u^*_h\|_{X_{p'}} \lesssim  \|\nabla \mathscr N_h(0)\|_{X_{p'}} \lesssim \Big\| \mathcal A_1^2(\delta h) \Big\|_{L^{p'}(\partial \Omega,\sigma)} \]
\end{itemize}
We conclude that, for $h\in T^{p',1,2}_{-1}(\Omega)  \cap L^2(\Omega)$ with zero average,
\[\|\nabla u^*_h\|_{X_{p'}} \lesssim \Big\| \mathcal A_1^2(\delta h) \Big\|_{L^{p'}(\partial \Omega,\sigma)} \]
whenever $u^*_h$ is the weak solution given by \eqref{formulaforub}, i.e. the weak solution to $L^* u = h$ with zero Neumann data and zero average. It means that the Poisson-Neumann-Regularity problem is solvable for $L^*$ in $L^{p'}$, and then, by Theorem \ref{ThwPN}, that the Poisson-Neumann problem \eqref{PNp}$_L$ is solvable. The Theorem follows, provided that \eqref{claimFPDKPa} and  \eqref{claimFPDKPb} above are proven.

 \medskip

\noindent{\bf Step 2(a): Estimate on $\|\nabla \mathscr N_h(0)\|_{X_{p'}}$.} Let $h\in T^{p',1,2}_{-1}(\Omega)  \cap L^2(\Omega)$ and has zero average. Indeed, in this case, since the Poisson-Neumann problem \eqref{PNp}$_{-\Delta}$ is solvable, the Poisson-Neumann-Regularity is also solvable in $L^{p'}$ - see Theorem \ref{ThwPN} - and so $\mathscr N_{h}(0)$, which is the weak-solution to $-\Delta u= h$ in $\Omega$ with zero Neumann data and zero average, verifies
\[ \Big\|\N^2\big(\nabla \big[\mathscr N_{h}(0)\big]\big)\Big\|_{L^{p'}(\partial \Omega,\sigma)} \lesssim \Big\| \mathcal A_1^2(\delta h) \Big\|_{L^{p'}(\partial \Omega,\sigma)}.\]
But Theorem \ref{ThS<N21} - which applies since $\omega \in A_\infty(\sigma)$ - we also have that 
\[ \Big\|\A\big(\delta\nabla \big[\mathscr N_{h}(0)\big]\big)\Big\|_{L^{p'}(\partial \Omega,\sigma)} \lesssim  \Big\|\N^2\big(\nabla \big[\mathscr N_{h}(0)\big]\big)\Big\|_{L^{p'}(\partial \Omega,\sigma)} +  \Big\| \mathcal A_1^2(\delta h) \Big\|_{L^{p'}(\partial \Omega,\sigma)}.\]
Altogether, we deduce that $\nabla \mathscr N_{h}(0) \in X_{p'}$ and
\[ \| \nabla \mathscr N_{h}(0) \|_{X_{p'}} \lesssim \Big\| \mathcal A_1^2(\delta h) \Big\|_{L^{p'}(\partial \Omega,\sigma)} < +\infty,\]
as desired.

\medskip

\noindent{\bf Step 2(b): The map $\mathscr N_h$ is contractive for the semi-norm $\|\nabla . \|_{X_{p'}}$.} Take $h\in T^{p',1,2}_{-1}(\Omega)  \cap L^2(\Omega)$ with zero average, and take $v,w \in \bar X_{p',1}$. 

To lighten the notation, write $\bar v $ and $\bar w$ for $ \mathscr N_h(v)$ and $ \mathscr N_h(w)$ respectively. Then 
\[ \bar v - \bar w= \int_\Omega \nabla_Y N_0(X,Y) \cdot (A^T-I) \nabla (v-w) \, dX.\]
As such, $\bar v - \bar w$ a solution to $-\Delta (\bar v - \bar w) = - \diver [(A^T-I) \nabla (v-w)]$. Since $(A-I)\nabla G = 0$ on $\Omega \setminus 2B_\Omega$, we have $(A^T-I) \nabla (v-w) \cdot \nabla G = 0$ on $\Omega \setminus 2B_\Omega$, and so Theorem \ref{ThS<N21} gives that 
\begin{multline*} \label{SN<NN}
 \big\|\A\big(\delta\nabla[\bar v - \bar w] \big) \big\|_{L^{p'}(\partial \Omega,\sigma)} \lesssim  \big\|\N^2 \big(\nabla[\bar v - \bar w] \big) \big\|_{L^{p'}(\partial \Omega,\sigma)} \\ +  \big\|\A \big(\delta \nabla([A^T-I]\nabla(v - w)) \big) \big\|_{L^{p'}(\partial \Omega,\sigma)} + \big\|\N^2\big([A^T-I]\nabla(v - w)\big)\big\|_{L^{p'}(\partial \Omega,\sigma)} 
 \end{multline*}
 The latest two terms verify
\[\big\|\N^2\big([A^T-I]\nabla(v - w)\big)\big\|_{L^{p'}(\partial \Omega,\sigma)} \leq \|A-I\|_{L^\infty(\Omega)} \big\|\N^2\big(\nabla[v - w]\big)\big\|_{L^{p'}(\partial \Omega,\sigma)} \leq \epsilon \big\|\N^2\big(\nabla[v - w]\big)\big\|_{L^{p'}(\partial \Omega,\sigma)}\]
 and, with the Carleson inequality (Theorem \ref{ThCarleson0}),
\begin{multline}
 \big\|\A \big(\delta \nabla([A^T-I]\nabla(v - w)) \big) \big\|_{L^{p'}(\partial \Omega,\sigma)}  \\
 \lesssim \|A-I\|_{L^\infty(\Omega)} \big\|\A\big(\delta\nabla^2[v - w] \big) \big\|_{L^{p'}(\partial \Omega,\sigma)} + \|\delta \nabla A\|_{CM_2^\infty} \big\|\N^2\big(\nabla[v - w] \big) \big\|_{L^{p'}(\partial \Omega,\sigma)} \lesssim \epsilon \|\nabla (v-w)\|_{X_{p'}},
 \end{multline}
 meaning that we have
 \[ \big\|\nabla (\bar v - \bar w)\big\|_{X_{p'}} \lesssim \big\|\N^2 \big(\nabla[\bar v - \bar w] \big) \big\|_{L^{p'}(\partial \Omega,\sigma)} +  \epsilon \|\nabla (v-w)\|_{X_{p'}}.\]
So, to prove that $\mathscr N_h$ is a contraction maps when $\epsilon$ is small enough, or more precisely \eqref{claimFPDKPb}, it remains to prove that 
\begin{equation} \label{claimFPDKPc}
\big\|\N^2 \big(\nabla[\bar v - \bar w] \big) \big\|_{L^{p'}(\partial \Omega,\sigma)} \lesssim \epsilon \|\nabla (v-w)\|_{X^{p'}}.
\end{equation}
\medskip

The proof of \eqref{claimFPDKPc} is done by duality. We use a compact set $K\subset \Omega$ to guarantee the {\em a priori} finiteness of $\|\N^2(\1_K \nabla[\bar v - \bar w] )\|_{L^{p'}}$. Since $L^\infty_c(\Omega)$ is dense in $T^{p,1,2}_{-1}(\Omega)$, and by duality (\cite[Lemma A.8]{Mando26}, see also \cite[Proposition 2.4]{MPT}), 
there exists $\vec F \in L^\infty_c(\Omega,\R^n)$ such that 
\begin{equation} \label{CondonvecF}
\|\A_1^2(\delta \vec F)\|_{L^p} \leq 1
\end{equation}
and
\begin{equation} \label{fbonN}
 \big\|\N^2 \big(\1_K \nabla[\bar v - \bar w] \big) \big\|_{L^{p'}(\partial \Omega,\sigma)} \lesssim \int_\Omega \nabla[ \bar v - \bar w] \cdot \1_K \vec F\, dX.
  \end{equation}
Without loss of generality, we can assume that $\1_K \vec F = \vec F$. Let now $u \in \widehat W^{1,2}(\Omega)$ be the solution to $-\Delta u = -\diver \vec F$ with zero Neumann data, or equivalently, the only function in $\widehat W^{1,2}(\Omega)$ verifying
\begin{equation} \label{defuPN}
 \int_\Omega \nabla \phi \cdot \nabla u \, dX = \int_\Omega \nabla \phi \cdot \vec F\, dX \quad \text{ for all } \phi \in \widehat W^{1,2}(\Omega).
 \end{equation}
By Theorem \ref{ThS<N11}, \eqref{PNp}$_{-\Delta}$, and then \eqref{CondonvecF}, we get that 
\begin{equation} \label{Condonu}
\|\A(\delta \nabla u)\|_{L^p(\partial \Omega,\sigma)} \lesssim  \|\N^2(u)\|_{L^p(\partial \Omega,\sigma)} + \|\A^2_1(\delta \vec F)\|_{L^p(\partial \Omega,\sigma)} \lesssim \|\A^2_1(\delta \vec F)\|_{L^p(\partial \Omega,\sigma)} \leq 1
\end{equation}
Recall that $(\bar v - \bar w) \in \widehat W^{1,2}(\Omega)$, since $v,w\in \widehat W^{1,2}(\Omega)$ and $\mathscr N_h$ maps $\widehat W^{1,2}(\Omega)$ into itself, so we can use \eqref{defuPN} with $\phi = \bar v - \bar w$ in \eqref{fbonN} to obtain
\begin{multline} \label{fbonN2} \big\|\N^2 \big(\1_K \nabla[\bar v - \bar w] \big) \big\|_{L^{p'}(\partial \Omega,\sigma)} \lesssim \int_\Omega \nabla[ \bar v - \bar w] \cdot \1_K \vec F\, dX = \int_\Omega \nabla[\bar v - \bar w] \cdot \nabla u\, dX \\ =  \int_\Omega (A^T-I)\nabla[v - w] \cdot \nabla u\, dX,
\end{multline}
because $\bar v - \bar w$ is the weak solution in $\widehat W^{1,2}(\Omega)$ to $-\Delta (\bar v - \bar w) = - \diver (A^T-I)\nabla (v-w)$
 with zero Neumann data, meaning that 
\begin{equation} \label{defbarv-w}
 \int_\Omega \nabla[\bar v - \bar w] \cdot \nabla \phi \, dX = \int_\Omega (A^T-I)\nabla [v-w] \cdot \nabla \phi \, dX \quad \text{ for all } \phi \in \widehat W^{1,2}(\Omega).
 \end{equation}
 
 For safety measure, we add a cut-off function $\phi\in C^\infty_c(\Omega)$ to estimate the right-hand side of \eqref{fbonN2}. Since $(A-I) \nabla G = 0$, we have $(A^T-I)\nabla[v - w] \cdot \nabla G = 0$, and so we can apply Lemma \ref{IPPwithNG} again, and we have 
\begin{multline} \label{randombdB}
\int_\Omega (A^T-I)\nabla[v - w] \cdot \nabla u\, \phi\,  dX \lesssim \int_{\Omega}  |A-I||\nabla u| \, |\phi| |\delta \nabla^2[v-w]| \big]\, dX  \\
+  \int_{\Omega}  \left( |\delta \nabla A| + |A-I| \Big[ \frac{\delta^2 \big|\nabla^2 G\big|}{G} +\1_{\delta |\nabla G|<cG} + |\delta \nabla \phi|\Big] \right) |\nabla u| |\nabla [v-w]| |\phi| \, dX
\end{multline}
Thanks to Theorems \ref{ThonG} and \ref{CorNG>0}, we have 
\[ f:= \frac{\delta^2 \big|\nabla^2 G\big|}{G} +\1_{\delta |\nabla G|<cG} + |\delta \nabla \phi| \in CM_2^\infty(K),\]
where $K$ depends on the CAD constants of $\Omega$ and the cut-off constant of $\phi$. Together with the facts that $\delta \nabla A \in CM_2^\infty(\epsilon)$ and $|A-I| \leq \epsilon$, it yields that
\[|\delta \nabla A| + |A-I| \Big[ \frac{\delta^2 \big|\nabla^2 G\big|}{G} +\1_{\delta |\nabla G|<cG} + |\delta \nabla \phi|\Big] \in CM_2^\infty(C\epsilon)\]
and then, with Theorem \ref{ThCarleson0}, the bound \eqref{randombdB} becomes
\begin{multline*} 
\int_\Omega (A^T-I)\nabla[v - w] \cdot \nabla u\, \phi\,  dX \\
\lesssim \|A-I\|_{L^\infty(\Omega)} \big\|\A(\delta \nabla^2[v-w])\big\|_{L^{p'}(\partial \Omega,\sigma)} \big\|\A(\delta \nabla u)\big\|_{L^{p}(\partial \Omega,\sigma)} \\
+ \big\|\A(f\nabla[v-w])\big\|_{L^{p'}(\partial \Omega,\sigma)} \big\|\A(\delta\nabla u)\big\|_{L^{p}(\partial \Omega,\sigma)} \\
\lesssim \epsilon \left( \big\|\A(\delta \nabla^2[v-w])\big\|_{L^{p'}(\partial \Omega,\sigma)} +  \big\|\N^2(\nabla[v-w])\big\|_{L^{p'}(\partial \Omega,\sigma)} \right) \big\|\A(\delta\nabla u)\big\|_{L^{p}(\partial \Omega,\sigma)}
\\ = \epsilon \|\nabla[v-w]\|_{X_{p'}} \|\A(\delta\nabla u)\|_{L^p(\partial \Omega,\sigma)}.
\end{multline*}
Taking $\phi \uparrow \1_\Omega$ with uniformly bounded cut-off constants (possible thanks to Corollary \ref{corcutoff})  leads to 
\[
\int_\Omega (A^T-I)\nabla[v - w] \cdot \nabla u\, dX \lesssim \epsilon \|\nabla[v-w]\|_{X_{p'}} \|\A(\delta\nabla u)\|_{L^p(\partial \Omega,\sigma)} \lesssim  \epsilon \|\nabla[v-w]\|_{X_{p'}}
\]
by \eqref{Condonu}. We reinsert this last bound in \eqref{fbonN2}, and we obtain
\[\big\|\N^2 \big(\1_K \nabla[\bar v - \bar w] \big) \big\|_{L^{p'}(\partial \Omega,\sigma)} \lesssim  \epsilon \|\nabla[v-w]\|_{X_{p'}}.\]
Taking $K\uparrow \Omega$ proves the desired claim \ref{claimFPDKPc}. The Theorem follows.

\subsection{Proof of Theorem \ref{ThLipPert}: The Stability of the Poisson-Neumann problem under small bi-Lipschitz transformations of domains.} \label{SSMainThProof}

In an ideal word, we want the following variant of  Lemma 4.2 in \cite{FLZ??}.

\begin{theoremc} \label{ThLrhoc}
Let $\Omega$ be a bounded Lipschitz domain, and write $G$ for the Green function (for $-\Delta$ in $\om$) with pole at center. There exist two constants $\epsilon_0 \in (0,1]$ and $C_0>0$ - both depending only the Lipschitz character of $\Omega_0$ -  such that the following holds. 

For any $\epsilon \in [0,\epsilon_0]$, take any bi-Lipschitz map $\Phi:\R^{n}\to\R^n$  satisfying $\|\nabla \Phi - I\|_\infty \leq \epsilon$.
There exists another bi-Lipschitz transformation $\rho=\rho_{\Phi}$ that maps $\Omega$ to $\Phi(\Omega)$ and satisfies $\rho(x) = \Phi(x)$ when $x\in \partial \Omega$ such that, if $L_\rho:= - \diver A_\rho \nabla$ is the conjugate of $-\Delta$ by $\rho^{-1}$, there exists a symmetric operator $L_0 := -\diver A_0 \nabla$ on $\Omega_0$ such that
\begin{enumerate}[(i)]
\item $\|A_0 - I\|_\infty \leq C_0 \epsilon$;
\item $(A_0 - I)\nabla G = 0$ in $\om_0 \setminus 2B_{\Omega_0}$,
\item $\delta \nabla A_0 \in CM_2^\infty(C_0\epsilon)$,
\item $|A_\rho - A_0| \in CM_2^\infty(C_0\epsilon)$, 
\end{enumerate}
where the two Carleson conditions are with respect to domain $\om$.
\end{theoremc}

Provided that the above conjecture is true, the proof of Theorem \ref{ThLipPert} would succinctly be

 \begin{tikzcd}
(PN_{p})_{-\Delta,\Omega} \arrow[rr,Rightarrow, bend left, "{\tiny  \begin{array}{c}  (A_0-I)\nabla G = 0 \\ \|A_0-I\|_\infty \leq C_0\epsilon \\ \delta \nabla A_0 \in CM_2^\infty(C_0\epsilon) \end{array}}" above, "{\text{\tiny \begin{tabular}{c} Thm \ref{MainTh} \end{tabular}}}" below] 
 & &  (PN_{p})_{L_0,\Omega} \arrow[rr,Rightarrow, bend left, "{\tiny A_0 - A_\rho \in CM_2^\infty(C_0 \epsilon)}" above, "{\text{\tiny \begin{tabular}{l} Thm \ref{ThCarlesonPert} \end{tabular}}}" below] 
 & &  (PN_{p})_{L_\rho,\Omega} \arrow[rr,Leftrightarrow, bend left, "{\tiny \rho \text{ bi-Lipschitz}}" above] 
 & & (PN_{p})_{-\Delta,\Phi(\Omega)}
\end{tikzcd}

\medskip

Although the conjecture above has a good chance to be true, it would require to significantly modify the proof of Lemma 4.2 in \cite{FLZ??}, which treats the case where $\Omega$ is a Lipschitz graph domain and $G$ is the Green function with pole at infinity. Indeed, the lemma in \cite{FLZ??} rely on the following property: 
\begin{equation} \label{keyprop}
\delta(X)\nabla G(X) \cdot e_n \geq c \, G(X)
\end{equation}
for some vector $e_n$ and constant $c>0$ independent of $X\in \Omega$. 
On the other hand, as long as we have the above property, the construction will be the same, and a theorem can be reached. So we will use the following weaker/local version instead.

\begin{theorem} \label{ThLrho}
Let $\Omega$ be a bounded Lipschitz domain, and write $G$ for the Green function (for $-\Delta$ in $\om$) with pole at center. There exist three constants $r_0>0$, $\epsilon_0 \in (0,1]$ and $C_0>0$ - all three depending only on the Lipschitz character of $\Omega$ -  such that the following holds. 

For any $x\in \partial \Omega_0$, any $r\in (0,r_0)$, any $\epsilon \in [0,\epsilon_0]$,  and any bi-Lipschitz map $\Phi:\R^{n}\to\R^n$  satisfying:
\begin{enumerate}[(a)]
\item $\nabla \Phi = I$ in $\Omega \setminus  B(x,4r)$,
\item $\|\nabla \Phi - I\|_\infty \leq \epsilon$,
\end{enumerate}
there exists another bi-Lipschitz transformation $\rho=\rho_{\Phi}$ that maps $\Omega$ to $ \Phi(\Omega)$, satisfies $\rho(x) = \Phi(x)$ when $x\in \partial \Omega$, and such that, if $L_\rho:= - \diver A_\rho \nabla$ is the conjugate of $-\Delta$ by $\rho^{-1}$, there exists a symmetric operator $L_0 := -\diver A_0 \nabla$ on $\Omega_0$ such that
\begin{enumerate}[(i)]
\item $\|A_0 - I\|_\infty \leq C_0 \epsilon$; \smallskip
\item $(A_0-I)\nabla  G = 0$ in $\om \setminus 2B_{\Omega}$; \smallskip
\item $\delta \nabla A_0 \in CM_2^\infty(C_0\epsilon)$; \smallskip
\item $|A_\rho - A_0| \in CM_2^\infty(C_0\epsilon)$.
\end{enumerate}
\end{theorem}

\bp The key idea of the theorem is the fact that when $\Phi$ is just a translation outside of a small area $4B\cap \Omega$, and since $\rho_\Phi$ is meant to be a better version of $\Phi$, we can take $\rho_\Phi = \Phi$ outside of $32B \cap \Omega$, while constructing $\rho_\Phi$ in $32B\cap \Omega$ locally like in \cite{FLZ??}. The proof relies heavily on the computations in \cite{FLZ??}, which are long and highly non-trivial, and we preferred to avoid repeating the argument (as it would take over 15 pages if we had to rewrite all the details).

\bigskip

Since $\Omega$ is bounded and Lipschitz, we can represent $\Omega$ locally as the graph of a Lipschitz function, meaning that there exists $r_2$ and $M$ such that, for any $x\in \partial \Omega$, there exists a system of coordinates $(y,t)_x \in \R^{n-1} \times \R$ of $\R^n$ and a $M$-Lipschitz function $g_x :\, \R^{n-1} \to \R$ in this system of coordinates such that
\[ \Omega \cap B(x,r_2) = \{(y,t)_x \in \R^n, \, t>g_x(y) \} \cap B(x,r_2).\]
Moreover, $G$ is a positive harmonic function in each of the sets $\{(y,t)_x \in \R^n, \, t>g_x(y) \} \cap B(x,r_2)$, so by \cite[Theorem 11.10]{CS05}, there exists $r_1>0$ and $c_1>0$ independent of $x$ such that we have the key property 
\begin{equation} \label{keyprop2}
\delta \nabla G \cdot (0_{\R^{n-1}},1)_x \geq c_1 G \qquad \text{ in } B(x,r_1) \cap \Omega.
\end{equation}
\medskip

We set $r_0 = r_1/100$, and take $\epsilon_0$ to be the $\epsilon_2$ in \cite[Lemma 4.2]{FLZ??}. With this choices of $r_0$ and $\epsilon_0$, we can construct 
\[ \bar \rho_{\Phi}:\, B(x,20r_0) \cap \{(y,t)_x \in \R^n, \, t>g_x(y) \} \to \Phi(\Omega)\]
like we construct $\rho_\Omega$ in \cite{FLZ??}, because the construction depends only on the values of $G$ and $\Phi$ in $B(x,100r_0) \cap \Omega_0$ where the key property \eqref{keyprop2} applies. 

\medskip

Let $\eta \in C^\infty_0(\R,[0,1])$ such that $\eta \equiv 1$ on $[0,1]$, $\eta \equiv 0$ on $[2,\infty)$, and $|\eta'| \leq 2$. Define $\eta_{x,r}$ as
\begin{equation} \label{defetaxr}\eta_{x,r}(X) := \eta\Big( \frac{X-x}{r} \Big)\end{equation} 
Define now for $X\in \Omega$
\[ \rho_\Phi = \eta_{x,16r_0} \bar\rho_\Phi + (1- \eta_{x,16r_0}) \Phi,\]
and 
\[ A_0 := \eta_{x,8r_0} \bar A_0 + (1-\eta_{x,8r_0}) I,\]
where $\bar A_0$ be constructed from $\bar \rho_\Phi$ as $A_0$ is constructed from $\rho_\Phi$ in \cite[Lemma 4.2]{FLZ??}.

\medskip

Note that the construction of $\rho_\Phi(X)$ depends only on $\Phi$ via the values of $\Phi$ around $X$.  As such, the values of $\rho_\Phi$ in 
\begin{equation} \label{RegionTampon}
R_x:= [B(x,32r_0) \setminus B(x,8r_0)] \cap \{(y,t)_x \in \R^n, \, g_x(y) + r_0 >t>g_x(y) \}
\end{equation}
relies on $\Phi$ only via its values in $B(x,36r_0) \setminus B(x,4r_0)$, where it is a simple translation. So in this region, the key quantity $\lambda$ from \cite[Theorem 3.23]{FLZ??} used to build $\bar \rho_\phi$ is constant equal to $\Lambda = \Phi(X) - X$, meaning that $\mathbf B$ and $\mathbf C$ in \cite[Theorem 3.23]{FLZ??} are zero. So in the region $R_x$, the quantity $\bar \rho_\Phi(X)$ defined as in \cite[(3.39)]{FLZ??} is simply $X + \lambda$, i.e. $\Phi(X)$, and the matrix $\bar A_0$ is the identity. We conclude that
\begin{equation} \label{RxPr}
\rho_\Phi = \bar \rho_\Phi = \Phi \ \text{ and } \ A_0 = \bar A_0 = I \ \text{ in $R_x$.} 
\end{equation}

\medskip

Let us verify the desired properties $(i)-(iv)$.
\begin{enumerate}[$(i)$]
\item comes from the fact that $\|\bar A_0 - I\|_\infty \lesssim \epsilon$ by \cite[Lemma 4.2 (ii)]{FLZ??}; \smallskip
\item holds because $A_0 -I = (\bar A_0 - I)\eta_{x,8r_0}$, and $(\bar A_0 - I)\nabla G = 0$ as shown in the end of the proof of \cite[Lemma 4.2]{FLZ??}. \smallskip
\item Note first that in the relevant region $B(x,32r_0) \cap \Omega$, we have $\delta \nabla \bar A_0 \in CM_2^\infty(C\epsilon)$. This comes from the fact that $\bar A_0$ is - with the notation of \cite{FLZ??} - $\det(J_0)P(J_0)^{-T}(J_0)^{-1}P^{-T}$. But $J_0$ satisfies $\|J_0-I\|_\infty \lesssim \epsilon$ and  $\delta \nabla J_0 \in CM_2(C\epsilon)$ by \cite[Lemma 3.36]{FLZ??}, and $P$ is the orthogonal matrix that changes the cartesian basis to the orthonormal basis $(v_1,\dots,v_n)$, where $\delta \nabla v_i \in CM_2^\infty(C\epsilon)$ by \cite[Corollary 3.9]{FLZ??}.

By \eqref{RxPr}, observe now that $\nabla A_0$ is either $\nabla \bar A_0$ or $0$ except in the Whitney region 
\[ W_x:= B(x,32r_0) \cap \{ (y,t)_x, \, t\geq g_x(y) + r_0\}.\]
So
\[ |\delta \nabla A_0| \leq \underbrace{|\delta \nabla \bar A_0| \1_{B(x,32r_0)\cap \Omega_0}}_{\in CM_2^\infty(C\epsilon)} + \underbrace{|\bar A_0 - I|}_{\leq C\epsilon} \underbrace{|\nabla \eta_{x,8r_0}|}_{\leq 1/(4r_0)} \underbrace{\1_{W_x}}_{\in CM_2^\infty(C)} \in CM_2^\infty(C\epsilon)\]
as desired.
\smallskip

\item By the property \eqref{RxPr}, we have that $\rho_\Phi$ is $\bar \rho_\Phi$ or a translation except in $W_x$. It means that the matrix $A_\rho$, which is defined as
\[ A_\rho = \det(\Jac_{\rho_\Phi}) (\Jac_{\rho_\Phi})^{-T} (\Jac_{\rho_\Phi})^{-1}\]
with $\Jac_{\rho_\Phi}$ being the Jacobian matrix of $\rho_\Phi$, is either $\bar A_\rho$ or the Identity. So 
\[ |A_\rho - A_0| \leq |\bar A_\rho - \bar A_0| + C|\Jac_{\bar \rho_\Phi} - I| \1_{W_x}.\]
The difference $ |\bar A_\rho - \bar A_0| \in CM_2^\infty(C\epsilon)$ by \cite[Lemma 4.2 (iii)]{FLZ??}, $\1_{W_x} \in CM_2^\infty$ because $W_x$ is a Whitney region, and $\|\Jac_{\bar \rho_\Phi} - I\|_\infty \leq \epsilon$ because $\bar \rho_\Phi$ is a modification of $\Phi$, so is still a small bi-Lipschitz transformation (see \cite[Theorem 3.42]{FLZ??}).
\end{enumerate} 
The theorem follows.
\ep

\noindent{\em Proof of Theorem \ref{ThLipPert}.} Let $\Omega \subset \R^n$ and $p\in (1,\infty)$ be as in the theorem. Let $G$ be the Green function with pole at center. Let 
\begin{itemize}
\item $\epsilon_0$ be as in Theorem \ref{ThLrho},
\item $\epsilon_1$ be the $\epsilon$ constant (that depends on $p$, $n$, the constant $C_p(-\Delta,\Omega)$ in \eqref{PNp}$_{-\Delta,\Omega}$, and the CAD constants of $\Omega$) given by Theorem \ref{MainTh} when $\delta = C_p(-\Delta,\Omega)$,
\item $\epsilon_2$ be the $\epsilon$ constant in Theorem \ref{ThCarlesonPert} when $L_0$ is uniformly elliptic with constant 2, the constant in \eqref{PNp}$_{L_0,\Omega}$ is smaller than $2C_p(-\Delta,\Omega)$, and $\delta = 2C_p(-\Delta,\Omega)$. It means that ultimately, $\epsilon_2$ depends on $p$, $n$, the constant $C_p(-\Delta,\Omega)$ in \eqref{PNp}$_{-\Delta,\Omega}$ and the CAD constants of $\Omega$.
\end{itemize}
We take then 
\[ \epsilon:= \frac15 \min\{\epsilon_0,\epsilon_1/C_0,\epsilon_2/C_0, 1/100\},\]
where $C_0$ is the other constant Theorem \ref{ThLrho}, and then $\Phi$ is bi-Lipschitz transformation in $\R^n$ that satisfies $\|\Jac(\Phi) -I \|_\infty \leq \epsilon$.

\medskip

We aim to prove that the Poisson-Neumann problem \eqref{PNp}$_{-\Delta,\Phi(\Omega)}$ is solvable. One extra difficulty is the fact that  Theorem \ref{ThLrho} deals with {\em localized} bi-Lipschitz transformations and cannot be applied directly with $\Phi$, and so we want to use a localized version of the Poisson-Neumann problem. Fortunately for us, Corollary \ref{CorLocp} is there to help.

Let $r_0$ is the third constant from Theorem \ref{ThLrho}. By Corollary \ref{CorLocp}, since $\Phi(\Omega)$ is a Lipschitz domain, we can prove that \eqref{PNp}$_{-\Delta,\Phi(\Omega)}$ is solvable by proving that \eqref{Locp}$_{-\Delta,\Phi(\Omega)}$ holds. By Proposition \ref{prLocp}, it is sufficient to show that there is a $C>0$ such that for any $y\in \Phi(\partial \Omega)$, any $r\leq r_0$, and any solution to $u\in W^{1,2}(B(y,4r)\cap \Phi(\Omega)) $ to $-\Delta u = 0$ with zero Neumann data, we have
\begin{equation} \label{Locp1}
\|\N^2(\nabla u\1_{B(y,r)})\|_{L^{p'}(\partial \Omega)} \lesssim \left(\fint_{B(y,4r) \cap \Phi(\Omega)} |\nabla u|^2 \, dX \right)^\frac12.
\end{equation}

\medskip

Fix now $y\in \Phi(\partial \Omega)$. Define $x:= \Phi^{-1}(y) \in \partial \Omega$, and then $\Phi_x$ as
\[ \Phi_x(X) = \Phi(X) \eta_{x,2r_0}(X) + (X + y-x)(1-\eta_{x,2r_0}(X)) \qquad \text{ on } \R^n,\]
where $\eta_{x,2r_0}$ is as in \eqref{defetaxr}. We can check that
\[ \nabla \Phi_x - I = (\nabla \Phi - I)\eta_{x,2r_0} + ([\Phi(X) - X] - [\Phi(x) - x]) \nabla \eta_{x,2r_0},\]
so, since $\Phi - I$ is a Lipschitz map, 
\[ |\nabla \Phi_x - I| \leq \underbrace{|\nabla \Phi - I|}_{\leq \epsilon} + \underbrace{| \nabla \Phi - I|}_{\leq \epsilon} \underbrace{|X-x|}_{\leq 4r_0 \text{ on } \supp \eta_{x,2r_0}} \underbrace{|\nabla \eta_{x,2r_0}|}_{\leq 1/r_0} \leq 5\epsilon \leq \epsilon_0.\]
Moreover, it satisfies $\nabla \Phi = I$ in $\Omega \setminus B(x,4r_0)$ so Theorem \ref{ThLrho} applies and we can construct the map $\rho_{x}:= \rho_{\Phi_{x}}$, and the operators $L_{\rho_{x}} = -\diver A_\rho \nabla$ and $L_0 = -\diver A_0 \nabla$ as in the theorem.

\medskip

By Theorem \ref{ThLrho}, we have $\|A_0 - I\|_\infty \leq \epsilon_1$, $\delta \nabla A_0 \in CM_2^\infty(\epsilon_1)$ and $(A_0-I) \nabla G = 0$ in $\Omega \setminus 2B_\Omega$. Since moreover the Poisson-Neumann problem \eqref{PNp}$_{-\Delta,\Omega}$ is solvable with constant $C_p(-\Delta,\Omega)$,  Theorem \ref{MainTh} and our choice for $\epsilon_1$ entail that the Poisson-Neumann problem \eqref{PNp}$_{L_0,\Omega}$ is solvable with constant $2C_p(-\Delta,\Omega)$.

By Theorem \ref{ThLrho}, we also have that $A_\rho - A_0 \in CM_2^\infty(\epsilon_2)$. So Theorem \ref{ThCarlesonPert} and our choice for $\epsilon_2$ imply that the Poisson-Neumann problem \eqref{PNp}$_{L_{\rho_{x}},\Omega}$ is solvable with constant $4C_p(-\Delta,\Omega)$.

Finally, Theorem \ref{ThLocp} yields that \eqref{Locp}$_{L_{\rho_x},\Omega}$ holds with a constant independent of $x$ (since the constant in \eqref{PNp}$_{L_{\rho_x},\Omega}$ is independent of $x$ and $r$). In particular, for any $r\in (0,r_0)$ and any solution $v\in W^{1,2}(B(x,2r) \cap \Omega)$ to $L_{\rho_{x}} v = 0$ in $B(x,3r) \cap \Omega$ with zero Neumann data, we have
\begin{equation} \label{Locp2} \|\N^2(\nabla v \1_{B(x,2 r)})\|_{L^{p'}(\partial \Omega,\sigma)} \leq C \left( \fint_{B(x,3r) \cap \Omega} |\nabla v|^2 \, dX \right)^\frac12 \end{equation}
with a constant $C$ independent of $x \in \partial \Omega$.

\medskip

To conclude, take $r \in (0,r_0/2)$ and a weak solution $u\in W^{1,2}(B(y,4r)\cap \Phi(\Omega))$  to $-\Delta u = 0$ in $B(y,4r)\cap \Phi(\Omega)$ with zero Neumann data.  The function $v:= u\circ \rho_{x}$ is thus a weak solution to $L_{\rho_{x}} v = 0$ in $\rho_x^{-1}(B(y,4r)\cap \Phi(\Omega))$ with zero Neumann data. However, $\rho_x$ and $\Phi_x$ is a small bi-Lipschitz transformation (recall that $\epsilon \leq 1/500$), so
\[\rho_{x,r}(B(x,3r) \cap \Omega) \subset B(y,4r) \cap \Phi_{x}(B(x,4r)) \subset B(y,4r) \cap \Phi(B(x,4r))\]
because $\Phi_x = \Phi$ on $B(x,4r)$, and thus 
\[ \rho_{x,r}^{-1}(B(y,4r)\cap \Phi(\Omega)) \supset B(x,3r) \cap \rho_{x,r}^{-1}(\Phi(\Omega)) = B(x,3r) \cap \Omega.\]
It means that $v$ is a weak solution to $L_{\rho_{x}} v = 0$ in $B(x,3r) \cap \Omega$ with zero Neumann data, and consequently $v$ satisfies \eqref{Locp2}, i.e.
\[ \|\N^2(\nabla [u\circ \rho_{x}] \1_{B(x,2 r)})\|_{L^{p'}(\partial \Omega,\sigma)} \leq C \left( \fint_{B(x,3r) \cap \Omega} |\nabla [u\circ \rho_x]|^2 \, dX \right)^\frac12.\]
But $\rho_x$ is a (small) bi-Lipschitz change of variable, so for all $X\in \Omega$,
\[ |\nabla [u\circ \rho_{x}](X)| \approx |\nabla u(\rho_x(X))|,\]
for all $z\in \partial \Omega$ and $w\in L^2_{loc}(\Omega)$,
\[ \N^2(w)(\rho_x(z)) \leq \N^2_{(3,2)}(w \circ \rho_{x})(z),\]
$\rho_x^{-1}(B(y,r)) \subset B(x,2r)$, and $\rho_{x}(B(x,3r)) \subset B(y,4r)$. Together with the equivalence of the $\N^2_{\lambda,\kappa}$ under $L^p$ norm given in Lemma \ref{lemcompNA}, this shows that 
\begin{multline*}
 \|\N^2(\nabla u \1_{B(y,r)})\|_{L^{p'}(\Phi(\partial \Omega),\sigma)} \lesssim \|\N^2_{(3,2)}(\nabla [u\circ \rho_x] \1_{B(x,2r)})\|_{L^{p'}(\partial \Omega,\sigma)} \\ 
 \lesssim  \|\N^2(\nabla [u\circ \rho_x] \1_{B(x,2r)})\|_{L^{p'}(\partial \Omega,\sigma)} \lesssim \left( \fint_{B(x,3r) \cap \Omega} |\nabla [u\circ \rho_x]|^2 \, dX \right)^\frac12 \\ \lesssim  \left( \fint_{B(y,4r) \cap \Phi(\Omega)} |\nabla u|^2 \, dX \right)^\frac12.
 \end{multline*}
The constants in the above inequality being independent of $y\in \Phi(\partial \Omega)$, $r\in (0,r_0)$, and the weak solution $u\in W^{1,2}(B(y,4r)\cap \Phi(\Omega))$  to $-\Delta u = 0$ in $B(y,4r)\cap \Phi(\Omega)$ with zero Neumann data, it means that the claim \eqref{Locp1} holds. The theorem follows.
\ep

%

\bibliographystyle{alpha}
\bibliography{reference}

\end{document}